\documentclass[a4paper,11pt]{article}
\usepackage[utf8]{inputenc}
\usepackage[T1]{fontenc}
\usepackage{lmodern}
\usepackage{amsmath,amssymb,amsthm,mathtools}
\usepackage{enumitem}
\usepackage[hidelinks]{hyperref}
\newtheorem{theorem}{Theorem}[section]
\newtheorem{lemma}[theorem]{Lemma}
\newtheorem{proposition}[theorem]{Proposition}
\newtheorem{corollary}[theorem]{Corollary}

\newtheorem{fact}[theorem]{Fact}
\newtheorem{definition}[theorem]{Definition}
\newtheorem{question}[theorem]{Question}
\theoremstyle{remark}
\newtheorem{remark}[theorem]{Remark}

\newcommand{\HC}{\mathrm{HC}}
\newcommand{\Code}{\operatorname{Code}}
\newcommand{\Term}{\operatorname{Term}}
\newcommand{\cf}{\operatorname{cf}}
\newcommand{\cof}{\operatorname{cof}}
\newcommand{\supp}{\operatorname{supp}}
\newcommand{\dom}{\operatorname{dom}}
\newcommand{\Sig}{\Sigma^2}
\newcommand{\PiT}{\Pi^2}
\newcommand{\UP}{\mathrm{UP}}
\newcommand{\Club}{\mathbb{C}}
\newcommand{\R}{\mathbb R}
\newcommand{\Ord}{\operatorname{Ord}}

\hypersetup{
  pdftitle={Forcing global Pi2n uniformization with and without projective determinacy},
  pdfauthor={Stefan Hoffelner},
  pdfsubject={Set theory: forcing and second-order uniformization},
  pdfkeywords={uniformization, second-order pointclasses, forcing,
    stationary sets, almost disjoint coding, projective determinacy}
}

\title{Forcing global $\Pi^2_n$ uniformization,\\
with and without projective determinacy}
\author{Stefan Hoffelner\thanks{This research was funded in whole or in
part by the Austrian Science Fund (FWF)
[\href{https://doi.org/10.55776/P37228}{10.55776/P37228}].
For open access purposes, the author has applied a CC BY public
copyright license to any author-accepted manuscript version arising
from this submission.}}
\date{TU Wien\\12 November 2026}

\begin{document}
\maketitle

\begin{abstract}
We construct a set-forcing extension of $L$ which has the same reals
as $L$ and satisfies the $\Pi^2_n$-uniformization property for every
$n\geq1$. Assuming $\omega$ Woodin cardinals and a measurable cardinal
above them, we also obtain a forcing extension of
$L(\mathbb R)^{M_\omega}$ which has the same reals and satisfies
projective determinacy together with the $\Pi^2_n$-uniformization
property for every $n\geq1$. 
\end{abstract}

\noindent{\small\emph{2020 Mathematics Subject Classification.}
Primary 03E35; Secondary 03E15, 03E45.\par}
\smallskip
\noindent{\small\emph{Keywords.} Uniformization, second-order
pointclasses, forcing, stationary sets, almost disjoint coding.\par}

\section{Introduction}\label{sec:intro}

A \emph{uniformization} of a relation
$A\subseteq\omega^\omega\times\omega^\omega$ is a function
$U\subseteq A$ with $\dom U=\dom A$. A pointclass has the
\emph{uniformization property} if each of its relations has a
uniformization in the pointclass. We study the second-order
pointclasses $\Sig_n$ and $\PiT_n$ over the set $\HC$ of hereditarily
countable sets: first-order quantifiers range over $\HC$, and class
quantifiers range over all subsets of $\HC$. We write $\Sig_n$-$\UP$
and $\PiT_n$-$\UP$ for the uniformization properties of these
lightface pointclasses; the precise definitions are given in
Section~\ref{sec:prelim}.

For projective pointclasses, definable uniformizations come from two
sources: scales, which projective determinacy propagates through the
projective hierarchy~\cite{Moschovakis1971}, and good wellorders of
the reals, as in $L$. Both sources have second-order counterparts.
Under $\mathsf{AD}+V=L(\mathbb R)$, the scale analysis of
Steel~\cite{SteelScales} and Martin and Steel~\cite{MartinSteel}
gives the scale property, and hence uniformization, for $\Sig_1$.
Woodin extended this to
$\mathsf{ZF}+\mathsf{AD}^++V=L(\mathcal P(\mathbb R))$; see Steel and
Trang~\cite[Theorem~2]{SteelTrang}. Under the same hypothesis
$\mathsf{AD}+V=L(\mathbb R)$, an argument of Kechris and Solovay,
similar to the classical argument of L\'evy, gives a $\PiT_1$ relation
which has no uniformization in $L(\mathbb R)$
\cite[Theorem~3.11 and the following discussion]{SteelGames}.

In the presence of CH, a determinacy hypothesis for games of length
$\omega_1$ gives the uniformization property for $\PiT_1$. Building
on~\cite{SteelLong}, Steel~\cite{SteelOpen} showed that
$\omega_1$-open-projective determinacy gives the scale property for
the open game quantifier on analytical trees: the games have length
$\omega_1$ and integer moves, and the trees of countable positions are
definable over $(\HC,\in)$ from parameters. Under CH this pointclass
is $\PiT_1$, so the hypothesis yields $\PiT_1$-$\UP$
\cite[Corollary~0.6 and Theorem~2.1]{SteelOpen}. By Neeman's
work~\cite{Neeman}, in the form quoted by
Steel~\cite[Theorem~0.1]{SteelOpen}, the determinacy hypothesis holds
if every real belongs to a countable, $(\omega_1+1)$-iterable mouse
satisfying $\mathsf{ZFC}^-+$``there is a measurable Woodin cardinal''.
Steel notes that deriving this mouse-existence hypothesis from a
large-cardinal assumption requires iterability results which are not
yet known \cite[footnote~1]{SteelOpen}.

Good wellorders give the uniformization property for $\Sig_n$ at every
level $n$. In $L$, the canonical wellorder of the subsets of
$\omega_1$ is $\Sig_1$-good, and choosing least witnesses gives
$\Sig_n$-$\UP$ for every $n\geq1$ (Proposition~\ref{prop:Sigma-ground}).
We show that a forcing over $L$ which adds no reals can produce the
dual situation, $\PiT_n$-$\UP$ for every $n\geq1$. Such a forcing
preserves $\HC$ and all projective statements, and it changes the
second-order pointclasses only through the new subsets of $\HC$ which
it adds. An earlier forcing of this kind is the construction of
Abraham and Shelah~\cite{AbrahamShelah}: without adding reals, they
obtain a $\Delta^2_2$ wellorder of the reals and, by the same method,
code prescribed relations on the reals.

\begin{theorem}\label{thm:main}
There is a cardinal-preserving set-forcing extension $L[G]$ of $L$
which has the same reals as $L$ and satisfies $\PiT_n$-$\UP$ for
every $n\geq1$.
\end{theorem}

The construction is compatible with projective determinacy.

\begin{theorem}\label{thm:PD}
Assume there are $\omega$ Woodin cardinals and a measurable cardinal
above them, and put $\mathcal W=L(\mathbb R)^{M_\omega}$. There is a
forcing $\mathbb P\in\mathcal W$ which preserves all cardinals of
$\mathcal W$, adds no reals, and forces $\mathsf{PD}+\mathsf{CH}$ and
$\PiT_n$-$\UP$ for every finite $n\geq1$.
\end{theorem}

The model $\mathcal W$ is its own $L(\mathbb R)$ and satisfies
$\mathsf{ZFC}+\mathsf{CH}+\mathsf{PD}$ (Lemma~\ref{lem:PD-ground}).
Since $\mathbb P$ adds no reals, it changes neither the projective
sets nor the strategies for games on $\omega$, and hence preserves
$\mathsf{PD}$. Like $L$, the model $\mathcal W$ satisfies
$\Sig_n$-$\UP$ for every $n\geq1$: it has the form $L[A]$ for a set
$A\subseteq\omega_1$ whose membership and nonmembership relations have
parameter-free $\Sig_1$ definitions, and its canonical wellorder of
the subsets of $\omega_1$ is $\Sig_1$-good relative to $A$
(Proposition~\ref{prop:Sigma-ground}). Since boldface
$\Sig_n$-uniformization and boldface $\PiT_n$-uniformization cannot
hold at the same time (Lemma~\ref{lem:dual-UP}), we obtain the
following.

\begin{corollary}\label{cor:patterns}
Under the hypothesis of Theorem~\ref{thm:PD}, the models
$\mathcal W\subseteq\mathcal W[G]$ both satisfy
$\mathsf{ZFC}+\mathsf{PD}+\mathsf{CH}$, have the same $\HC$, and
satisfy, for every finite $n\geq1$,
\[
 \begin{aligned}
 \mathcal W&\models\Sigma^2_n\text{-}\UP\ \land\
                 \neg\boldsymbol\Pi^2_n\text{-}\UP,\\
 \mathcal W[G]&\models\PiT_n\text{-}\UP\ \land\
                 \neg\boldsymbol\Sigma^2_n\text{-}\UP.
 \end{aligned}
\]
The analogous statement holds for $L$ and the extension $L[G]$ of
Theorem~\ref{thm:main}.
\end{corollary}

Thus, relative to the hypothesis of Theorem~\ref{thm:PD}, each of the
two combinations in Corollary~\ref{cor:patterns} is consistent with
$\mathsf{PD}+\mathsf{CH}$, and a forcing which changes neither the
reals nor any projective statement turns the first into the second.

Large cardinals restrict the definable wellorders on which
$\Sig_n$-$\UP$ rests in our ground models. A theorem of Woodin rules
out a $\Sig_1$ wellorder of the reals under CH and a measurable Woodin
cardinal \cite[p.~2]{AbrahamShelah}. L\"ucke, Schindler and
Schlicht~\cite{LSS} obtain a good $\Sigma_1(\omega_1)$-definable
wellorder of $H_{\omega_2}$ in $M_1$ \cite[Theorem~5.2]{LSS}, and they
show that infinitely many Woodin cardinals with a measurable cardinal
above them preclude a $\Sigma_1(\omega_1)$-definable uniformization of
the club filter on $\omega_1$ \cite[Theorem~4.19]{LSS}. L\"ucke and
M\"uller~\cite[Theorem~1.6]{LuckeMueller} show that the same
hypothesis, or Woodin's axiom $(*)$, precludes a $\Sigma_1$-definable
wellorder of any subset of $\mathcal P(\omega_1)$ of size greater than
$\aleph_1$, with parameters in $\HC\cup\{\omega_1\}$.

The models of Theorems~\ref{thm:main} and~\ref{thm:PD} have no
measurable Woodin cardinal. Indeed, both satisfy CH and have a
$\Delta^2_1$ wellorder of the reals: in $L[G]$ the constructible
wellorder of the reals, and in $\mathcal W[G]$ the wellorder $\prec$
of Section~\ref{sec:PD}, which is $\Delta^2_1$ there by
Lemmas~\ref{lem:PD-ground} and~\ref{lem:PD-translation}. Steel's
theorem, on the other hand, gives $\PiT_1$-$\UP$ from CH and the
mouse-existence hypothesis above. We therefore ask the following.

\begin{question}\label{q:measurable-Woodin}
Is $\PiT_2$-$\UP$ consistent with $\mathsf{CH}$ and the existence of a
measurable Woodin cardinal? More generally, is the conjunction of
$\PiT_n$-$\UP$ for all $n\geq1$ consistent with these assumptions?
\end{question}

\paragraph{The method.}
Fix a recursive injection $\ell$ which codes a triple $(e,x,y)$,
consisting of a natural number and two reals, by a real $\ell(e,x,y)$.
Let $A_{n,e}(x,y)$ be the $\PiT_n$ formula with G\"odel number $e$.
In the final extension we define the uniformizing function by
\[
 U_{n,e}(x,y)\iff
 A_{n,e}(x,y)\land\neg\Code_n(\ell(e,x,y)).
\]
Each real $t$ in the range of $\ell$ indexes a member $S^n_t$ of a
canonical family of stationary subsets of $\omega_n$ in $L$. Given a
sequence $\vec T=(T_n)_{n\geq1}$ of sets of such reals, we construct a
forcing $F_{\vec T}$ in whose extension
\[
 \Code_n(t)\iff S^n_t\text{ is nonstationary in }\omega_n
             \iff t\in T_n.
\]
For $n=1$, $\Code_1(t)$ asserts that there is a club in $\omega_1$
disjoint from $S^1_t$. For $n\geq2$, a club in $\omega_n$ is too large
to be a class over $\HC$. Here $\Code_n(t)$ asserts that there is a
set $X\subseteq\omega_1$ from which iterated almost disjoint decoding
produces a set $C_X\subseteq\omega_n$\footnote{For
coding by stationary sets see Friedman and Magidor~\cite{FriedmanMagidor},
and for localization see Fischer and Friedman~\cite{FischerFriedman}
and Fischer, Friedman and Zdomskyy~\cite{FFZ}. A related application
is given in~\cite{HoffelnerSeparation}.}
with the following property: every suitable model $N$ which is
likewise coded by a subset of $\omega_1$, and which contains the
corresponding set decoded from $X$, satisfies that $C_X\cap\omega_n^N$
is a club in $\omega_n^N$ disjoint from $S^n_t\cap\omega_n^N$. The
ordinals $\omega_n^N$ are cofinal in $\omega_n$, and all these models
compute initial segments of the same set $C_X$; hence such a set $X$
yields a club in $\omega_n$ disjoint from $S^n_t$ in the final model.
Section~\ref{sec:decoders} expresses this property by a
parameter-free $\Sig_n$ formula. Thus $U_{n,e}$ is $\PiT_n$, and it
uniformizes $A_{n,e}$ provided that every nonempty section
$\{y:A_{n,e}(x,y)\}$ contains exactly one $y$ with
$\ell(e,x,y)\notin T_n$.

The difficulty is that the sections themselves depend on $\vec T$.
The forcing adds no reals, so every candidate $y$ is already
in $L$, but new subsets of $\HC$ can change the truth value of
$A_{n,e}(x,y)$.
To determine in $L$ which statements $A_{n,e}(x,y)$ a given forcing
$F_{\vec T}$ forces, we reduce $\PiT_n$ truth in the full extension to
$\PiT_1(H_{\omega_n})$ truth in the intermediate extension by the club
shooting forcing. Cone homogeneity of the club shooting forcing then
shows that each such statement with real parameters from the ground
model is decided by the trivial condition.

We start with $T_n=I$ for every $n$, where $I$ is the range of $\ell$,
and treat $n=1,2,\ldots$ in increasing order. At step $n$, take a pair
$(e,x)$ to which no value has been assigned yet and for which
$F_{\vec T}\Vdash A_{n,e}(x,y)$ for some $y$. Choose one such $y$,
assign it to $(e,x)$, and remove only $\ell(e,x,y)$ from $T_n$.
Shrinking $T_n$ while all other $T_m$ are fixed preserves every
$\PiT_n$ statement which is already forced, so $y$ remains in its
section. A section which was empty can become nonempty, so we repeat
this step, taking intersections of the sets $T_n$ at limit stages,
until every nonempty section has been assigned a value. Each step
assigns a value to a new pair $(e,x)$, so by CH the recursion stops
before $\omega_2$. Once $T_1,\ldots,T_n$ are fixed, changing the sets
$T_m$ for $m>n$, even all of them simultaneously, does not change which
$\PiT_n$ statements are forced. Thus each nonempty section in the final
model contains its assigned value $y$, and this $y$ is the only element
of the section with $\ell(e,x,y)\notin T_n$.

For Theorem~\ref{thm:PD} we write $\mathcal W=L[A]$, where
$A\subseteq\omega_1^{\mathcal W}$ codes the enumeration of the reals
in the order of construction of $M_\omega$, and relativize the
construction to $(\HC,\in,A)$. Membership and nonmembership in $A$
have parameter-free $\Sigma^2_1$ definitions which remain correct in
all relevant extensions; they allow us to eliminate $A$ from the
definitions of the uniformizing functions without increasing their
complexity.

\paragraph{Organization.}
Section~\ref{sec:prelim} collects preliminaries.
Section~\ref{sec:source} constructs the stationary sets and the club
shooting forcing. Subsection~\ref{subsec:level-two} sketches the case
$n=2$, which already shows one of the main ideas: we pass from larger
models to smaller ones by removing codes, and this process is carried
out entirely in the ground model.
Section~\ref{sec:master} defines the coding forcing and proves its
basic properties. Section~\ref{sec:names} shows that its bookkeeping
handles all relevant names and gives a first-order characterization of
the class of forcings used in the reduction of complexity.
Sections~\ref{sec:truthdrop} and~\ref{sec:decoders} give the reduction
of complexity and the formulas $\Code_n$.
Sections~\ref{sec:comparison} and~\ref{sec:selection} compare the
forcings $F_{\vec T}$ for different sequences $\vec T$ and complete the
proof of Theorem~\ref{thm:main}.
Section~\ref{sec:PD} relativizes the construction to $L[A]$, proves
Theorem~\ref{thm:PD} and derives Corollary~\ref{cor:patterns}.

\section{Preliminaries}\label{sec:prelim}

\subsection{Second-order pointclasses}

Let $H$ be a transitive set. We refer to subsets of $H$ as
\emph{classes over $H$}. A formula is $\Sig_n(H)$ if it has
the form
\[
 \exists X_1\forall X_2\cdots QX_n\,\varphi,
\]
where the class quantifiers range over $\mathcal P(H)$ and $\varphi$
is first order over $(H,\in)$, with the class variables as
additional predicates. Define $\PiT_n(H)$ dually. If no domain
is mentioned, it is $\HC$.

The pointclasses $\PiT_n$ and $\Sig_n$ are lightface: their
definitions use no parameters. Their boldface counterparts
$\boldsymbol\Pi^2_n$ and $\boldsymbol\Sigma^2_n$ allow real
parameters.

\begin{definition}
$\PiT_n$-$\UP$ means that every $\PiT_n$ relation
$A\subseteq\omega^\omega\times\omega^\omega$ has a $\PiT_n$
uniformization. Define $\Sigma^2_n$-$\UP$ analogously.
\end{definition}

Until Section~\ref{sec:PD}, the extensions constructed for the main theorem
add no reals and preserve $\omega_1$. Thus they have the same $\HC$ as $L$
and satisfy CH. Their collections of classes over $\HC$ may differ.

\subsection{Cardinals}
Until Section~\ref{sec:PD}, the ground model is $L$. We write
$\omega_n$ and $\aleph_\alpha$ for its cardinals, with $\omega_0=\omega$.
As proved below, the forcing extensions used in the construction preserve
all cardinals, so this notation has the same interpretation in each of them.

\subsection{A criterion for simultaneous uniformization}

Fix a recursive injection $\ell$ which maps every triple $(e,x,y)$,
where $e<\omega$ and $x,y$ are reals, to a real
$\ell(e,x,y)\in2^\omega$, and let
\[
 I=\{\ell(e,x,y):e<\omega,\ x,y\in\R^L\}
\]
be its range in $L$. We call the elements of $I$ \emph{codes} and
regard them as subsets of $\omega$. Thus $|I|=\omega_1$ in $L$.

For $\vec R,\vec T\in(\mathcal P(I))^\omega$, write
$\vec R\leq\vec T$ if $R_n\subseteq T_n$ for every $n\geq1$. The following criterion is proved
in Section~\ref{sec:selection}.

\begin{theorem}
\label{thm:abstract}
Suppose that $\vec T\mapsto F_{\vec T}$ is a function, definable in
$L$, which assigns a set forcing $F_{\vec T}\in L$ to every
$\vec T\in(\mathcal P(I))^\omega\cap L$, and that the following hold
for all such sequences $\vec R,\vec T$ and every $n\geq1$.
\begin{enumerate}[label=(A\arabic*)]
\item $F_{\vec T}$ preserves cardinals and adds no reals, and for
every second-order formula $\varphi$ and all reals $a$ in $L$, either
$1_{F_{\vec T}}\Vdash\varphi(a)$ or
$1_{F_{\vec T}}\Vdash\neg\varphi(a)$.
\item There is a parameter-free $\Sig_n$ formula $\Code_n$, not
depending on $\vec T$, such that
$F_{\vec T}\Vdash(\Code_n(t)\leftrightarrow t\in T_n)$ for every
$t\in I$.
\item If $\vec R\leq\vec T$ and $R_m=T_m$ for all $m\ne n$, then every
$\PiT_n$ statement with real parameters which is forced by
$F_{\vec T}$ is forced by $F_{\vec R}$.
\item If $R_m=T_m$ for all $m\leq n$, then $F_{\vec R}$ and
$F_{\vec T}$ force the same $\PiT_n$ statements with real parameters;
here $R_m$ and $T_m$ may differ for infinitely many $m>n$.
\end{enumerate}
Assume moreover that the relation ``$1_{F_{\vec T}}\Vdash\varphi(a)$''
is definable in $L$, as a relation between G\"odel numbers of
second-order formulas $\varphi$, sequences $\vec T$ and finite
sequences $a$ of reals.
Then there is $\vec T_*\in L$ such that $F_{\vec T_*}$ forces
$\PiT_n$-$\UP$ for every $n\geq1$.
\end{theorem}

\subsection{Forcing terminology}

We use the forcing and iteration terminology of \cite[Chapters 14--16]{Jech}.
Let $\rho$ be a regular cardinal and let $R$ be a forcing.
\begin{itemize}
\item $R$ is \emph{$<\rho$-closed} if every descending sequence of
length $<\rho$ has a lower bound, and \emph{$<\rho$-distributive} if it
adds no new sequences of ordinals of length $<\rho$.
\item A set of conditions is \emph{centered} if each of its finite
subsets has a common lower bound. $R$ is \emph{$<\rho$-compact} if
every centered set of fewer than $\rho$ conditions has a lower bound.
\item $R$ is \emph{$\rho$-linked} if it is the union of $\rho$ many
sets of pairwise compatible conditions. A function $c:R\to\rho$ such
that any two conditions with the same value are compatible is called a
\emph{linking function}.
\item $R$ is \emph{cone homogeneous} if for any two conditions $p,q$
there are $p'\leq p$ and $q'\leq q$ such that the cones
$\{r\in R:r\leq p'\}$ and $\{r\in R:r\leq q'\}$ are isomorphic.
\item Following \cite{CDMMS}, $R$ is
\emph{countably parallel-closed} if, whenever
$\langle p_n:n<\omega\rangle$ and $\langle q_n:n<\omega\rangle$
are descending sequences in $R$ such that $p_n$ and $q_n$ are
compatible for every $n<\omega$, there is $r\in R$ with
$r\leq p_n,q_n$ for all $n<\omega$.
\item For a regular uncountable cardinal $\rho$, $R$ satisfies the
\emph{stationary $\rho^+$-chain condition} if for every sequence
$\langle p_\alpha:\alpha<\rho^+\rangle$ in $R$ there are a club
$C\subseteq\rho^+$ and a regressive function
$f:C\cap\{\alpha<\rho^+:\cf(\alpha)=\rho\}\to\rho^+$ such that
$p_\alpha$ and $p_\beta$ are compatible whenever $\alpha,\beta\in C$
have cofinality $\rho$ and $f(\alpha)=f(\beta)$.
\end{itemize}
A \emph{projection} $\pi:P\to Q$ preserves the order and the largest
condition and, whenever $q\leq\pi(p)$, admits $p'\leq p$ with
$\pi(p')\leq q$. For the restriction maps considered below we shall
prove more: whenever $q\leq\pi(p)$, there is $p'\leq p$ with
$\pi(p')=q$ which agrees with $p$ outside the coordinates to which we
restrict.

An isomorphism between two cones of a forcing $R$ induces a
bijection between the corresponding names which maps $\check a$ to
$\check a$ for every set $a$ of the ground model. This gives the
following standard fact.

\begin{fact}\label{fact:cone}
Let $R$ be cone homogeneous, let $\varphi$ be a formula, and let $a$
be a set in the ground model. Then either
$1_R\Vdash\varphi(\check a)$ or $1_R\Vdash\neg\varphi(\check a)$.
\end{fact}
\begin{proof}
Otherwise there are conditions $p,q$ with $p\Vdash\varphi(\check a)$
and $q\Vdash\neg\varphi(\check a)$. Choose $p'\leq p$ and $q'\leq q$
with isomorphic cones. Then $p'\Vdash\varphi(\check a)$ implies
$q'\Vdash\varphi(\check a)$, a contradiction.
\end{proof}

The following fact combines two results. The chain condition is
the preservation theorem \cite[Theorem 1.2]{CDMMS}; the closure is
the standard fact that a $<\rho$-support iteration of
$<\rho$-closed forcings is $<\rho$-closed.

\begin{fact}\label{fact:CDMMS}
Let $\rho$ be regular uncountable with $\rho^{<\rho}=\rho$.
A $<\rho$-support iteration whose iterands are forced to be
$<\rho$-compact and to satisfy the stationary $\rho^+$-chain
condition is $<\rho$-closed and satisfies the stationary
$\rho^+$-chain condition, hence the $\rho^+$-chain condition.
\end{fact}

Here compactness gives the closure of the iterands and the countable
parallel closure which are required by
\cite[Theorem 1.2]{CDMMS}. Indeed, a $<\rho$-compact
forcing is $<\rho$-closed. Moreover, if
$\langle p_n:n<\omega\rangle$ and $\langle q_n:n<\omega\rangle$ are
descending sequences such that $p_n$ and $q_n$ are compatible for
every $n$, then $\{p_n:n<\omega\}\cup\{q_n:n<\omega\}$ is centered,
so countable compactness gives a common lower bound. Similarly,
$\rho$-linkedness implies the stationary $\rho^+$-chain condition: if
$c$ is a linking function, then $\alpha\mapsto c(p_\alpha)$ is
regressive on the club of all $\alpha\geq\rho$. Both implications are
recorded in \cite[Lemma 1.5]{CDMMS}. We therefore obtain the
following consequence, which is the form used below.

\begin{corollary}\label{cor:CDMMS}
Let $\rho$ be regular uncountable with $\rho^{<\rho}=\rho$. Every
$<\rho$-support iteration whose iterands are forced to be
$<\rho$-compact and $\rho$-linked is $<\rho$-closed and satisfies the
stationary $\rho^+$-chain condition; in particular it is
$\rho^+$-cc.
\end{corollary}

We shall verify the preservation of cardinals by means of the
following lemma.

\begin{lemma}\label{lem:cardinals-general}
Let $\mathbb P$ be a forcing. Suppose that for every $n\geq1$ there
are a $<\omega_n$-distributive forcing $Q_n$ and a $Q_n$-name
$\dot R_n$ for an $\omega_n$-cc forcing such that $\mathbb P$ is
forcing equivalent to $Q_n*\dot R_n$. Then $\mathbb P$ preserves all
cardinals $\leq\aleph_{\omega+1}$. If $\mathbb P$ is also
$\aleph_{\omega+2}$-cc, then it preserves all cardinals.
\end{lemma}
\begin{proof}
Let $G$ be $\mathbb P$-generic over $V$, fix $n\geq1$, and let
$M=V[G_{Q_n}]$ be the corresponding intermediate model. Since $Q_n$
adds no sequences of ordinals of length $<\omega_n$, all cardinals
$\leq\omega_n$ of $V$ are cardinals in $M$, $\omega_n$ is regular in
$M$, and
\[
 \cf^M(\aleph_{\omega+1})\geq\omega_n.
\]
The forcing $R_n=\dot R_n^{G_{Q_n}}$ is $\omega_n$-cc in $M$ and thus
preserves all cardinals $\geq\omega_n$ of $M$. Hence $\omega_n$ is a
cardinal in $V[G]$. As $n$ was arbitrary, all cardinals
$\leq\aleph_\omega$ are preserved.

Now let $f:\eta\to\aleph_{\omega+1}$ be a function in $V[G]$, where
$\eta<\omega_n$, and let $\dot f\in M$ be an $R_n$-name for $f$. For
$\zeta<\eta$ let $X_\zeta\in M$ be the set of all $\alpha$ such that
some condition of $R_n$ forces $\dot f(\zeta)=\check\alpha$. By the
$\omega_n$-chain condition, in $M$,
\[
 |X_\zeta|<\omega_n,\qquad
 X=\bigcup_{\zeta<\eta}X_\zeta\in[\aleph_{\omega+1}]^{<\omega_n},
 \qquad\sup X<\aleph_{\omega+1},
\]
and the range of $f$ is contained in $X$. Consequently, in $V[G]$,
\[
 \cf(\aleph_{\omega+1}^V)\geq\omega_n\qquad(n\geq1).
\]
Since $\aleph_\omega$ is singular in $V[G]$, the cofinality of
$\aleph_{\omega+1}^V$ in $V[G]$ is a regular cardinal greater than
$\aleph_\omega$. If $\aleph_{\omega+1}^V$ were not a cardinal in
$V[G]$, it would have cardinality $\aleph_\omega$ there, and hence
cofinality at most $\aleph_\omega$. Thus $\aleph_{\omega+1}^V$
remains a regular cardinal. Finally, the $\aleph_{\omega+2}$-chain
condition preserves all cardinals $\geq\aleph_{\omega+2}$.
\end{proof}

\subsection{Club shooting and almost disjoint coding}

For a regular $\rho$ and $S\subseteq\rho$ with $\rho\setminus S$ stationary, $\Club(\rho\setminus S)$ is the forcing of closed bounded subsets of $\rho\setminus S$, ordered by end extension.

Fix a regular $\rho$ with $\rho^{<\rho}=\rho$, let $\sigma=\rho^+$, and fix a constructible almost disjoint family $\langle B_\eta:\eta<\sigma\rangle$ of unbounded subsets of $\rho$, that is, a family with pairwise bounded intersections. For $D\subseteq\sigma$, the forcing $\mathbb{A}_\rho(D)$ consists of pairs $(s,F)$ with $s\in 2^{<\rho}$ and $F\in[D]^{<\rho}$, where $(s',F')\leq(s,F)$ iff $s'\supseteq s$, $F'\supseteq F$, and $s'(\alpha)=0$ whenever $\alpha\in\dom s'\setminus\dom s$ and $\alpha\in B_\eta$ for some $\eta\in F$. We call $s$ the \emph{stem} of $(s,F)$. If $G$ is generic and
$Z=\{\alpha<\rho:s(\alpha)=1\text{ for some }(s,F)\in G\}$, then
\[
 \eta\in D\iff Z\cap B_\eta\text{ is bounded in }\rho\qquad(\eta<\sigma).
\]
In this situation we say that $Z$ \emph{codes $D$ by almost disjoint
coding}. The forcing $\mathbb{A}_\rho(D)$ is $<\rho$-compact and $\rho$-linked, since conditions with the same stem are compatible.

\subsection{Absoluteness}

\begin{lemma}\label{lem:small}
Let $\rho$ be regular, let $R$ be $<\rho$-distributive, and let $P$
have a dense subset of size $<\rho$. Then $R$ is still
$<\rho$-distributive after forcing with $P$. That is, if
$G_P\times G_R$ is generic for the product $P\times R$, formed in
$V$, then $V[G_P][G_R]$ and $V[G_P]$ have the same sequences
of ordinals of length $<\rho$.
\end{lemma}
\begin{proof}
For $\rho=\omega$ the assertion is immediate. Otherwise, a fixed
dense suborder of $P$ of size $<\rho$ remains $\rho$-cc in $V^R$.
Apply \cite[Lemma~4.7]{StejskalovaEaston}.
\end{proof}

\begin{fact}[Silver]\label{fact:closedabs}
Let $\rho$ be regular with $\rho^{<\rho}=\rho$ and let $R$ be $<\rho$-closed. Then $\Sig_1(H_\rho)$ statements with parameters in $\mathcal P(H_\rho)\cap V$ are absolute between $V$ and $V^R$.
\end{fact}
\begin{proof}
Since $|H_\rho|=\rho$ and $H_\rho^{V^R}=H_\rho^V$,
code $(H_\rho,\in,\vec A)$ by a structure on $\rho$ in $V$.
Apply Silver's absoluteness lemma for existential second-order
formulas \cite[Lemma~3.5]{FuchsClosed}; see also
\cite[p.~298, (I6)]{Kunen}.
\end{proof}

\subsection{Term forcing}\label{sec:term}

Let $\rho$ be regular uncountable, $\sigma=\rho^+$, $\rho^{<\rho}=\rho$ and $\sigma^{<\rho}=\sigma$. Let $P$ be $\rho$-cc with a fixed enumeration of a dense subset $P_0$ of cardinality at most $\rho$, and let $\dot Q$ be a $P$-name with $P\Vdash\dot Q=\mathbb{A}_\rho(\dot B)$, where $\dot B$ is a $P$-name for a subset of $\sigma$.

\begin{definition}
The \emph{term forcing} $\Term(P,\dot Q)$ consists of all
$P$-names $\tau=(\dot s,\dot F)$ with $1_P\Vdash\tau\in\dot Q$,
where $\dot s$ and $\dot F$ are nice $P_0$-names for the two
components of a condition. Put $\tau\leq\upsilon$ iff
$1_P\Vdash\tau\leq\upsilon$.
\end{definition}

\begin{theorem}\label{thm:term}
$\Term(P,\dot Q)$ is $<\rho$-compact, $\rho$-linked and has a dense subset of size at most $\sigma$. The map $(p,\tau)\mapsto(p,\tau)$ is a projection
\[
 P\times\Term(P,\dot Q)\twoheadrightarrow P*\dot Q.
\]
\end{theorem}
\begin{proof}
The $\rho$-cc of $P$ and the regularity of $\rho$ give
\[
 1_P\Vdash\dot\delta<\rho
 \quad\Longrightarrow\quad
 \exists\beta<\rho\ (1_P\Vdash\dot\delta<\beta).
\]
Thus the lengths of the stems and the sizes of the second components
of elements of $\Term(P,\dot Q)$ are bounded below $\rho$. Using nice
names for stems and for enumerations of second components, choose sets
$\mathcal S,\mathcal N$ of representatives for stems and conditions
with
\[
 \begin{aligned}
 |\mathcal S|&\leq(\rho\cdot|P_0|)^{<\rho}=\rho,\\
 |\mathcal N|&\leq(\sigma\cdot|P_0|)^{<\rho}=\sigma.
 \end{aligned}
\]
Every element of $\Term(P,\dot Q)$ is equivalent to a member of
$\mathcal N$, so $\mathcal N$ is dense. For two elements with the same
name for the stem,
\[
 1_P\Vdash(\dot s,\dot F\cup\dot F')
       \leq(\dot s,\dot F),(\dot s,\dot F').
\]
Hence the map which sends $(\dot s,\dot F)$ to the index of the
representative of $\dot s$ in $\mathcal S$ is a linking function.

Let $\{\tau_\alpha:\alpha<\eta\}$ be centered in
$\Term(P,\dot Q)$, where $\eta<\rho$. Then $1_P$ forces that
$\{\tau_\alpha:\alpha<\eta\}$ is centered in $\dot Q$. Since $P$
preserves the regularity of $\rho$, the compactness of $\dot Q$ gives
\[
 1_P\Vdash\exists q\in\dot Q\;
              \forall\alpha<\eta\ (q\leq\tau_\alpha).
\]
By the maximum principle there is
$\tau\in\Term(P,\dot Q)$ with $\tau\leq\tau_\alpha$ for all
$\alpha<\eta$.

Put $\pi(p,\tau)=(p,\tau)$. Given
$(p',\dot q)\leq_{P*\dot Q}(p,\tau)$, choose a maximal antichain
$A\subseteq P$ containing $p'$ and a name $\tau'$ such that
\[
 p'\Vdash\tau'=\dot q,\qquad
 a\Vdash\tau'=\tau\quad(a\in A\setminus\{p'\}).
\]
Then $\tau'\in\Term(P,\dot Q)$, $1_P\Vdash\tau'\leq\tau$, and
\[
 (p',\tau')\leq_{P\times\Term(P,\dot Q)}(p,\tau),\qquad
 \pi(p',\tau')\leq_{P*\dot Q}(p',\dot q).
\]
The map $\pi$ preserves the order and the largest condition,
so it is a projection.
\end{proof}

\begin{remark}
In our applications, the term forcing is defined in a model in which
$\rho^{<\rho}=\rho$ holds, whereas the $P$-generic filter which
evaluates its elements is generic over a larger model. All
applications below use the conditions and the order of the term
forcing as defined in the smaller model.
\end{remark}

\section{Stationary sets and club shooting}\label{sec:source}

Fix $n\geq1$. In this section we define stationary subsets of
$\omega_n$ which are indexed by the reals in $I$, and the forcing
which adds clubs disjoint from some of them.

The construction is related to the stationary sets used in
\cite[Section 2]{HoffelnerReduction} and
\cite[Section 2, Steps 0--1]{FFZ}. We give the definitions and
proofs needed here for every finite $n$. For $n=1$, only the
diamond sequence and its stationary sets are needed.

\subsection{Constructible definitions}\label{sec:source-spec}

Fix the standard good $\Sigma_1$ wellorder $<_L$ of $L$, which
orders sets first by the level of the constructible hierarchy at which
they appear and then by their least definition;
see \cite[Chapter 13]{Jech}. Write $\mathrm{Init}(p,J)$ for
$p\in L$ and $J=\{q:q<_L p\}$, and fix a bounded formula
$q_{\rm init}$ with
\[
 \mathrm{Init}(p,J)\iff\exists w\,q_{\rm init}(p,J,w).
\]
For $p\in L_{\omega_n}$, both $J$ and a witness $w$ belong to $L_{\omega_n}$.

For $\delta<\omega_n$, we say that a strictly increasing cofinal map
$\omega_{n-1}\to\delta$ \emph{witnesses} $\cf(\delta)=\omega_{n-1}$.
We say that $\cf(\delta)\ne\omega_{n-1}$ is witnessed if $\delta=0$, or
by a cofinal map $\gamma\to\delta$ with $\gamma<\omega_{n-1}$.
Here a map is cofinal in $\delta$ if every $\eta<\delta$ is at most
some value of the map; thus a map with domain $1$ is a witness for a
successor ordinal.
These witnesses lie in $H_{\omega_n}^L=L_{\omega_n}$. Collecting them for
$\delta\leq a$ gives a $\Sigma_1$ definition of
\[
 A_a=(a+1)\cap\cof(\omega_{n-1})\qquad(a<\omega_n).
\]
We use these witnesses in the bounded formulas below.

\subsection{\texorpdfstring{A canonical $\Sigma_1$-definable diamond sequence}{A canonical Sigma1-definable diamond sequence}}

Put $\rho=\omega_{n-1}$ and $A=\omega_n\cap\cof(\rho)$.

\begin{lemma}\label{lem:canonical-diamond}
In $L$, there is a canonical $\diamondsuit(A)$-sequence
$\langle d_\delta:\delta\in A\rangle$ whose graph has a
functional $\Sigma_1$ definition $\Phi_D(\omega_{n-1},\delta,d)$,
uniformly in $n$. More precisely, there is a fixed bounded
formula $q_D$ such that
\[
 \Phi_D(\omega_{n-1},\delta,d)\iff\exists w\,q_D(\omega_{n-1},\delta,d,w),
\]
and, for every $\delta\in A$, a witness for
$\Phi_D(\omega_{n-1},\delta,d_\delta)$ belongs to $L_{\omega_n}$.
\end{lemma}
\begin{proof}
For $\gamma\in A$ let $(d_\gamma,c_\gamma)$ be the $<_L$-least
pair $(u,c)$ with $u\subseteq\gamma$, $c$ club in $\gamma$, and
\[
 \forall\xi\in c\cap A\quad u\cap\xi\ne d_\xi.
\]
Such a pair exists: the range of a continuous cofinal sequence in
$\gamma$ of length $\rho$, with successor ordinals at its initial and
successor steps, is a club in $\gamma$ which is disjoint from
$A\cap\gamma$. Jensen's condensation argument gives
$\diamondsuit(A)$ \cite[Lemma~6.5]{Jensen}. The clubs of order type
$\rho$ just mentioned do not affect this argument. Indeed, let
$(X,C)$ be the $<_L$-least counterexample, let $M$ be an elementary
submodel which contains it, with $\gamma=M\cap\omega_n\in A$, and let
$L_{\bar\beta}$ be the transitive collapse of $M$. Then $\gamma$ is
regular in $L_{\bar\beta}$, so no club in $\gamma$ of order type
$\rho$ belongs to $L_{\bar\beta}$. All $<_L$-predecessors of
$(X\cap\gamma,C\cap\gamma)$ belong to $L_{\bar\beta}$, and the
displayed condition is absolute for $L_{\bar\beta}$. Hence
$(d_\gamma,c_\gamma)=(X\cap\gamma,C\cap\gamma)$ as in the usual
argument, and $\gamma\in C$ gives a contradiction.

A witness for $\Phi_D(\rho,\delta,d)$ consists of the recursion
through $\delta$, the cofinality witnesses for $A_\delta$, and the
$<_L$-initial segments and the witnesses for leastness at each
stage. The displayed condition and its failure for each predecessor
are bounded. The recursion determines $d=d_\delta$ uniquely, also in
cardinal-preserving extensions: the cofinality partition below
$\omega_n$ and the constructible initial segments are absolute.
Each component has hereditary size at most $\rho$ in $L$, and there
are at most $\rho$ stages, so the witness belongs to $L_{\omega_n}$.
\end{proof}

For $t\subseteq\omega$, put
\[
 \mathrm{Tr}(d,t)\iff \omega\in d\ \wedge\ d\cap\omega=t,
\]
\[
 S_t=\{\delta\in A:\mathrm{Tr}(d_\delta,t)\},\qquad
 S_*=\{\delta\in A:d_\delta=\varnothing\}.
\]
The set $S_t$ is stationary, because the diamond sequence guesses the
set $t\cup\{\omega\}$ on a stationary set; similarly $S_*$ is
stationary, because the sequence guesses the empty set on a stationary
set. Distinct sets $S_t$ are disjoint, and all of them are disjoint
from $S_*$. Being a subset of $\omega$, the real $t$ is fixed by every
transitive collapse. When $n$ varies we write $S^n_t$ and $S^n_*$ for
these sets.

\subsection{\texorpdfstring{A $\Sigma_1$ formula for clubs disjoint from $S_t$}{A Sigma1 formula for clubs disjoint from St}}

For ordinal parameters $\rho,\nu$, let
$\psi_{\rm club}(\rho,\nu,C,t)$ assert the existence of a set
$K$ and functions $U,F$ satisfying the following conditions:
\begin{enumerate}
\item $C\subseteq\nu$ is club in $\nu$, and $K\subseteq\nu$.
\item $\dom U=\nu$. For $\delta\in K$, $U(\delta)$ is a
strictly increasing cofinal map from $\rho$ into $\delta$.
For $\delta\in\nu\setminus K$, either $\delta=0$ or
$U(\delta)$ is a cofinal map from some $\gamma<\rho$ into
$\delta$.
\item $\dom F=K$. Writing $F(\delta)=(a_\delta,w_\delta)$,
we require
\[
 \forall\delta\in K\quad
 q_D(\rho,\delta,a_\delta,w_\delta).
\]
\item $\forall\delta\in C\cap K\ \neg\mathrm{Tr}(a_\delta,t)$.
\end{enumerate}
The formula following $\exists K\,\exists U\,\exists F$ is
bounded, so $\psi_{\rm club}$ is $\Sigma_1$. For $\rho=\omega_{n-1}$
and $\nu=\omega_n$, the witnesses for the cofinalities give $K=A$, and
the functionality of $\Phi_D$ gives
$a_\delta=d_\delta$ for every $\delta\in K$. Consequently,
\[
 \psi_{\rm club}(\omega_{n-1},\omega_n,C,t)
 \iff C\text{ is club in }\omega_n\text{ and }C\cap S_t=\varnothing.
\]

\subsection{The constructible almost disjoint family}

Enumerate $\mathcal P(\omega_{n-1})^L$ in the order $<_L$ as
$\langle U_\eta:\eta<\omega_n\rangle$, identifying subsets with
their characteristic functions. Let $h(s)$ be the $<_L$-index
of $s$ among the constructible binary strings of length
$<\omega_{n-1}$; these strings are ordered by $<_L$ in order type
$\omega_{n-1}$. Put
\[
 B_\eta=\{h(U_\eta\restriction\zeta):\zeta<\omega_{n-1}\}.
\]
Since $<_L$ is a good $\Sigma_1$ wellorder and $L_{\omega_n}$
satisfies Collection, the function $\eta\mapsto B_\eta$ has a
functional $\Sigma_1$ definition $\Phi_B(\omega_{n-1},\eta,B)$ over
$L_{\omega_n}$.
Each $B_\eta$ is unbounded in $\omega_{n-1}$, and for
$\eta\ne\eta'$ the intersection $B_\eta\cap B_{\eta'}$ is bounded. We
always use this fixed family for almost disjoint coding at
$\omega_{n-1}$. When $n$ varies we write $B^{n-1}_\eta$ for $B_\eta$;
thus $\langle B^r_\eta:\eta<\omega_{r+1}\rangle$ is an almost
disjoint family of subsets of $\omega_r$.

\subsection{Coding clubs}\label{subsec:local-clubs}

Let $n\geq2$. A club in $\omega_n$ disjoint from $S_t$ will be
recognized by its initial segments in transitive models of the
following kind.

\begin{definition}\label{def:suitable}
Let $\mathsf{ZFC}^-$ denote ZFC without the Power Set axiom, with
Collection and Choice. A transitive model $N\models\mathsf{ZFC}^-$
is \emph{suitable for $n$} if $\omega_{n-1}^N=\omega_{n-1}$,
\[
 N\models \text{``$\omega_{n+1}$ exists and $\omega_{n+1}=\omega_{n+1}^{L}$''},
\]
and the following hold in $N$, with $\beta=\omega_n^N$.
\begin{enumerate}
\item over the set $L_\beta$, $\Phi_D$ is total and functional at every $\delta\in\beta\cap\cof(\omega_{n-1})$;
\item over $L_\beta$, $\Phi_B$ is total and functional at every $\eta<\beta$.
\end{enumerate}
\end{definition}

Since $N$ is transitive, its $L_\beta$ is the true $L_\beta$. As
$\Phi_D$ and $\Phi_B$ are $\Sigma_1$ and define functions in $L$,
upward absoluteness shows that a model $N$ which is suitable for $n$
computes $d_\delta$ and $B_\eta$ correctly for all
$\delta\in\beta\cap\cof(\omega_{n-1})$ and all $\eta<\beta$. In every
cardinal-preserving extension of $L$, the transitive collapse of an
elementary submodel of a sufficiently large $H_\Theta$ which has size
$\omega_{n-1}$ and contains $\omega_{n-1}$ as a subset is suitable for
$n$, by elementarity.

Put
\[
 \operatorname{Even}(P)=\{\alpha:2\cdot\alpha\in P\},\qquad
 \operatorname{Odd}(P)=\{\alpha:2\cdot\alpha+1\in P\}.
\]

Let $C\subseteq\omega_n\setminus S_t$ be club. In $L[C]$ choose
canonically a club $E\subseteq\operatorname{Lim}(C)$ with
\[
 \min E>\omega_{n-1},\qquad
 d(c):=\min(E\setminus(c+1))\geq c+\omega_{n-1}
       \quad(c\in E\cup\{0\}).
\]
For $\omega_{n-1}\leq d<\omega_n$ let $w_d\subseteq\omega_{n-1}$ be
the $<_L$-least code of a wellorder of $\omega_{n-1}$ of type $d$,
using the ordinal pairing obtained by ordering pairs first by
their maximum and then lexicographically.
Define the set $P=P^n_t(C)\subseteq\omega_n$ by
\[
 \begin{aligned}
 \operatorname{Even}(P)&=C,\\
 \operatorname{Odd}(P)&=
   \bigcup_{c\in E\cup\{0\}}\{c+\xi:\xi\in w_{d(c)}\}.
 \end{aligned}
\]
This defines the map $C\mapsto P^n_t(C)$ uniformly in $L[C]$. Thus
$P^n_t(C)$ codes $C$ on its even part, and its odd part codes the
successor function of the club $E$.

\begin{lemma}\label{lem:direct-padding}
Let $C\subseteq\omega_n\setminus S_t$ be club and put $P=P^n_t(C)$.
If $N$ is suitable for $n$, $\beta=\omega_n^N<\omega_n$, and
$t,P\cap\beta\in N$, then
\[
 N\models\psi_{\rm club}(\omega_{n-1},\beta,C\cap\beta,t).
\]
\end{lemma}
\begin{proof}
We first show that $\beta\in E$. Otherwise there are successive
elements $c<d$ of $E\cup\{0\}$ with $c<\beta<d$. Since $\beta$ is a
cardinal of $N$ above $\omega_{n-1}$, we have
\[
 c+\omega_{n-1}<\beta,\qquad
 w_d=\{\xi<\omega_{n-1}:c+\xi\in
                    \operatorname{Odd}(P\cap\beta)\}\in N.
\]
The ordinal pairing is absolute to $N$, which computes the order
type coded by $w_d$ correctly. In $N$ this order type has
cardinality $\omega_{n-1}$, giving
\[
 d<\omega_n^N=\beta<d.
\]

Thus $\beta\in E\subseteq\operatorname{Lim}(C)$, so
$C\cap\beta=\operatorname{Even}(P\cap\beta)$ is club in $\beta$
also in $N$. Put $\rho=\omega_{n-1}$. In $N$, take
$K=\{\delta<\beta:\cf(\delta)=\rho\}$ and choose the cofinality
witnesses $U$. For $\delta\in K$, a strictly increasing cofinal
map from the true $\rho$ into $\delta$ witnesses $\cf(\delta)=\rho$
in the universe and hence in $L$, since $\delta<\omega_n$.
Suitability gives the correct $d_\delta$ and a witness to
$\Phi_D(\rho,\delta,d_\delta)$ over $L_\beta$. By Collection we
obtain $F$ in $N$. Since $C\cap S_t=\varnothing$, the resulting
$K,U,F$ witness the assertion in $N$.
\end{proof}

We shall apply Lemma~\ref{lem:direct-padding} to models $N$ which
contain a set $Z\subseteq\omega_{n-1}$ coding $P$ by almost
disjoint coding. For such a model $N$, the correctness of
$\Phi_B$ in $N$ and the fact that $N$ has the true $\omega_{n-1}$ give
\[
 \{\eta<\beta:Z\cap B_\eta\text{ is bounded in }\omega_{n-1}\}
        =P\cap\beta\in N,
\]
so the lemma applies. For $t\in I$ we have $t\in L_\beta\subseteq N$.
For the converse direction we use the following observation.

\begin{lemma}\label{lem:club-from-segments}
Let $C\subseteq\omega_n$ and suppose that the set
\[
 D=\{\beta<\omega_n:C\cap\beta\text{ is a club in }\beta
                         \text{ disjoint from }S_t\}
\]
is unbounded in $\omega_n$. Then $C$ is a club in $\omega_n$
disjoint from $S_t$.
\end{lemma}
\begin{proof}
Since $D$ is unbounded and $C\cap\beta$ is unbounded in $\beta$ for
every $\beta\in D$, the set $C$ is unbounded in $\omega_n$. Every
element of $C$ belongs to $C\cap\beta$ for some $\beta\in D$, so
$C\cap S_t=\varnothing$. Finally, if $0<\delta<\omega_n$ and
$\sup(C\cap\delta)=\delta$, choose $\beta\in D$ with $\beta>\delta$.
Then $\delta$ is a limit point of $C\cap\beta$ below $\beta$, so
$\delta\in C\cap\beta\subseteq C$.
\end{proof}

In the application, $C=\operatorname{Even}(P)$ for the set $P$ coded
by some $Z\subseteq\omega_{n-1}$, and the elements of $D$ will be the
ordinals $\omega_n^N$ of transitive collapses $N$ of elementary
submodels of a sufficiently large $H_\Theta$ which have size
$\omega_{n-1}$ and contain $Z$ as an element and $\omega_{n-1}$ as a
subset. These ordinals are cofinal in $\omega_n$, the model $N$
contains the true $L_\beta$ for $\beta=\omega_n^N$, and the collapsing
map fixes $Z$ and every $B_\eta$ with $\eta<\beta$.

\subsection{The product of the club shooting forcings}\label{subsec:sourceproduct}

Recall that
\[
 I=\{\ell(e,x,y):e<\omega,\ x,y\in\R^L\}\subseteq2^\omega
\]
is the set of codes of triples consisting of a natural number and two
constructible reals.

For a sequence $\vec T=\langle T_n:n\geq1\rangle\in L$ with $T_n\subseteq I$, let $C_n(T_n)$ be the product of the forcings $\Club(\omega_n\setminus S^n_t)$, $t\in T_n$, with countable support for $n=1$ and full support for $n\geq2$. Let
\[
 C_{\vec T}=\prod_{n\geq1}C_n(T_n)\quad\text{(full support)},\qquad B_{\vec T}=L[G_{C_{\vec T}}].
\]
All these products are formed in $L$; in particular the factors
$C_n(T_n)$ are not redefined in intermediate extensions. For $\vec R\leq\vec T$, $C_{\vec R}$ is a complete subforcing of $C_{\vec T}$ via the natural projection of the product.

\begin{lemma}\label{lem:row-source}
For every $n\geq1$ and $T\subseteq I$, the forcing $C_n(T)$ has
cardinality at most $\omega_n$, is $<\omega_n$-distributive, and
preserves the stationarity of $S^n_*$ and of every $S^n_t$ with
$t\notin T$. For $n\geq2$ it is also
$<\omega_{n-1}$-closed.
\end{lemma}
\begin{proof}
Put $\rho=\omega_{n-1}$. The size bound follows from GCH and the
support bound $\rho$. For $n\geq2$, a descending sequence of length
$<\rho$ has a lower bound by taking unions and adding their suprema
at each coordinate, since $S^n_t\subseteq\operatorname{cof}(\rho)$.

Choose $M\prec H_\Theta$ containing the parameters and the initial
condition, with $\rho\subseteq M$, $|M|=\rho$, $M^{<\rho}\subseteq M$,
and $\delta=M\cap\omega_n$ in the stationary set to be preserved.
For distributivity use $S^n_*$. Such models are obtained at limits
of cofinality $\rho$ of a continuous elementary chain whose successor
models are closed under $<\rho$-sequences. Fix a cofinal sequence
$\langle\alpha_\eta:\eta<\rho\rangle$ in $\delta$. Construct a descending
sequence in $M$. At step $\eta$ meet the next open dense set in $M$,
or force a point of the named club in $(\alpha_\eta,\delta)$ for the
stationarity argument. Then end-extend the condition, at every
coordinate of its support, to a common maximum
$\gamma_\eta\in M\cap S^n_*$ above $\alpha_\eta$ and above all
current maxima. At limits below $\rho$ take lower bounds in $M$.

At the end, take the union at each coordinate in the union of the
supports and add $\delta$.
This is a lower bound because $\delta$ is outside every $S^n_t$ for
$t\in T$, and the maxima at each such coordinate are cofinal in
$\delta$. The union of the supports is the support of a condition:
it is countable if $n=1$, and the product has full support if
$n\geq2$. The lower bound meets all the given dense sets, or
forces $\delta$ into the named club, respectively.
\end{proof}

For $n\geq1$ write
\[
 Q_{\geq n}=\prod_{m\geq n}C_m(T_m),\qquad
 Q_{>n}=\prod_{m>n}C_m(T_m),\qquad
 Q_{<n}=\prod_{m<n}C_m(T_m).
\]

\begin{lemma}\label{lem:source}\label{lem:club-arithmetic-verified}
For every $n\geq1$, the forcing $Q_{\geq n}$ is
$<\omega_n$-distributive. The forcing
$C_{\vec T}$ preserves all cardinals and the GCH and adds no countable sequences of ordinals. The set $S^n_*$ and every $S^n_t$ with $t\notin T_n$ remain stationary. $C_{\vec T}$ is cone homogeneous.
\end{lemma}
\begin{proof}
Fix $n$. By Lemma~\ref{lem:row-source}, each factor $C_m(T_m)$
with $m>n$ is $<\omega_{m-1}$-closed, so $Q_{>n}$ is
$<\omega_n$-closed. Hence $Q_{>n}$ adds no elements of $H_{\omega_n}$
and preserves the stationarity of $S^n_*$ and of all sets $S^n_t$.
In particular $C_n(T_n)$ has the same conditions in the extension by
$Q_{>n}$, and the proof of Lemma~\ref{lem:row-source} applies in this
extension. Thus
\[
 Q_{\geq n}\text{ is }<\omega_n\text{-distributive},
 \qquad |Q_{<n}|<\omega_n,
\]
and $S^n_*$ and every $S^n_t$ with $t\notin T_n$ are stationary in the
extension by $Q_{\geq n}$. Since $Q_{<n}$ has size $<\omega_n$, it
preserves the stationarity of these sets.

The case $n=1$ shows that $C_{\vec T}=Q_{\geq1}$ adds no countable
sequences of ordinals. For the preservation of cardinals we apply
Lemma~\ref{lem:cardinals-general} with $Q_n=Q_{\geq n}$ and
$R_n=Q_{<n}$; note that $Q_{<n}$ is $\omega_n$-cc because of its size,
and that, in $L$,
\[
 |C_{\vec T}|\leq\aleph_\omega^\omega\leq2^{\aleph_\omega}=\aleph_{\omega+1},
\]
so $C_{\vec T}$ is $\aleph_{\omega+2}$-cc.

Next we verify the GCH in the extension. For $i\geq1$, the forcing
$Q_{>i}=Q_{\geq i+1}$ adds no subsets of $\omega_i$, and
$|Q_{\leq i}|\leq\omega_i$, where $Q_{\leq i}=Q_{<i+1}$. Counting nice
$Q_{\leq i}$-names for subsets of $\omega_i$ in $L$ therefore gives
\[
 2^{\omega_i}\leq\bigl((2^{\omega_i})^{\omega_i}\bigr)^L
                 =\omega_{i+1}.
\]
CH holds because no reals are added.
In the extension fix bijections
$e_i:\omega_{i+1}\to\mathcal P(\omega_i)$ ($i\geq1$).
The map
\[
 X\longmapsto
 \langle e_i^{-1}(X\cap\omega_i):i\geq1\rangle
 \qquad(X\subseteq\aleph_\omega)
\]
is injective and takes values in $(\prod_{i\geq1}\omega_{i+1})^L$,
since no countable sequences of ordinals are added.
Thus $2^{\aleph_\omega}\leq(\aleph_\omega^\omega)^L\leq\aleph_{\omega+1}$,
hence $2^{\aleph_\omega}=\aleph_{\omega+1}$.
Finally, for every cardinal $\tau\geq\aleph_{\omega+1}$,
nice names coded by subsets of $\tau\times C_{\vec T}$ give
$2^\tau\leq(2^{\tau\cdot\aleph_{\omega+1}})^L=(\tau^+)^L$.

For cone homogeneity, let $p,q\in C_{\vec T}$. Extend them to
$p'\leq p$ and $q'\leq q$ which have the same support and such that,
for every $n$, there is an ordinal $\gamma_n\in S^n_*$ which is the
maximum of $p'(n,t)$ and of $q'(n,t)$ for every $t$ in the support at
$n$. This is possible because $S^n_*$ is unbounded in $\omega_n$ and
disjoint from every $S^n_t$. The map which, for a condition below
$p'$, replaces the initial segment $p'(n,t)$ by $q'(n,t)$ at every
coordinate $(n,t)$ of the support of $p'$ is an isomorphism between
the cones below $p'$ and $q'$.
\end{proof}

For later reference we summarize the properties of the stationary
sets.

\begin{proposition}\label{hyp:source}
In $L$, for every $n\geq1$ there are, uniformly definably in $n$:
\begin{enumerate}[label=(S\arabic*)]
\item\label{S:partition} pairwise disjoint stationary sets $S^n_t\subseteq\omega_n\cap\operatorname{cof}(\omega_{n-1})$, $t\in I$, and a stationary set $S^n_*\subseteq\omega_n\cap\operatorname{cof}(\omega_{n-1})$ disjoint from all of them, defined from the canonical constructible diamond sequence on $\omega_n\cap\operatorname{cof}(\omega_{n-1})$;
\item\label{S:graph} for $n=1$, the function $(t,\alpha)\mapsto S^1_t\cap\alpha$ ($\alpha<\omega_1$) is first-order definable over $\HC$ in every outer model with the same reals, by the definition in Lemma~\ref{lem:canonical-diamond}, whose witnesses belong to $\HC$;
\item\label{S:padding} for $n\geq2$, a definable map
$C\mapsto P^n_t(C)\subseteq\omega_n$, computed in $L[C]$ from a club
$C$ disjoint from $S^n_t$, such that
$\operatorname{Even}(P^n_t(C))=C$ and such that the following holds: if
$N$ is suitable for $n$,
$\beta=\omega_n^N<\omega_n$, and
$P^n_t(C)\cap\beta\in N$, then
\[
 N\models\psi_{\rm club}\bigl(\omega_{n-1},\beta,
          \operatorname{Even}(P^n_t(C)\cap\beta),t\bigr);
\]
\item\label{S:stationary} in the extension by the forcing $C_{\vec T}$ of Subsection~\ref{subsec:sourceproduct}, every $S^n_t$ with $t\notin T_n$ remains stationary.
\end{enumerate}
\end{proposition}
These properties follow from the definitions above,
Lemma~\ref{lem:direct-padding}, and Lemmas~\ref{lem:row-source} and~\ref{lem:source}.

\subsection{The case \texorpdfstring{$n=2$}{n=2}}
\label{subsec:level-two}

We illustrate the main idea of our construction by sketching how to
obtain $\Pi^2_2$ uniformization. A coding forcing enlarges the
universe, which fits well with $\Sigma$ statements, whereas $\Pi$
statements behave well when the universe shrinks.

Recall that $\vec T=\langle T_n:n\geq1\rangle$ determines the
forcing $C_{\vec T}$: for $t\in T_n$, the factor
$\Club(\omega_n\setminus S^n_t)$ of $C_n(T_n)$ adds a club disjoint
from $S^n_t$. In the recursion for $n=2$ we vary only $T_2$ and keep
every $T_m$ with $m\ne2$ fixed.

We begin with $T_2=I$, so that the first model is the largest one. We
then define the uniformizing function by successively removing factors
$\Club(\omega_2\setminus S^2_t)$ from $C_2(T_2)$, that is, by
shrinking $T_2$; in this way we pass to smaller and smaller models.
Values which were chosen at earlier stages remain in their sections
throughout this process.

Put $F_{\vec T}=C_{\vec T}*\dot D_{\vec T}$, where $D_{\vec T}$
is the coding forcing constructed below. This forcing depends on the
entire sequence $\vec T$. Proposition~\ref{prop:exact} gives a single
$\Sigma^2_2$ formula $\Code_2$, independent of $\vec T$, such that
for every such sequence
\[
 F_{\vec T}\Vdash
       \bigl(\Code_2(t)\longleftrightarrow t\in T_2\bigr)
       \qquad(t\in I).
\]
The aim is to arrange that every nonempty section contains exactly one
$y$ with $\ell(e,x,y)\notin T_2$; this $y$ will be the value of the
uniformizing function.

As long as there is a pair $(e,x)$ to which no value has been assigned
and whose section is forced to be nonempty by the current forcing
$F_{\vec T}$, choose the $<_L$-least such pair and the $<_L$-least
$y$ with $F_{\vec T}\Vdash A_{2,e}(x,y)$, assign $y$ to $(e,x)$, and
replace
\[
 \vec T\longmapsto
 \bigl\langle T_1,\ T_2\setminus\{\ell(e,x,y)\},\ T_3,\ldots\bigr\rangle.
\]
The forcing adds no reals, so every candidate $y$ belongs to $L$.
The reduction of complexity below and cone homogeneity allow all these
choices to be made in $L$.

The chosen value remains in its section because the relevant
statements are downward absolute. Indeed, suppose that
$R_2\subseteq T_2$ and $R_m=T_m$ for $m\ne2$, and that the generic
filter for $C_{\vec R}$ is the restriction of the generic filter for
$C_{\vec T}$. Then
\[
 B_{\vec R}\subseteq B_{\vec T},\qquad
 H_{\omega_2}^{B_{\vec R}}=H_{\omega_2}^{B_{\vec T}}.
\]
Proposition~\ref{prop:compile} expresses $\Pi^2_2$ truth in the full
extension as $\Pi^2_1$ truth over $H_{\omega_2}$ in $B_{\vec T}$.
Passing from $B_{\vec T}$ to $B_{\vec R}$ leaves $H_{\omega_2}$
unchanged and removes some of its subsets. A $\Pi^2_1(H_{\omega_2})$
statement therefore remains true. Consequently,
\[
 F_{\vec T}\Vdash A_{2,e}(x,y)
 \quad\Longrightarrow\quad
 F_{\vec R}\Vdash A_{2,e}(x,y),
\]
as proved in Proposition~\ref{prop:ownaxis}. Every value which was
chosen earlier is therefore preserved when further codes are removed.

A section which was empty can become nonempty when we pass to a
smaller model. We therefore repeat the choice until every nonempty
section has been assigned a value, taking intersections of the sets
$T_2$ at limit stages. There are only $\omega_1$ pairs $(e,x)$, and
each step assigns a value to a new one, so the recursion stops before
$\omega_2$. For the resulting sequence $\vec T$,
\[
 U_{2,e}(x,y)\quad\Longleftrightarrow\quad
 A_{2,e}(x,y)\wedge\neg\Code_2(\ell(e,x,y))
\]
is a $\Pi^2_2$ uniformization of $A_{2,e}$: each nonempty section
contains its chosen value $y$, and this is the only element of the
section with $\ell(e,x,y)\notin T_2$.
Subsequent changes of the sets $T_m$ for $m>2$ preserve this
uniformization, by Proposition~\ref{prop:higheraxis} and the
equivalence between $\Code_2(t)$ and $t\in T_2$.
Section~\ref{sec:selection} carries out this recursion for every
finite $n$.

\section{The coding forcing and its quotients}\label{sec:master}

Fix $\vec T$ and recall that
\[
 B_{\vec T}=L[G_{C_{\vec T}}],
\]
where $G_{C_{\vec T}}$ is $C_{\vec T}$-generic over $L$ and
$C_{\vec T}=\prod_{n\geq1}C_n(T_n)$ is the product of club shooting
forcings defined in Subsection~\ref{subsec:sourceproduct}.
Write $B=B_{\vec T}$. In this section the ground model is $B$. By
Lemma~\ref{lem:source}, $B$ has the same cardinals as $L$ and
satisfies the GCH.

The forcing $D=D_{\vec T}$ is an iteration of length
$\aleph_{\omega+2}$ whose coordinates are assigned levels $i\geq1$.
The iterand at a coordinate of level $i$ is a $<\omega_i$-compact,
$\omega_i$-linked forcing on $\omega_{i+1}$, and it is defined in the
extension by only a part of the preceding coordinates. There are two
kinds of iterands. First, at level $i$ we force with the iterations
of such forcings which are named over suitable subforcings (the
$\Gamma_i$-iterations defined below), and with their restrictions
below conditions. Secondly, we force with term forcings for almost
disjoint coding; these are used to code sets which belong to
extensions by coordinates of lower levels. A bookkeeping function
enumerates the names for all these forcings. We first describe the
general form of such an iteration, in which each coordinate depends
only on some of the earlier coordinates, together with its conditions
and order. We then define $D$ by recursion along the bookkeeping
function, and finally prove its basic properties.

\subsection{\texorpdfstring{$\Gamma_i$-iterations}{Gamma-i-iterations}}

\begin{definition}\label{def:small-iterations}
Fix $i\geq1$. A \emph{$\Gamma_i$-iteration} is a
$<\omega_i$-support iteration
\[
 \langle P_\alpha,\dot Q_\alpha:\alpha<\delta\rangle,
 \qquad |\delta|\leq\omega_{i+1},
\]
such that for every $\alpha<\delta$ there is a $P_\alpha$-name
$\dot{\leq}_\alpha$ for a binary relation on $\omega_{i+1}$ with
\[
 P_\alpha\Vdash
 ``\dot Q_\alpha=(\check\omega_{i+1},\dot{\leq}_\alpha)
 \text{ is a }<\omega_i\text{-compact and }\omega_i\text{-linked preorder}.''
\]
We take $0$ as the largest condition of each iterand.

Let $\Gamma_i$ consist of the last forcings $P_\delta$ of these
iterations, given with their iteration codes. The codes are described
in Subsection~\ref{subsec:eligibility-justification}.
\end{definition}

By Corollary~\ref{cor:CDMMS}, in every model satisfying
$\omega_i^{<\omega_i}=\omega_i$, every $\Gamma_i$-iteration is
$<\omega_i$-closed and satisfies the $\omega_{i+1}$-chain condition.
If $q\in P_\delta$, the cone $P_\delta\restriction q$ is forcing
equivalent to a $\Gamma_i$-iteration. At stage $\alpha$, restrict
the iterand below $q(\alpha)$ and re-index it with largest condition
$0$. By induction this gives the cone below $q\restriction\alpha$
at each stage. At limits, the union of the support of a condition
with that of $q$ still has size $<\omega_i$.

\subsection{Iterations with partial memory}\label{subsec:memory}

The forcing $D$ is an iteration in which the iterand at a coordinate
$\xi$ is defined from the generic filter on a set
$m(\xi)\subseteq\xi$ of earlier coordinates only, the \emph{memory}
of $\xi$. Such iterations are often called iterations with partial
memory.
Our use of memory is inspired by the
admissible support pairs in~\cite[Subsection~2.6]{HoffelnerPiUniformization}.
In addition we assign levels to the coordinates, which determine the
sizes of the supports and the degree of closure. We first describe the
general form of these iterations and their conditions. The bookkeeping
which chooses the iterands follows in
Subsection~\ref{subsec:forcing-definition}.

\begin{definition}\label{def:memory-scheme}
Let $\gamma$ be an ordinal. An \emph{iteration with partial memory of
length $\gamma$} assigns to every $\xi<\gamma$
\begin{itemize}
\item a \emph{level} $i(\xi)\geq1$,
\item a \emph{memory} $m(\xi)\subseteq\xi$, and
\item an \emph{iterand}: a $D_{m(\xi)}$-name $\dot Q_\xi$ for a
$<\omega_{i(\xi)}$-compact, $\omega_{i(\xi)}$-linked preorder
$(\omega_{i(\xi)+1},\dot\leq_\xi)$ with largest element $0$,
together with a name for a linking function.
\end{itemize}
Writing $I_k=\{\xi<\gamma:i(\xi)=k\}$, we require for every $\xi<\gamma$:
\begin{enumerate}[label=(\alph*)]
\item\label{mem:closed} $\eta\in m(\xi)$ implies $m(\eta)\subseteq m(\xi)$;
\item\label{mem:up} $m(\xi)\subseteq\bigcup_{k\geq i(\xi)}I_k$;
\item\label{mem:small} $|m(\xi)\cap I_k|\leq\omega_{k+1}$ for every $k$.
\end{enumerate}
A set $M\subseteq\gamma$ is \emph{closed} if $m(\xi)\subseteq M$ for
every $\xi\in M$.
\end{definition}

Thus the memory of a coordinate contains only coordinates of its own
level or of higher levels, and only few of them at each level. Every
initial segment of $\gamma$ is closed, since $m(\xi)\subseteq\xi$.
By~\ref{mem:closed} every memory is closed, and by~\ref{mem:up} so is
$\bigcup_{k\geq j}I_k$ for every $j$. The forcings $D_M$, for closed
$M$, are given by the next definition; we also write
$D\restriction M$ for $D_M$. Since $m(\xi)\subseteq\xi$, the name
$\dot Q_\xi$ refers only to a forcing which has already been defined.
As for ordinary iterations, Definitions~\ref{def:memory-scheme}
and~\ref{def:coding-conditions} form a single simultaneous recursion
on~$\gamma$.

In the next definition, the value $p(\xi)$ of a condition is taken
from a set $\mathcal N_\xi$ of names which is small enough for the
size and chain condition estimates below, and which is closed under
canonical lower bounds, so that the forcing is sufficiently closed.

\begin{definition}\label{def:coding-conditions}
Fix an iteration with partial memory of length $\gamma$, and let
$\xi<\gamma$ have level $i$. For $\eta<\omega_i$ and a sequence
$\vec\tau=\langle\tau_\zeta:\zeta<\eta\rangle$ of $D_{m(\xi)}$-names
forced by $1_{D_{m(\xi)}}$ to be ordinals below $\omega_{i+1}$,
let $b_{\xi,\eta}(\vec\tau)$ be a fixed $D_{m(\xi)}$-name for
\[
 \min\{b<\omega_{i+1}:\forall\zeta<\eta\ (b\mathrel{\dot\leq_\xi}\tau_\zeta)\}
\]
if $\vec\tau$ is $\dot\leq_\xi$-descending, and for $0$ otherwise.
(Whether $\vec\tau$ is descending is decided in the
$D_{m(\xi)}$-generic extension, after the names $\tau_\zeta$ and
$\dot\leq_\xi$ have been evaluated.)
Let $\mathcal N_\xi$ be the least set of $D_{m(\xi)}$-names which
contains $\check a$ for every $a<\omega_{i+1}$ and is closed under the
operations $b_{\xi,\eta}$, $\eta<\omega_i$.

For closed $M\subseteq\gamma$, the conditions of $D_M$ are the
functions $p$ with
\[
 \begin{gathered}
 \supp(p)\subseteq M,\qquad
 p(\xi)\in\mathcal N_\xi\quad(\xi\in\supp(p)),\\
 |\supp(p)\cap I_k|<\omega_k\quad(k\geq1).
 \end{gathered}
\]
Setting $p(\xi)=\check0$ for $\xi\notin\supp(p)$, put
\[
 q\leq p\quad\Longleftrightarrow\quad
 \forall\xi\in M\quad
 q\restriction m(\xi)\Vdash_{D_{m(\xi)}}
   q(\xi)\mathrel{\dot\leq_\xi}p(\xi).
\]
\end{definition}

\begin{remark}\label{rem:conditions}
The minimum in the definition of $b_{\xi,\eta}$ exists because a
descending sequence is centered and $\dot Q_\xi$ is
$<\omega_i$-compact. Since $\omega_i$ is regular, the closure defining
$\mathcal N_\xi$ is reached after $\omega_i$ steps. Equivalently,
$\mathcal N_\xi$ consists of the names given by well-founded terms of
size $<\omega_i$ built from the constants $\check a$ and the
operations $b_{\xi,\eta}$.
Two features of the definition are used throughout.
\begin{enumerate}
\item Closure under the operations $b_{\xi,\eta}$, together with the
bound on the supports at each level, gives coordinatewise lower bounds
of short descending sequences in the ground model
(Lemma~\ref{lem:native-structural}).
\item Under GCH,
$|\mathcal N_\xi|\leq\omega_{i+1}^{<\omega_i}=\omega_{i+1}$,
which gives the bound on the size of the forcing and the chain
condition (Lemma~\ref{lem:size-cc-verified}).
\end{enumerate}
All values are ordinals below $\omega_{i+1}$, so
$1_{D_{m(\xi)}}\Vdash\tau\in\dot Q_\xi$ for every
$\tau\in\mathcal N_\xi$. For closed $M$, the order of $D_M$ refers only
to the forcings $D_{m(\xi)}$ with $\xi\in M$. Hence $D_M$ is the
suborder of $D_\gamma$ consisting of the conditions whose support is
contained in $M$.
\end{remark}

Formally, $D_M$ is a preorder; we pass to its separative quotient
without further mention.

\subsection{The bookkeeping}\label{subsec:forcing-definition}

The forcing $D$ is $D_{\aleph_{\omega+2}}$ for an iteration with
partial memory of length $\aleph_{\omega+2}$ which is defined by
recursion below. Its coordinates come in blocks.
Each block consists either of the iterands of a $\Gamma_i$-iteration,
for some $i$, or of a single term forcing.

We code conditions and names by sets of ordinals. A forcing of size
at most $\omega_{i+1}$ is identified with a preorder on
$\omega_{i+1}$ whose largest condition is $0$; several ordinals may
represent the same condition. A nice name is then coded by the
ordinals representing the conditions and the values which occur in
it. The details of this coding are given in
Subsection~\ref{subsec:eligibility-justification}.

Let $\mathcal B$ be the set of all sets of ordinals below
$\aleph_{\omega+2}$ of size at most $\aleph_{\omega+1}$. By GCH,
\[
 |\mathcal B|\leq\aleph_{\omega+2}^{\aleph_{\omega+1}}
                  =\aleph_{\omega+2}.
\]
Fix a function $F:\aleph_{\omega+2}\to\mathcal B$ such that
\[
 \sup F^{-1}\{b\}=\aleph_{\omega+2}\qquad(b\in\mathcal B).
\]
The elements of $\mathcal B$ are regarded as candidates for codes of
the quadruples $(i,M,p,\dot X)$ considered in
Definition~\ref{def:coding-recursion}; whether $F(\beta)$ actually
codes such a quadruple is checked only at stage $\beta$ of the
recursion.

For $p\in D_M$ and a name $\dot R$ which $p$ forces to be a forcing,
let $\dot R^{\,p}$ be a name for the forcing which is equal to
$\dot R$ if $p$ belongs to the generic filter for $D_M$, and trivial
otherwise. We use this notation for names over $D_M$, over
$D_N$ for closed sets $N\supseteq M$, and over two-step iterations
$D_M*\dot P$. In the last two cases $p$ is identified with its
canonical image in the larger forcing, and the generic filter for
$D_M$ is the one induced by the generic filter for the larger
forcing.

\begin{definition}\label{def:coding-recursion}\label{def:requests}
We define an iteration with partial memory of length
$\aleph_{\omega+2}$ by recursion, appending one block of coordinates
at each stage $\beta<\aleph_{\omega+2}$. Let $\gamma_\beta$ be the
length of the iteration defined before stage $\beta$; thus
$\gamma_0=0$ and $\gamma_\lambda=\sup_{\beta<\lambda}\gamma_\beta$ for
limit ordinals $\lambda$.

At stage $\beta$, suppose that $F(\beta)$ codes a quadruple
$(i,M,p,\dot X)$ with the following properties.
\begin{enumerate}[label=(\roman*)]
\item\label{rec:M} $i\geq1$, and $M\subseteq\gamma_\beta$ is a closed
set such that $M\subseteq\bigcup_{k\geq i}I_k$ and
$|M\cap I_k|\leq\omega_{k+1}$ for all $k\geq i$.
\item\label{rec:name} $p\in D_M$, and $\dot X$ is a $D_M$-name.
\item\label{rec:cases} One of the following two cases applies.
\begin{description}
\item[Case 1 \textnormal{(iterations).}] $p$ forces that
$\dot X$ is a code for a $\Gamma_i$-iteration
$\langle P_\alpha,\dot R_\alpha:\alpha<\delta\rangle$ of fixed length
$\delta$ with $|\delta|\leq\omega_{i+1}$.
\item[Case 2 \textnormal{(term forcings).}] $p$ forces that
$\dot X=(P,\dot E)$, where $P$ is $\omega_i$-cc with a given dense
subset of size at most $\omega_i$, $P\Vdash\dot E\subseteq\omega_{i+1}$,
and $\omega_i^{<\omega_i}=\omega_i$ and
$\omega_{i+1}^{<\omega_i}=\omega_{i+1}$. In this case put $\delta=1$
and let $\dot R_0$ be a name for the term forcing
$\Term(P,\mathbb{A}_{\omega_i}(\dot E))$, defined in
the extension by $D_M$ and identified with a preorder on
$\omega_{i+1}$.
\end{description}
\end{enumerate}
If moreover $\delta>0$, append the coordinates
$\xi_\alpha=\gamma_\beta+\alpha$, for $\alpha<\delta$, with
\[
 i(\xi_\alpha)=i,\qquad
 m(\xi_\alpha)=M\cup\{\xi_\zeta:\zeta<\alpha\},\qquad
 \dot Q_{\xi_\alpha}=\dot R_\alpha^{\,p},
\]
and put $\gamma_{\beta+1}=\gamma_\beta+\delta$. Here $\dot R_\alpha$
is a name over $D_M*\dot P_\alpha^{\,p}$, where $P_\alpha$ is the
iteration of the iterands $\dot R_\zeta$, $\zeta<\alpha$, and $P_0$
is trivial. By induction on
$\alpha$, the map
\[
 r\longmapsto\bigl(r\restriction M,\ \langle r(\xi_\zeta):\zeta<\alpha\rangle\bigr)
\]
shows that $D_{m(\xi_\alpha)}$ is forcing equivalent to
$D_M*\dot P_\alpha^{\,p}$; this is verified by the argument in
the proof of Proposition~\ref{prop:native-full}, applied to this
block over $D_M$. We use this equivalence to regard
$\dot R_\alpha^{\,p}$ as a $D_{m(\xi_\alpha)}$-name.

In all other cases, append a single coordinate of level $1$ with empty
memory and trivial iterand, and put
$\gamma_{\beta+1}=\gamma_\beta+1$.

Finally, put $D=D_{\aleph_{\omega+2}}$.
\end{definition}

In Case~2, the forcing $P$ and the name $\dot E$ are given by their
codes. Coordinates of levels below $i$ may occur in these codes, but
only as ordinals; they need not belong to $M$. In both cases $p$
forces every $\dot R_\alpha$ to be a $<\omega_i$-compact,
$\omega_i$-linked preorder on $\omega_{i+1}$ with a linking function
(Definition~\ref{def:small-iterations}, Theorem~\ref{thm:term}).
Whether $F(\beta)$ codes a quadruple with the properties
\ref{rec:M}--\ref{rec:cases} is a statement about the forcing relation
of $D_M$, which has already been defined at stage $\beta$.

The requirements of Definition~\ref{def:memory-scheme} hold for every
block which is appended in this way.
\begin{itemize}
\item[\ref{mem:closed}] holds because $M$ is closed and the memories
within a block increase.
\item[\ref{mem:up}] holds because the block has level $i$ and
$M\subseteq\bigcup_{k\geq i}I_k$.
\item[\ref{mem:small}] holds because the block adds at most
$|\delta|\leq\omega_{i+1}$ coordinates of level $i$.
\end{itemize}
In particular $|m(\xi)|\leq\aleph_\omega$ for every $\xi$. Every block
has cardinality at most $\omega_{i+1}<\aleph_\omega$ for some $i$.
Since $\aleph_{\omega+2}$ is regular, it follows that
\[
 \gamma_\beta<\aleph_{\omega+2}\quad(\beta<\aleph_{\omega+2}),
 \qquad \sup_{\beta<\aleph_{\omega+2}}\gamma_\beta
       =\aleph_{\omega+2}.
\]

The support of a condition has cardinality at most $\aleph_\omega$ and may meet
infinitely many levels. For $j\geq1$, let $D_{\geq j}$ be the
restriction of $D$ to $\bigcup_{k\geq j}I_k$, and put
$D_{>j}=D_{\geq j+1}$. We prove below that $D_M$ is a complete
subforcing of $D$ for every closed set $M$. Given a $D$-generic
filter $G$ over $B$, we write $G_M$ for its restriction to $D_M$, and
\[
 W_j=B[G_{D_{\geq j}}]\qquad(j\geq1).
\]

The construction can be carried out uniformly in $\vec T$, in $L$.
Indeed, every element of $\mathcal B$ has a nice $C_{\vec T}$-name in
$L$, which can be coded by a set of ordinals of size at most
$\aleph_{\omega+1}$, because $|C_{\vec T}|\leq\aleph_{\omega+1}$.
Since $\aleph_{\omega+2}^{\aleph_{\omega+1}}=\aleph_{\omega+2}$ in
$L$, there are at most $\aleph_{\omega+2}$ such names. Take the
$<_L$-least enumeration of these names of length $\aleph_{\omega+2}$
in which every name occurs unboundedly often, and let $F(\beta)$ be
the value of the $\beta$-th name.
At each stage, whether Case~1 or Case~2 applies is expressed by the
forcing relation of the part of the forcing which has already been
defined. The same definable recursion therefore
produces a $C_{\vec T}$-name $\dot D_{\vec T}$ for every
$\vec T\in L$.

The results of the remainder of this section use only the properties
of the iterands, the supports and the memories stated in
Definitions~\ref{def:memory-scheme} and~\ref{def:coding-conditions}.
They therefore hold for every iteration with partial memory, in
particular for the iteration constructed up to any stage of the
recursion.

\subsection{Complete subforcings and closure}\label{sec:quotients}

\begin{lemma}\label{lem:native-structural}
Consider an iteration with partial memory, and let $K\subseteq M$ be
closed sets of coordinates.
\begin{enumerate}
\item The restriction map $p\mapsto p\restriction K$ is a projection
from $D_M$ onto $D_K$. More precisely, if $p\in D_M$ and
$q\leq p\restriction K$ in $D_K$, then
$p'=q\cup(p\restriction(M\setminus K))$ is a condition in $D_M$ with
$p'\leq p$ and $p'\restriction K=q$. In particular $D_K$ is a
complete subforcing of $D_M$.
\item Suppose that all coordinates in $M$ have level at least $j$.
Let $\langle p_\zeta:\zeta<\eta\rangle$ be a descending sequence in
$D_M$, where $\eta<\omega_j$, and let $r\in D_K$ be a lower bound of
$\langle p_\zeta\restriction K:\zeta<\eta\rangle$. Then
$\langle p_\zeta:\zeta<\eta\rangle$ has a lower bound $q\in D_M$ with
$q\restriction K=r$. In particular $D_M$ is $<\omega_j$-closed.
\end{enumerate}
\end{lemma}
\begin{proof}
We argue by induction on the length of the iteration, simultaneously
for all closed sets of coordinates.

(1) Let $p'$ be as in the statement. For $\xi\in K$ we have
$m(\xi)\subseteq K$, so $p'\restriction m(\xi)=q\restriction m(\xi)$
forces $q(\xi)\mathrel{\dot\leq_\xi}p(\xi)$, because
$q\leq p\restriction K$. For $\xi\in M\setminus K$ we have
$p'(\xi)=p(\xi)$, and the required inequality holds because
$\dot\leq_\xi$ is forced to be reflexive. Thus
$p'\leq p$ and $p'\restriction K=q$. The restriction map preserves
the order, again because $m(\xi)\subseteq K$ for $\xi\in K$. The
order of $D_M$ is transitive: if $r\leq q\leq p$, then
$r\restriction m(\xi)\leq q\restriction m(\xi)$ for every $\xi\in M$,
so $r\restriction m(\xi)$ forces both
$r(\xi)\mathrel{\dot\leq_\xi}q(\xi)$ and
$q(\xi)\mathrel{\dot\leq_\xi}p(\xi)$, and $\dot\leq_\xi$ is forced to
be transitive.

(2) Define $q$ by $q\restriction K=r$ and, for $\xi\in M\setminus K$
in the union of the supports of the conditions $p_\zeta$,
\[
 q(\xi)=b_{\xi,\eta}
             \bigl(\langle p_\zeta(\xi):\zeta<\eta\rangle\bigr).
\]
Since $\eta<\omega_j$ and all coordinates in $M$ have level at
least $j$, the regularity of the cardinals $\omega_k$ shows that
$|\supp(q)\cap I_k|<\omega_k$ for every $k$ and that
$q(\xi)\in\mathcal N_\xi$ for every $\xi$. Thus $q\in D_M$. We show
by induction on $\xi\in M$ that
$q\restriction m(\xi)\leq p_\zeta\restriction m(\xi)$ and
$q\restriction m(\xi)\Vdash q(\xi)\mathrel{\dot\leq_\xi}p_\zeta(\xi)$
for all $\zeta<\eta$. For $\xi\in K$ this holds because
$m(\xi)\subseteq K$ and $r$ is a lower bound of the restrictions to
$K$. Let $\xi\in M\setminus K$. By the induction hypothesis,
$q\restriction m(\xi)$ is below all
$p_\zeta\restriction m(\xi)$. Hence it forces the sequence
$\langle p_\zeta(\xi):\zeta<\eta\rangle$ to be descending, so that
$b_{\xi,\eta}$ yields a lower bound of this sequence. Therefore
$q\leq p_\zeta$ for all $\zeta<\eta$, and $q\restriction K=r$.
\end{proof}

\begin{lemma}\label{lem:quotient-closed}
\begin{enumerate}
\item Let $\tau$ be regular, let $R$ be a $<\tau$-closed forcing and
let $G_R$ be $R$-generic over $V$. Then every $X\subseteq G_R$ of
size $<\tau$ in $V[G_R]$ has a lower bound in $G_R$.
\item Let $K\subseteq M$ be closed sets of coordinates of levels at
least $j$. Then $D_M/G_K$ is $<\omega_j$-closed in $B[G_K]$.
\end{enumerate}
\end{lemma}
\begin{proof}
(1) We have $X\in V$, since $R$ adds no sequences of ordinals of
length $<\tau$. Choose $r\in G_R$ with
$r\Vdash\check X\subseteq\dot G_R$. For each $a\in X$, the set
$\{q\leq r:q\leq a\}$ is dense below $r$. By $<\tau$-closure, the
intersection of these sets is dense below $r$ and hence meets $G_R$.

(2) Let $\langle p_\zeta:\zeta<\eta\rangle$
be a descending sequence in $D_M/G_K$, where $\eta<\omega_j$. Since
$D_K$ is $<\omega_j$-closed, the sequence
belongs to $B$, and by (1) there is $r\in G_K$ with
$r\leq p_\zeta\restriction K$ for all $\zeta<\eta$.
Lemma~\ref{lem:native-structural}(2) gives a lower bound $q$ of the
sequence with $q\restriction K=r$. Then $q\in D_M/G_K$.
\end{proof}

\subsection{The quotients at level \texorpdfstring{$i$}{i}}
\label{subsec:level-quotient}

Fix $i\geq1$ and put $\kappa=\omega_i$, $\sigma=\omega_{i+1}$,
\[
 U_i=\bigcup_{k>i}I_k,\qquad R=D_{>i}=D_{U_i}.
\]
By Lemma~\ref{lem:native-structural}, $R$ is
$<\sigma$-closed. Let $G_R$ be $R$-generic over $B$ and work in
$W=B[G_R]$; thus $W=W_{i+1}$. We call $D_{\geq i}/G_R$ the
\emph{quotient of $D$ at level $i$}.

For an ordinal $\gamma$ which is at most the length of the iteration,
let $S_\gamma=S_{i,\gamma}\in B$ be the set of all conditions of $D$
whose support is contained in $I_i\cap\gamma$. In $W$, define a
preorder on $S_\gamma$ by
\[
 q\preceq_{G_R}p
 \quad\Longleftrightarrow\quad
 \exists r\in G_R\ (r\cup q\leq_D p).
\]
By Lemma~\ref{lem:native-structural}(1), the preorder
$(S_\gamma,\preceq_{G_R})$ is forcing equivalent to the quotient
$D_{U_i\cup(I_i\cap\gamma)}/G_R$: a condition $p$ of this quotient is
equivalent to its restriction $p\restriction I_i\in S_\gamma$, since
$p\restriction U_i\in G_R$.

More generally, for $C\subseteq I_i$ in $W$ with $U_i\cup C$ closed,
use the conditions supported on $C$ with order $\preceq_{G_R}$.
We denote this forcing by $D_{U_i\cup C}/G_R$ up to forcing
equivalence. A restriction to $C$ belongs to $B$, since it is coded
by a sequence of length $<\kappa$ and $R$ adds no such sequences.
The amalgamation of Lemma~\ref{lem:native-structural}(1) therefore
gives a complete subforcing also for $C\notin B$.

The lower-bound construction of Lemma~\ref{lem:native-structural}(2)
also works for $\preceq_{G_R}$. For a short descending sequence
and a lower bound of its restrictions to $C$, first take a common
lower bound in $G_R$ of the witnesses to the order comparisons,
using Lemma~\ref{lem:quotient-closed}(1). Then use the prescribed
values on $C$ and the coordinatewise lower bounds off $C$.
All the data on the union of the level-$i$ supports belong to $B$.
Thus Lemma~\ref{lem:quotient-closed}(2) applies over $W$ to these
restrictions.

\begin{lemma}\label{lem:iterand-preserved}
Let $\xi\in I_i$ and put $A=U_i\cup(I_i\cap\xi)$. Then $A$ is closed
and contains $m(\xi)$, and the evaluation of $\dot Q_\xi$ by
$G_{m(\xi)}$ is $<\kappa$-compact and $\kappa$-linked in $B[G_A]$,
with the same linking function and the same lower bounds.
\end{lemma}
\begin{proof}
The set $A$ is closed by~\ref{mem:up}, and $m(\xi)\subseteq A$. By
Lemma~\ref{lem:quotient-closed}(2), $B[G_A]$ is an extension of
$B[G_{m(\xi)}]$ by a $<\kappa$-closed forcing.
Hence every family of fewer than $\kappa$ conditions of the iterand,
which is a family of ordinals below $\sigma$, already belongs to
$B[G_{m(\xi)}]$. If such a family is centered in $B[G_A]$, it is
centered in $B[G_{m(\xi)}]$, where it has a lower bound. Similarly,
conditions with the same value of the linking function are
compatible in $B[G_{m(\xi)}]$ and hence in $B[G_A]$.
\end{proof}

\begin{proposition}\label{prop:native-full}
For every $\gamma$, the preorder
$(S_\gamma,\preceq_{G_R})$ is forcing equivalent in $W$ to the
$<\kappa$-support iteration $P_\gamma$ of the iterands
$\dot Q_\xi$, $\xi\in I_i\cap\gamma$; more precisely, the
coordinatewise interpretation of the elements of $S_\gamma$ as
conditions of $P_\gamma$ preserves and reflects the order and has
dense range. The forcing $P_\gamma$ satisfies the stationary
$\sigma$-chain condition.
\end{proposition}
\begin{proof}
We argue by induction on $\gamma$, constructing $P_\gamma$ and the
interpretation of $S_\gamma$ in $P_\gamma$ simultaneously. Let
$\xi\in I_i\cap\gamma$. The forcing $D_{m(\xi)}$ is a complete
subforcing of $D_{U_i\cup(I_i\cap\xi)}$, and by the induction
hypothesis the quotient of the latter by $G_R$ is forcing equivalent
to $P_\xi$. Hence $\dot Q_\xi$ can be regarded as a $P_\xi$-name in
$W$. By Lemma~\ref{lem:iterand-preserved}, it is forced to be
$<\kappa$-compact and $\kappa$-linked. A condition $p\in S_\gamma$ is
interpreted as the condition of $P_\gamma$ whose value at $\xi$ is the
name $p(\xi)$, regarded as a $P_\xi$-name.

\emph{The order.} By induction, the interpretation preserves
$\preceq_{G_R}$. Conversely, let $p,q\in S_\gamma$ with
$p\leq_{P_\gamma}q$, and put $A=\supp(p)\cup\supp(q)$. Fix
$\xi\in A$, write
\[
 K=m(\xi)\cap U_i,\qquad C_\xi=m(\xi)\cap I_i,
\]
and let $\theta_\xi$ be the statement
$p(\xi)\leq_{\dot Q_\xi}q(\xi)$ of the forcing language of
$D_{m(\xi)}$. The set $K$ is closed, and $G_K=G_R\cap D_K$ is
$D_K$-generic over $B$. In $B$, let $E$ be the set of all
$s\in D_K$ for which one of the following holds:
\begin{enumerate}[label=(\roman*)]
\item $s\cup(p\restriction C_\xi)\Vdash_{D_{m(\xi)}}\theta_\xi$;
\item there is $a\in D_{m(\xi)}$ with
$a\leq s\cup(p\restriction C_\xi)$, $a\restriction K=s$ and
$a\Vdash_{D_{m(\xi)}}\neg\theta_\xi$.
\end{enumerate}
The set $E$ is dense in $D_K$. Indeed, if $s_0\in D_K$ does not
satisfy (i), choose $a\leq s_0\cup(p\restriction C_\xi)$ which forces
$\neg\theta_\xi$; then $s=a\restriction K\leq s_0$ satisfies (ii).
Choose $s\in G_K\cap E$.

We claim that $s$ satisfies (i). Otherwise let $a$ be as in (ii), and
let $p^*\in S_\xi$ be the condition which agrees with $a$ on $C_\xi$
and with $p$ on $(I_i\cap\xi)\setminus C_\xi$. Since
$a\leq s\cup(p\restriction C_\xi)$ and $a\restriction K=s$, we have
$s\cup p^*\leq_Dp\restriction\xi$, and therefore
$p^*\preceq_{G_R}p\restriction\xi$. By the induction hypothesis for
$\xi$, the interpretation of $p^*$ is below that of
$p\restriction\xi$ in $P_\xi$, and every $P_\xi$-generic filter over
$W$ which contains it induces a $D_{m(\xi)}$-generic filter over $B$
which contains $s\cup(p^*\restriction C_\xi)=a$. In the
corresponding extension $\theta_\xi$ fails. On the other hand,
$p\leq_{P_\gamma}q$ implies that $p\restriction\xi$ forces
$\theta_\xi$ over $P_\xi$, a contradiction.

Hence for every $\xi\in A$ there is $s_\xi\in G_R$ such that
\[
 s_\xi\cup(p\restriction C_\xi)
 \Vdash_{D_{m(\xi)}}p(\xi)\leq_{\dot Q_\xi}q(\xi).
\]
Since $|A|<\kappa$, Lemma~\ref{lem:quotient-closed}(1) gives
$r\in G_R$ below all $s_\xi$. Then $r\cup p\leq_Dq$, so
$p\preceq_{G_R}q$.

\emph{Density.} We show that for every $p\in P_\gamma$ there is
$q\in S_\gamma$ whose interpretation is below $p$. If
$\gamma=\xi+1$ and $\xi\in I_i$, first extend $p\restriction\xi$ to a
condition which decides $p(\xi)=\check c$ for some $c<\sigma$. By
the induction hypothesis there is $q\in S_\xi$ below this
condition, and then $q\cup\{(\xi,\check c)\}\in S_\gamma$ is below
$p$; recall that $1\Vdash\check c\in\dot Q_\xi$. If $\gamma=\xi+1$
and $\xi\notin I_i$, there is nothing to show. If $\gamma$ is a
limit ordinal with $\cf(\gamma)\geq\kappa$, then
$\sup\supp(p)<\gamma$, and the induction hypothesis applies.

Now let $\gamma$ be a limit ordinal with $\tau=\cf(\gamma)<\kappa$.
Choose a continuous increasing sequence
$\langle\gamma_\eta:\eta\leq\tau\rangle$ with
$\gamma_0=0$ and $\gamma_\tau=\gamma$, and construct
\[
 \begin{gathered}
 q_\eta\in S_{\gamma_\eta},\qquad
 q_\eta\leq_{P_{\gamma_\eta}}p\restriction\gamma_\eta
 \quad(\eta\leq\tau),\\
 q_\zeta\preceq_{G_R}q_\eta\quad(\eta<\zeta\leq\tau).
 \end{gathered}
\]
Start with the trivial condition $q_0$. At successor steps apply the
induction hypothesis to the condition
$q_\eta\cup(p\restriction[\gamma_\eta,\gamma_{\eta+1}))$ of
$P_{\gamma_{\eta+1}}$; the resulting condition $q_{\eta+1}$ is below
$q_\eta$ in $P_{\gamma_{\eta+1}}$ and hence
$q_{\eta+1}\preceq_{G_R}q_\eta$ by the first part of the proof.
At a limit $\nu\leq\tau$, use Lemma~\ref{lem:quotient-closed}(1) to
choose $r\in G_R$ such that $r\cup q_\zeta\leq_Dq_\eta$ for all
$\eta<\zeta<\nu$. Then
$\langle r\cup q_\eta:\eta<\nu\rangle$ is a descending sequence in
$D$, and it belongs to $B$ because $R$ adds no sequences of
length $<\sigma$ of elements of $B$.
Lemma~\ref{lem:native-structural}(2) yields a condition $s$ with
\[
 s\restriction U_i=r,\qquad
 s\leq_Dr\cup q_\eta\quad(\eta<\nu).
\]
Put $q_\nu=s\restriction I_i\in S_{\gamma_\nu}$. Then
$q_\nu\preceq_{G_R}q_\eta$ for all $\eta<\nu$, and since
$\gamma_\nu=\sup_{\eta<\nu}\gamma_\eta$ we obtain
$q_\nu\leq_{P_{\gamma_\nu}}p\restriction\gamma_\nu$.
In particular, $q_\tau\leq_{P_\gamma}p$.

\emph{The chain condition.} Since $R$ is $<\sigma$-closed,
$\kappa^{<\kappa}=\kappa$ and $\sigma=\kappa^+$ hold in $W$. The
iterands of $P_\gamma$ are forced to be $<\kappa$-compact and
$\kappa$-linked, so Corollary~\ref{cor:CDMMS} applies to $P_\gamma$.
\end{proof}

In Definition~\ref{def:coding-recursion}, a $\Gamma_i$-iteration is
appended as a block consisting of its iterands, so the proposition
applies directly to these coordinates.

\begin{definition}\label{def:PC}
Let $i\geq1$ and work in $W=W_{i+1}$. Let $\mathcal C_i$ be the set
of all $C\subseteq I_i$ such that $|C|\leq\omega_{i+1}$ and
$U_i\cup C$ is closed. For $C\subseteq I_i$ such that $U_i\cup C$ is
closed, put
\[
 P_C=D_{U_i\cup C}/G_{D_{>i}},
\]
let $G_C$ be the $P_C$-generic filter induced by the $D$-generic
filter, and let $N_C=W[G_C]$.
\end{definition}

\begin{remark}\label{rem:relativecore}
A set $C\in\mathcal C_i$ need not be closed and need not belong to
$B$: the memories of its elements contain coordinates of higher
levels, and $C$ is only required to belong to $W$. As shown above,
$P_C$ is a complete subforcing of the quotient at level $i$. All
statements below about the sizes of the forcings $P_C$ and about the
genericity of filters on them are made in $W$.
\end{remark}

\subsection{Size and preservation of cardinals}

\begin{lemma}\label{lem:size-cc-verified}
The forcing $D$ of
Definitions~\ref{def:coding-recursion} and~\ref{def:coding-conditions}
satisfies
\[
 |D|\leq\aleph_{\omega+2},\qquad D\text{ is }\aleph_{\omega+2}\text{-cc}.
\]
\end{lemma}
\begin{proof}
Since $B$ satisfies GCH,
\[
 \begin{gathered}
 \aleph_\omega^{\aleph_\omega}=2^{\aleph_\omega}=\aleph_{\omega+1},\\
 \aleph_{\omega+1}^{\aleph_\omega}=\aleph_{\omega+1},\qquad
 \aleph_{\omega+2}^{\aleph_\omega}=\aleph_{\omega+2},\\
 |\mathcal N_\xi|\leq\omega_{i+1}^{<\omega_i}
       =\omega_{i+1}<\aleph_\omega\quad(\xi\in I_i).
 \end{gathered}
\]
Since $|\supp(p)|\leq\aleph_\omega$ for every condition $p$, we obtain
$|D|\leq\aleph_{\omega+2}^{\aleph_\omega}\cdot\aleph_\omega^{\aleph_\omega}=\aleph_{\omega+2}$.

Given $\langle p_\alpha:\alpha<\aleph_{\omega+2}\rangle$ in $D$, put
\[
 \begin{gathered}
 A_\alpha=\supp(p_\alpha),\qquad
 U_\alpha=\bigcup_{\beta<\alpha}A_\beta,\\
 S=\{\alpha<\aleph_{\omega+2}:\cf(\alpha)=\aleph_{\omega+1}\}.
 \end{gathered}
\]
For $\alpha\in S$, choose $f(\alpha)<\alpha$ with
$A_\alpha\cap U_\alpha\subseteq U_{f(\alpha)}$; this is possible
since $|A_\alpha\cap U_\alpha|\leq\aleph_\omega<\cf(\alpha)$.
Fodor's lemma gives a stationary $S'\subseteq S$ with
$f[S']=\{\gamma\}$. As
$|U_\gamma|\leq\aleph_{\omega+1}$ and
$|[U_\gamma]^{\leq\aleph_\omega}|\leq\aleph_{\omega+1}^{\aleph_\omega}=\aleph_{\omega+1}<\aleph_{\omega+2}$,
there are $H\in[S']^{\aleph_{\omega+2}}$ and $a$ such that
$A_\alpha\cap U_\alpha=a$ for all $\alpha\in H$. Hence
\[
 a\subseteq A_\alpha\cap A_\beta
   \subseteq A_\beta\cap U_\beta=a
 \quad(\alpha<\beta,\ \alpha,\beta\in H).
\]
There are at most $\aleph_\omega^{|a|}\leq\aleph_{\omega+1}$ restrictions
$p_\alpha\restriction a$. Choose $K\in[H]^{\aleph_{\omega+2}}$ and $u$ with
$p_\alpha\restriction a=u$ for all $\alpha\in K$.
For distinct $\alpha,\beta\in K$, the union $r=p_\alpha\cup p_\beta$
is a function and
\[
 |\supp(r)\cap I_i|
 \leq|A_\alpha\cap I_i|+|A_\beta\cap I_i|
 <\omega_i\quad(i\geq1).
\]
Since $r(\xi)=p_\epsilon(\xi)$ whenever
$\epsilon\in\{\alpha,\beta\}$ and $\xi\in A_\epsilon$,
every inequality required for $r\leq_Dp_\epsilon$ is an instance of
reflexivity and is forced by the trivial condition; hence
$r\in D$ and $r\leq_Dp_\alpha,p_\beta$.
\end{proof}

\begin{lemma}\label{lem:cardinals-verified}
The forcing $D$ preserves all cardinals and adds no reals.
\end{lemma}
\begin{proof}
By Lemma~\ref{lem:native-structural}, $D=D_{\geq1}$ is countably
closed, so it adds no reals. For the preservation of cardinals we
apply Lemma~\ref{lem:cardinals-general}. Let $j\geq1$. Then $D$ is
forcing equivalent to $D_{\geq j}*(D/\dot G_{D_{\geq j}})$, and
$D_{\geq j}$ is $<\omega_j$-closed by
Lemma~\ref{lem:native-structural}. The quotient
$D/G_{D_{\geq j}}$ is the iteration of the quotients at the levels
$j-1,j-2,\ldots,1$. By Proposition~\ref{prop:native-full}, the
quotient at level $k$ is $\omega_{k+1}$-cc and hence $\omega_j$-cc
for $k<j$. Each of them preserves $\omega_j$, so their finite
iteration is $\omega_j$-cc. Finally $D$ is $\aleph_{\omega+2}$-cc by
Lemma~\ref{lem:size-cc-verified}.
\end{proof}

The next lemma shows that a $\Gamma_i$-iteration of the ground model
remains a $\Gamma_i$-iteration, up to forcing equivalence, after a
sufficiently closed forcing.

\begin{lemma}\label{lem:old-iteration-verified}
Let $\kappa>\omega$ be regular, $\sigma=\kappa^+$ and
$\kappa^{<\kappa}=\kappa$. Suppose $R$ is $<\kappa$-closed and
preserves $\sigma$. Let
$\langle Q_\alpha,\dot S_\alpha:\alpha<\delta\rangle\in V$ be a
$<\kappa$-support iteration such that
\[
 Q_\alpha\Vdash
 ``\dot S_\alpha=(\check\sigma,\dot{\leq}_\alpha)
 \text{ is }<\kappa\text{-compact and }\kappa\text{-linked},''
\]
given together with names for linking functions and for lower bounds
of centered families. Let $G_R$ be $R$-generic over $V$, and let
$\langle\widehat Q_\alpha,\dot S_\alpha:\alpha<\delta\rangle$ be the
$<\kappa$-support iteration defined in $V[G_R]$ from the same names
for the orders of the iterands. Then, in $V[G_R]$:
\begin{enumerate}
\item for every $\alpha\leq\delta$, the coordinatewise interpretation
$e_\alpha:Q_\alpha^V\to\widehat Q_\alpha$ preserves and reflects the
order and has dense range; in particular $Q_\delta^V$ is forcing
equivalent to $\widehat Q_\delta$;
\item the iterands of $\widehat Q_\delta$ are forced to be
$<\kappa$-compact and $\kappa$-linked, with the same linking functions
and lower bounds;
\item if $\kappa=\omega_i$, $|\delta|\leq\sigma$, and
$D_\alpha\subseteq Q_\alpha$ are dense sets of size at most $\sigma$
in $V$, then
$\langle\widehat Q_\alpha,\dot S_\alpha:\alpha<\delta\rangle$ is a
$\Gamma_i$-iteration with dense sets $e_\alpha[D_\alpha]$ of size at
most $\sigma$. In particular $Q^V_\delta$ is forcing equivalent to
a member of $\Gamma_i$.
\end{enumerate}
\end{lemma}
\begin{proof}
Write $W=V[G_R]$. Since $R$ and every $Q_\alpha$ add no sequences of
ordinals of length $<\kappa$,
\[
 R\Vdash Q_\alpha^V\text{ is }<\kappa\text{-closed},\qquad
 Q_\alpha\Vdash R^V\text{ is }<\kappa\text{-closed}.
\]
We construct the maps $e_\alpha$ by induction on $\alpha$. For
ground-model names, the forcing relation for statements about the
orders of the iterands is absolute between $V$ and $W$, so the maps
$e_\alpha$ preserve and reflect the order. In the product
$R^V\times Q_\alpha^V$, the closure of $R^V$ over
$V[G_{Q_\alpha}]$ gives
\[
 [\sigma]^{<\kappa}\cap V[G_{Q_\alpha}][G_R]
   =[\sigma]^{<\kappa}\cap V[G_{Q_\alpha}].
\]
The order of the iterand and its linking function on $\sigma$ are
unchanged, and whether finitely many conditions have a common lower
bound is absolute. Thus compactness and linkedness are preserved,
which proves (2).

For the density of the range of $e_{\alpha+1}$, let
$p\in\widehat Q_{\alpha+1}$. Choose $p'\leq p\restriction\alpha$
which decides $p(\alpha)=\check c$ for some $c<\sigma$, and
$q\in Q_\alpha^V$ with $e_\alpha(q)\leq p'$. Then
$e_{\alpha+1}(q^\frown\langle\check c\rangle)\leq p$.
Let $\gamma$ be a limit ordinal and $p\in\widehat Q_\gamma$. If
$\cf^W(\gamma)\geq\kappa$, then
$\sup\supp(p)<\gamma$, and the induction hypothesis applies. Otherwise
let $\tau=\cf^W(\gamma)<\kappa$, choose a continuous increasing
sequence $\langle\gamma_\eta:\eta\leq\tau\rangle$ with $\gamma_0=0$
and $\gamma_\tau=\gamma$, and construct
\[
 \begin{gathered}
 q_\eta\in Q_{\gamma_\eta}^V,\qquad
 e_{\gamma_\eta}(q_\eta)\leq p\restriction\gamma_\eta
 \quad(\eta\leq\tau),\\
 q_\zeta\leq q_\eta\quad(\eta<\zeta\leq\tau).
 \end{gathered}
\]
Start with the trivial condition $q_0$, and at successor steps apply
the induction hypothesis to the condition
$e_{\gamma_\eta}(q_\eta)\cup
(p\restriction[\gamma_\eta,\gamma_{\eta+1}))$.
At limits $\nu\leq\tau$, we have
$\langle q_\eta:\eta<\nu\rangle\in V$, since $R$ adds no such
sequences; let $q_\nu$ be a lower bound in the $<\kappa$-closed
forcing $Q^V_{\gamma_\nu}$. Comparing coordinates, we obtain
$e_\gamma(q_\tau)\leq p$. Since the maps $e_\alpha$ reflect the order
and have dense range, they also preserve incompatibility. This proves
(1).

For (3), note that
$W\models\kappa^{<\kappa}=\kappa\land\sigma=\kappa^+$, so that
Corollary~\ref{cor:CDMMS} applies to $\widehat Q_\delta$ in $W$ and
gives the $\sigma$-chain condition. For each
dense set $D_\alpha\subseteq Q_\alpha$, the set
$e_\alpha[D_\alpha]$ is dense in $\widehat Q_\alpha$, and
$|\delta|,|e_\alpha[D_\alpha]|\leq\sigma$ in $W$. Restricting
$e_\delta$ to the conditions below a given condition gives the
assertion for the restrictions of $Q^V_\delta$.
\end{proof}

\section{The bookkeeping and a first-order characterization of \texorpdfstring{$\Gamma_i$}{Gamma-i}}\label{sec:names}

The forcing $D$ and its bookkeeping function $F$ are now fixed. We
first show that the names which are needed later are names over
$D_M$ for bounded closed sets $M$ of coordinates, and are therefore
handled at some stage of the recursion. We then show that the
$\Gamma_i$-iterations can be characterized by a first-order formula
over $H_{\omega_{i+1}}$; this is used in the reduction of complexity
in Section~\ref{sec:truthdrop}.

\subsection{Names depend on few coordinates}

Throughout, for a fixed $i\geq1$ we write $U_i=\bigcup_{k>i}I_k$ and
$W=W_{i+1}=B[G_{D_{>i}}]$, and we use the notation
$\mathcal C_i$, $P_C$, $G_C$ and $N_C=W[G_C]$ of
Definition~\ref{def:PC}.

We shall use Easton's lemma in the following form.

\begin{lemma}\label{lem:easton-locality}
Let $\sigma$ be regular uncountable, let $R$ be
$<\sigma$-closed and let $P$ be $\sigma$-cc, and let $G_R\times G_P$
be generic for $R\times P$ over $V$. Then every sequence of ordinals
of length $<\sigma$ in $V[G_R][G_P]$ belongs to $V[G_P]$, and hence
\[
 H_\sigma^{V[G_P]}=H_\sigma^{V[G_R][G_P]}.
\]
Consequently, first-order statements about $H_\sigma$ with
predicates from $V[G_P]$ have the same truth value in $V[G_P]$ and
in $V[G_R][G_P]$.
\end{lemma}
\begin{proof}
Let $r_0\in R$ and let $\dot f$ be an $R\times P$-name for a
function from some $\eta<\sigma$ to the ordinals. For every
$\zeta<\eta$ in turn, we construct by recursion on $\xi$ conditions
$p_{\zeta,\xi}\in P$ and $r_{\zeta,\xi}\leq r_0$ and ordinals
$\alpha_{\zeta,\xi}$ such that, for fixed $\zeta$, the conditions
$p_{\zeta,\xi}$ are pairwise incompatible, the conditions
$r_{\zeta,\xi}$ are
descending throughout the whole construction, and
\[
 (r_{\zeta,\xi},p_{\zeta,\xi})\Vdash
 \dot f(\zeta)=\check\alpha_{\zeta,\xi}.
\]
As long as there is a condition $q$ which is incompatible with all
$p_{\zeta,\xi'}$ chosen so far, take a lower bound $r$ of all
conditions of $R$ chosen so far, and extend $(r,q)$ to a condition
which decides $\dot f(\zeta)$. By the chain condition, the
construction for $\zeta$ stops after $\theta_\zeta<\sigma$ steps,
with a maximal antichain
$A_\zeta=\{p_{\zeta,\xi}:\xi<\theta_\zeta\}$. Since $\sigma$ is
regular, $\sum_{\zeta<\eta}\theta_\zeta<\sigma$, so the required
lower bounds exist by the closure of $R$, and there is a lower bound
$r_*$ of all the conditions $r_{\zeta,\xi}$. The antichains $A_\zeta$
and the ordinals $\alpha_{\zeta,\xi}$ define a $P$-name $\dot g\in V$
with $(r_*,1_P)\Vdash\dot f=\dot g$. Since $r_0$ was arbitrary, this
proves the first assertion. Every element of $H_\sigma$ is coded by a
well-founded relation on an ordinal below $\sigma$; this gives the
equality of the two sets, and the last assertion follows.
\end{proof}

\begin{lemma}\label{lem:memory-size-full}
Let $A\in B$ be a set of coordinates with
$|A\cap I_k|\leq\omega_{k+1}$ for every $k$. Then
$\bar A=A\cup\bigcup_{\xi\in A}m(\xi)$ is the least closed set
containing $A$, and it satisfies the same bounds. It has size at most
$\aleph_\omega$ and is bounded in $\aleph_{\omega+2}$. If $A$ is
contained in an ordinal $\alpha$, then so is $\bar A$.
\end{lemma}
\begin{proof}
Each $m(\xi)$ is closed, so $\bar A$ is closed. By~\ref{mem:up}, a
coordinate of level $k$ can belong to $m(\xi)$ only if $\xi$ has one
of the finitely many levels $1,\ldots,k$. There are at most
$\omega_{k+1}$ such coordinates $\xi\in A$, and each of them
contributes at most $\omega_{k+1}$ coordinates of level $k$. Hence
$|\bar A\cap I_k|\leq\omega_{k+1}$ and $|\bar A|\leq\aleph_\omega$.
Since $\aleph_{\omega+2}$ is regular, $\bar A$ is bounded in
$\aleph_{\omega+2}$. The last assertion holds because
$m(\xi)\subseteq\xi$.
\end{proof}

\begin{proposition}\label{prop:localization-full}
Let $i\geq1$ and $\alpha\leq\aleph_{\omega+2}$, and let
$D_{\alpha,>i}$ be the restriction of $D$ to the coordinates
$\xi<\alpha$ of level greater than $i$. Let
$u\in B[G_{D_{\alpha,>i}}]$ be a sequence of ordinals of length at
most $\omega_{i+1}$. Then there is a closed set $H\in B$ with
\[
 H\subseteq\{\xi<\alpha:i(\xi)>i\},\qquad
 |H\cap I_k|\leq\omega_{k+1}\quad(k>i),
\]
such that $u\in B[G_H]$. Moreover, the range of $u$ is contained in a
set in $B$ of size at most $\omega_{i+1}$.
\end{proposition}
\begin{proof}
Put $R=D_{\alpha,>i+1}$. By Lemma~\ref{lem:native-structural}, $R$ is
$<\omega_{i+2}$-closed, and by Proposition~\ref{prop:native-full}
(for level $i+1$) the quotient $D_{\alpha,>i}/G_R$ is forcing
equivalent to $(S_{i+1,\alpha},\preceq_{G_R})$ and is
$\omega_{i+2}$-cc. Recall that $S_{i+1,\alpha}\in B$.

In $B[G_R]$, choose a nice name for $u$ over this quotient. By the
chain condition it consists of at most $\omega_{i+1}$ triples, each
consisting of an argument, a condition in $S_{i+1,\alpha}$ and an
ordinal value. Replace the conditions by their indices in a fixed
enumeration of $S_{i+1,\alpha}$ which belongs to $B$. The resulting
sequence of ordinals has length at most $\omega_{i+1}$. Since $R$ is
$<\omega_{i+2}$-closed, this sequence belongs to
$B$, although its evaluation $u$ need not.

The supports of these conditions together contain at most
$\omega_{i+1}$ coordinates, all of level $i+1$. Let $H$ be the
closure of their union under the sets $m(\xi)$; by
Lemma~\ref{lem:memory-size-full}, $H$ has the required properties.
Each of these conditions is a condition of $D$ with support contained
in $H$, and it belongs to the generic filter for the quotient if and
only if it belongs to $G_H$. Thus the name can be evaluated using
$G_H$, and $u\in B[G_H]$. The set of all ordinal values which occur in
the name belongs to $B$ and has size at most $\omega_{i+1}$; it
contains the range of $u$.
\end{proof}

\begin{lemma}\label{lem:static-square-full}
Let $C\in\mathcal C_i$, let $Q\in\Gamma_i$ in $N_C$, and let $q\in Q$.
There is a bounded closed set $M\in B$ with
\[
 M\subseteq\bigcup_{k\geq i}I_k,\qquad
 |M\cap I_k|\leq\omega_{k+1}\ (k\geq i),\qquad
 C\subseteq K=M\cap I_i.
\]
There are $p\in G_M$ and a $D_M$-name $\dot R$ such that $p$ forces
$\dot R$ to code a $\Gamma_i$-iteration of positive length.
The last forcing of the iteration coded by $\dot R^{G_M}$ in
$B[G_M]$ is forcing equivalent in $N_K$ to $Q\restriction q$.
\end{lemma}
\begin{proof}
Put $\rho=\omega_i$, $\sigma=\omega_{i+1}$ and $W=W_{i+1}$.
By Proposition~\ref{prop:localization-full}, $C$ is covered by
a set in $B$ of size at most $\sigma$. Intersect this set with $I_i$
and close it under $\xi\mapsto m(\xi)\cap I_i$ to obtain $K\in B$
with $|K|\leq\sigma$ and $U_i\cup K$ closed. In $B$, fix
enumerations of $H_\sigma^B$ and of the conditions supported on $K$.
There are at most $\sigma$ such conditions by the bounds on
$\mathcal N_\xi$ and the supports. Write $P=P_K$ for the quotient,
represented as the set of these conditions with the order
$\preceq_{G_{D_{>i}}}$ of Subsection~\ref{subsec:level-quotient}.
The $<\rho$-closure and $\sigma$-cc
of $P$ give $(\ast)$ in $N_K$ by counting nice names.

By Lemma~\ref{lem:coretruth}(a), whose proof uses only the structural
properties of the quotients, $P/G_C$ is forcing equivalent to a
$<\rho$-closed, $\sigma$-cc forcing. In $N_K=W[G_K]$, take the
iteration which is defined from the names for the iterands of $Q$, as
in Lemma~\ref{lem:old-iteration-verified}, and restrict it below the
image of $q$; by the remark on cones following
Definition~\ref{def:small-iterations}, the result is again given by
a $\Gamma_i$-iteration. Choose a code $R\subseteq\sigma$
for the resulting iteration, adding a trivial iterand if needed,
and choose $\chi$ with
$\operatorname{Wit}(R,\chi)$ by Proposition~\ref{prop:eligibility-exact}.

The $\sigma$-cc of $P$ gives nice names for $R$ and $\chi$ coded
together by a sequence $u$ of at most $\sigma$ ordinals. Indeed,
a nice name for an element of $H_\sigma^{N_K}$ has a code in
$H_\sigma^W=H_\sigma^B$. For each such code use antichains of size at
most $\rho$ to decide membership in the two classes, using the fixed
enumerations for the indices of names and conditions. By
Proposition~\ref{prop:localization-full} and
Lemma~\ref{lem:memory-size-full}, choose a closed $H\in B$, $H\subseteq U_i$,
with the required bounds, containing $m(\xi)\cap U_i$ for every
$\xi\in K$, such
that $u\in Z=B[G_H]$. Put $M=H\cup K$.

In $Z$, the quotient $D_M/G_H$ is represented by the same preorder
$P$ as in $W$. For conditions $s,t$ supported on $K$, a witness
$r\in G_{D_{>i}}$ to $r\cup s\leq_D t$ can be restricted to $H$,
since $m(\xi)\subseteq M$ for $\xi\in K$. Conversely, a witness in
$G_H$ gives one in $G_{D_{>i}}$. Put $R'=D_{>i}/G_H$. It is
$<\sigma$-closed in $Z$, and $P$ is $\sigma$-cc there because its
order is the same in $W$. Applying Lemma~\ref{lem:easton-locality}
to $P\times R'$ over $Z$ gives
\[
 H_\sigma^{B[G_M]}=H_\sigma^{Z[G_K]}=H_\sigma^{N_K}.
\]
The sequence $u$ defines the same $R,\chi$ in $B[G_M]$, so the
first-order assertion $\operatorname{Wit}(R,\chi)$ holds there.
Proposition~\ref{prop:eligibility-exact} and the forcing theorem give
$p$ and a nice name $\dot R$ for a subset of $\sigma$.

Since $P$ adds no sequences of length $<\rho$ over $Z$, $R'$ is
still $<\rho$-closed after forcing with $P$. It preserves $\sigma$
by Lemma~\ref{lem:easton-locality}. By
Lemma~\ref{lem:old-iteration-verified}, the iteration given by $R$
in $B[G_M]$ is forcing equivalent to the iteration defined in $N_K$
from the same names. The antichain names in $R$ translate along the
dense embeddings of that lemma, so the latter is the iteration chosen
in $N_K$.
\end{proof}

\begin{lemma}\label{lem:term-bookkeeping-verified}
Fix $j\geq1$ and work in $W_j=B[G_{D_{\geq j}}]$.
Let $C$ be a set of coordinates of levels below $j$, with
$|C\cap I_k|\leq\omega_{k+1}$ for $k<j$, such that
$C\cup\bigcup_{k\geq j}I_k$ is closed. Represent the quotient
$D_{C\cup\bigcup_{k\geq j}I_k}/G_{D_{\geq j}}$ by the conditions
supported on $C$ with the quotient order. There are at most
$\omega_j$ such conditions by the bounds on the supports and
the $\mathcal N_\xi$. Enumerate them to give a preorder $P$ on
$\omega_j$, and let $\dot E$ be a nice $P$-name for a subset of
$\omega_{j+1}$. Then there is a closed set
$H\in B$ with
\[
 H\subseteq\bigcup_{k\geq j}I_k,\qquad
 |H\cap I_k|\leq\omega_{k+1}\ (k\geq j),\qquad\sup H<\aleph_{\omega+2},
\]
such that $P$ and $\dot E$ belong to $B[G_H]$. Moreover, there is a
stage $\beta$ of the recursion with $\sup(H\cup C)<\gamma_\beta$ at
which Case~2 of Definition~\ref{def:coding-recursion} applies with
$i=j$, $M=H$ and a name for $(P,\dot E)$; thus the term forcing
$\Term(P,\mathbb{A}_{\omega_j}(\dot E))$, as defined in $B[G_H]$, is
appended at stage $\beta$.
\end{lemma}
\begin{proof}
The forcing $P$ is $\omega_j$-cc, since the quotient is a finite
iteration of the quotients at the levels $j-1,\ldots,1$ restricted to
$C$, each of which is $\omega_j$-cc by
Proposition~\ref{prop:native-full}. The order of $P$ is a subset of
$\omega_j\times\omega_j$, and each antichain $A_\alpha$ in the nice
name $\dot E$ has size $<\omega_j$. Thus a sequence of ordinals $u$
coding $P$ and $\dot E$ has length
\[
 |u|\leq\omega_j+\sum_{\alpha<\omega_{j+1}}|A_\alpha|\leq\omega_{j+1}.
\]
In $B[G_{D_{>j}}]$, a nice name for $u$ over the quotient at level
$j$ consists of at most $\omega_{j+1}$ conditions and ordinals. Apply
Proposition~\ref{prop:localization-full} to a sequence of ordinals
coding this name, add the supports of the conditions which occur in
the name (these conditions are contained in a set in $B$ of size at
most $\omega_{j+1}$, by the last assertion of
Proposition~\ref{prop:localization-full}), and take the closure, using
Lemma~\ref{lem:memory-size-full}. These supports contain at most
$\omega_{j+1}$ coordinates of level $j$. This gives the required set
$H$, with $u\in Z=B[G_H]$ and $|H|\leq\aleph_\omega$; since
$\aleph_{\omega+2}$ is regular, $\sup H<\aleph_{\omega+2}$.

The coordinates in $C$, which have levels below $j$, occur in $u$
only as ordinals; the generic filter on these coordinates is not
needed to define $P$ and $\dot E$. A coordinate $\xi$ at which the
term forcing is appended may therefore have the memory
\[
 m(\xi)=H\subseteq\bigcup_{k\geq j}I_k.
\]
In $Z$, $P$ is the same preorder on $\omega_j$, and it is
$\omega_j$-cc there because it is $\omega_j$-cc in the larger model
$W_j$. Moreover $\omega_j^{<\omega_j}=\omega_j$ and
$\omega_{j+1}^{<\omega_j}=\omega_{j+1}$ hold in $Z$, since $D_H$ is
$<\omega_j$-closed. Theorem~\ref{thm:term} therefore applies to
\[
 T=\Term(P,\mathbb{A}_{\omega_j}(\dot E))^Z.
\]
In all later arguments we use this forcing $T$ with its conditions and
order as defined in $Z$, also after the generic filter for $P$ has
been added.

Choose $p\in G_H$ and a $D_H$-name $\dot X$ for $(P,\dot E)$ such
that $p$ forces the properties required in Case~2. The quadruple
$(j,H,p,\dot X)$ is coded by an element of $\mathcal B$, which is
equal to $F(\beta)$ for unboundedly many $\beta$. Since $H\cup C$ is
bounded, there is such a $\beta$ with $\sup(H\cup C)<\gamma_\beta$.
At this stage Case~2 applies and the term forcing $T$ is appended.
\end{proof}

\subsection{The bookkeeping handles all relevant names}\label{subsec:bookkeeping-exhaustive}

\begin{proposition}\label{prop:schedule-full}
\begin{enumerate}
\item Let $C\in\mathcal C_i$, $Q\in\Gamma_i$ in $N_C$, and $q\in Q$.
There are $M,p,\dot R$ as in Lemma~\ref{lem:static-square-full}
and a stage $\beta$ with $M\subseteq\gamma_\beta$ at which
Case~1 applies with $\dot X=\dot R$.
\item Let $j\geq1$ and let $C$, $P$ and $\dot E$ be as in
Lemma~\ref{lem:term-bookkeeping-verified}. Then the term forcing
$\Term(P,\mathbb{A}_{\omega_j}(\dot E))$ is appended at a stage after
all coordinates in $C$ and all coordinates on which $P$ and $\dot E$
depend. The generic filter induced on $P$ is generic over $W_j$, which
contains the generic filter for the term forcing.
\end{enumerate}
\end{proposition}
\begin{proof}
For a closed set $M\in B$ with $|M\cap I_k|\leq\omega_{k+1}$ for all
$k\geq1$, we have
\[
 |M|\leq\aleph_\omega,\qquad
 |D_M|\leq\aleph_\omega^{\aleph_\omega}=\aleph_{\omega+1}.
\]
A nice $D_M$-name given by maximal antichains therefore consists of
at most $\aleph_{\omega+1}$ pairs of conditions and ordinals. These
are coded by ordinals which represent conditions, coordinates below
$\aleph_{\omega+2}$, and ordinal values below $\omega_{i+1}$. Thus
every pair $(p,\dot X)$ of a condition and such a name, together with
$i$ and $M$, is coded by an element of the set $\mathcal B$.

(1) Take $M,p,\dot R$ from Lemma~\ref{lem:static-square-full}.
Strengthen $p$ in $G_M$ to decide the length of the iteration.
The code of $(i,M,p,\dot R)$ occurs unboundedly often, so it occurs
at a stage $\beta$ with $\sup M<\gamma_\beta$. The recursion
appends the iterands of the iteration given by this code.

(2) Lemma~\ref{lem:term-bookkeeping-verified}
provides a bounded closed set $H$ of coordinates of levels at least
$j$ and a stage at which Case~2 applies and the term forcing is
appended. By Lemma~\ref{lem:native-structural}, applied to the
restriction to the closed set $C\cup\bigcup_{k\geq j}I_k$, the filter
induced on $P$ is generic over $W_j$, and $W_j$ contains the generic
filter for the term forcing. Hence the projection of
Theorem~\ref{thm:term} from the product to the two-step iteration can
be applied.
\end{proof}

\subsection{Subforcings of the quotients}

\begin{lemma}\label{lem:coretruth}
Let $i\geq1$, work in $W=W_{i+1}$, and let $C\in\mathcal C_i$.
\begin{enumerate}[label=(\alph*)]
\item For every $X\subseteq H_{\omega_i}$ in $W_i$ there is
$C'\in\mathcal C_i$ with $C\subseteq C'$ and $X\in N_{C'}$. Moreover,
for all $C\subseteq C'$ in $\mathcal C_i$, the quotient
$P_{C'}/G_C$ is forcing equivalent, in $N_C$, to a
$\Gamma_i$-iteration.
\item Let $Q\in\Gamma_i$ in $N_C$ and $q\in Q$. There are
$C_0,C'\in\mathcal C_i$ with $C\subseteq C_0\subseteq C'$ and a
$Q$-generic filter $H$ over $N_{C_0}$ with $q\in H$ and
$N_{C'}=N_{C_0}[H]$.
\end{enumerate}
\end{lemma}
\begin{proof}
(a) Choose a nice name for $X$ over the quotient at level $i$. Since
$H_{\omega_i}$ has cardinality $\omega_i$ and every antichain of the
quotient has cardinality at most $\omega_i$, the name consists of at
most $\omega_i$ conditions. Let $C'$ be the closure of the union of
$C$ and the supports of these conditions under the maps
$\xi\mapsto m(\xi)\cap I_i$. As in the proof of
Lemma~\ref{lem:memory-size-full}, $|C'|\leq\omega_{i+1}$. Then
$C'\in\mathcal C_i$ and $X\in N_{C'}$.

Now let $C\subseteq C'$ be elements of $\mathcal C_i$. To describe
$P_{C'}/G_C$, we order the coordinates of $C$ first and then the
coordinates of $C'\setminus C$ in increasing order. Since
$U_i\cup C'$ is closed, every element of $m(\xi)$, for
$\xi\in C'\setminus C$, belongs to $U_i\cup C$ or precedes $\xi$ in
this order. For $\xi\in C'\setminus C$, the model obtained at this stage
lies between $B[G_{m(\xi)}]$ and $W_i$. By
Lemma~\ref{lem:quotient-closed}, $W_i$ adds no sequences of length
$<\omega_i$ over $B[G_{m(\xi)}]$, so the iterand remains
$<\omega_i$-compact and keeps its linking function.
The proof of Proposition~\ref{prop:native-full}, carried out for
conditions whose restriction to $C$ belongs to $G_C$,
shows that the quotient is forcing equivalent to the
$<\omega_i$-support iteration of the iterands at the coordinates in
$C'\setminus C$. This iteration has length of cardinality at most
$\omega_{i+1}$ and a dense subset of that cardinality, so it is a
$\Gamma_i$-iteration.

(b) Take $M,p,\dot R$ and the block from
Proposition~\ref{prop:schedule-full}(1). Let $C_0=M\cap I_i$
and let $C'$ be its union with the block. Then $C_0,C'\in\mathcal C_i$
by the bounds on $M$ and by the definition of the memories of the
block. By (a), the quotient from $N_{C_0}$ to $N_{C'}$ is given by
the iterands of the block. Since $p\in G_M$, these are the
iterands specified by the code $\dot R^{G_M}$, interpreted in $N_{C_0}$.
Lemma~\ref{lem:static-square-full} gives a forcing equivalence with
$Q\restriction q$ there. Under this equivalence, the generic filter
on the block corresponds to a $Q$-generic filter $H$ over $N_{C_0}$
with $q\in H$, and $N_{C'}=N_{C_0}[H]$.
\end{proof}

\subsection{A first-order characterization of
\texorpdfstring{$\Gamma_i$}{Gamma-i}-iterations}
\label{subsec:eligibility-justification}

Fix a regular uncountable cardinal $\rho$ and put $\sigma=\rho^+$.
We consider the following hypothesis on the universe:
\[
 (\ast)\qquad\rho^{<\rho}=\rho\quad\text{and}\quad 2^\rho=\sigma.
\]
It implies $|H_\sigma|=\sigma$ and $\sigma^\rho=\sigma$. We shall
apply the results of this subsection with $\rho=\omega_i$ in the
models $W_{i+1}$ and $N_C$, for $C\in\mathcal C_i$, where $(\ast)$
holds. We do not use them in $W_i$: forcing with the whole quotient
at level $i$ increases $2^\rho$.

In this subsection, a $\Gamma$-iteration is an iteration as in
Definition~\ref{def:small-iterations}, with $\rho$ in place of
$\omega_i$. To quantify over these iterations, we code them, and
names for them, by sets of ordinals. Let
$\langle P_\alpha,\dot Q_\alpha:\alpha<\delta\rangle$ be a
$\Gamma$-iteration, and assume $(\ast)$. Choose a wellorder
$(J,\prec)$ with
\[
 J\subseteq\sigma,\qquad \operatorname{otp}(J,\prec)=\delta.
\]
Recall that $|\delta|\leq\sigma$; the ordinal $\delta$ itself may
be larger than $\sigma$. Using the $\sigma$-chain condition and
$\sigma^\rho=\sigma$, one constructs by recursion, using nice names,
dense subsets and enumerations
\[
 e_\alpha:\sigma\twoheadrightarrow D_\alpha\subseteq P_\alpha,
 \qquad D_\alpha\text{ dense in }P_\alpha
 \quad(\alpha\leq\delta).
\]
By means of these enumerations, the conditions which occur in a name
can be replaced by ordinals. The enumerations need not be injective.
Similarly, each iterand is a preorder on $\sigma$ with largest
condition $0$.

The \emph{code} of the iteration consists of $(J,\prec)$, of names
for the orders of the iterands, for their linking functions and for
lower bounds of centered families, and of the enumerations
$e_\alpha$. A nice name is coded by the ordinals which represent the
conditions of its antichains, together with the corresponding
values. A name for an ordinal which is the value of a condition at
one coordinate consists of at most $\rho$ conditions. Since a
condition has support of size $<\rho$, its code has hereditary size
at most $\rho$ and belongs to $H_\sigma$. In the iteration code,
each value $e_\alpha(\xi)$ is given by the code of this condition,
that is, by its support and, for each coordinate of the support, by
the antichain and the values of the corresponding nice name. The
stages are indexed by the elements of $J$, and the last stage by one
additional index. The names for the orders and for the functions are
coded in the same way; where an argument is a set of size $<\rho$,
its elements are listed. Coding finite tuples of ordinals below $\sigma$ by
ordinals below $\sigma$, we take the iteration code to be a subset
of $\sigma$. A condition used in a forcing statement can be
strengthened to one in $D_\delta$. The iteration codes and the codes for names of subsets of $H_\sigma$
are classes over $H_\sigma$.

Let $\mathcal T$ be the following theory in the language of set
theory with additional constants for $\rho$, $\sigma$ and a code $R$.
It consists of $\mathsf{ZFC}^-$ (that is, ZFC without Power Set, with
Collection and Choice), of the statements that $\rho$ is regular
uncountable, $\sigma=\rho^+$, $\rho^{<\rho}=\rho$ and
$2^\rho=\sigma$, and of the statements that $R$ is the code of a
$\Gamma$-iteration: the relation $\prec$ is a wellorder, the coded
sequence is a $<\rho$-support iteration along this wellorder, the
sets enumerated by the maps $e_\alpha$ are dense, and the names for
the iterands are forced to be $<\rho$-compact preorders on $\sigma$
for which the given names are names for linking functions and for
lower bounds. All these are first-order statements about set forcing
and the components of the code. If $R$ is the code of a
$\Gamma$-iteration, then every elementary submodel of a sufficiently
large $H_\theta$ which contains $\rho$, $\sigma$ and $R$ satisfies
$\mathcal T$.

\begin{lemma}\label{lem:eligibility-prefix}
Let $M$ be a transitive model of $\mathcal T$ such that
$H_\sigma\subseteq M$ and $H_\sigma,\sigma\in M$, and such that the
constants $\rho$, $\sigma$ are interpreted by the true $\rho$,
$\sigma$. Let $\delta$ be the length of the iteration coded by $R$ in
$M$. For $\alpha\leq\delta$ let $P^M_\alpha$ be the $\alpha$-th stage
of this iteration as defined in $M$, and let $P_\alpha$ be the
$\alpha$-th stage of the $<\rho$-support iteration defined in $V$
from the names in $R$. Then for every $\alpha\leq\delta$:
\[
 \begin{gathered}
 P_\alpha^M\subseteq P_\alpha,\quad\text{and}\quad
 M\models p\leq_{P^M_\alpha}q\iff p\leq_{P_\alpha}q
 \quad(p,q\in P^M_\alpha),\\
 D\in M,\ M\models D\text{ is dense in }P_\alpha^M
 \ \Longrightarrow\ D\text{ is dense in }P_\alpha,\\
 P_\alpha\text{ is }<\rho\text{-closed and satisfies the stationary }\sigma\text{-chain condition},\\
 H_\sigma^{M[G]}=H_\sigma^{V[G]}
 \quad(G\subseteq P_\alpha\text{ generic over }V).
 \end{gathered}
\]
\end{lemma}
\begin{proof}
The relation $\prec$ on $J\subseteq\sigma$ is a wellorder in $V$: a
descending sequence $s\in{}^\omega J$ would belong to
$H_\sigma\subseteq M$, whereas $M$ satisfies that $\prec$ is
well-founded. We prove the four assertions by induction along this
wellorder. At stage $\alpha$, the induction hypothesis
$H_\sigma^{M[G_\alpha]}=H_\sigma^{V[G_\alpha]}$ gives
\[
 [\sigma]^{<\rho}\cap M[G_\alpha]
   =[\sigma]^{<\rho}\cap V[G_\alpha].
\]
The order of the iterand and its linking function on $\sigma$ are
the same in $M[G_\alpha]$ and $V[G_\alpha]$. Whether finitely many
conditions have a common lower bound is absolute, so every centered
family of size $<\rho$ in $V[G_\alpha]$ belongs to $M[G_\alpha]$ and
has a lower bound there; and conditions with the same value of the
linking function are compatible. Thus the iterands are
$<\rho$-compact and $\rho$-linked in $V[G_\alpha]$, and
Corollary~\ref{cor:CDMMS} applies to the iterations $P_\alpha$.

For the assertion on the order, let $\beta<\alpha$. By the induction
hypothesis, $P^M_\beta$ is a dense subset of $P_\beta$ with the
induced order, and every dense subset of $P^M_\beta$ in $M$ is dense
in $P_\beta$. Hence a $P^M_\beta$-name in $M$ is a $P_\beta$-name,
which gives $P^M_\alpha\subseteq P_\alpha$, and every filter $G$
which is $P_\beta$-generic over $V$ induces a filter which is
$P^M_\beta$-generic over $M$; the two filters evaluate the names in
$M$ in the same way. Let $p,q\in P^M_\alpha$, and let $\theta$ be
the statement $p(\beta)\leq q(\beta)$ about the order of the iterand
at $\beta$. Its truth value is the same in $M[G]$ and in $V[G]$.
Since every condition of $P_\beta$ has an extension in $P^M_\beta$,
it follows that $p\restriction\beta$ forces $\theta$ over $P^M_\beta$
in $M$ if and only if it forces $\theta$ over $P_\beta$ in $V$. As
the order of the iteration is defined coordinatewise from these
forcing statements, $p\leq q$ holds in $M$ if and only if it holds
in $V$.

For the density assertion, consider first a successor stage. A nice
name for an ordinal below $\sigma$ is given by a maximal antichain of
size at most $\rho$ in the dense set of the preceding stage, together
with ordinal values. Replacing the conditions by the ordinals which
represent them, we obtain a code in $H_\sigma\subseteq M$. If only
part of a maximal antichain is given, extend it to a maximal antichain
and assign the value $0$, which represents the largest condition, to
the new elements; the resulting name is forced by the trivial
condition to be a condition of the iterand. At a limit stage, code
the value at every coordinate of the support in this way. For the
resulting code $c_p$ of a condition $p$,
\[
 |\supp(p)|<\rho,\qquad
 |\operatorname{tc}(c_p)|\leq\rho,\qquad
 c_p\in H_\sigma\subseteq M.
\]
The condition coded by $c_p$ in $M$ has the same values at all
coordinates, since every ordinal below $\sigma$ represents a
condition of the iterand. Thus
$P_\alpha^M$ is dense in $P_\alpha$. Since the two orders agree on
$P^M_\alpha$, every set which is dense in $P^M_\alpha$ in the sense
of $M$ is dense in $P_\alpha$.

Finally, let $x\in H_\sigma^{V[G]}$. Choose a well-founded extensional
relation $E\subseteq\gamma^2$ and a point $e<\gamma<\sigma$ with
$x=\pi_E(e)$, where $\pi_E$ is the transitive collapse.
A nice name for $E$ consists of at most $\rho$ conditions; its code
belongs to $H_\sigma^V$, hence to $M$, and it has the same value in
$M[G]$. The relation $E$ is well-founded in $M[G]$, and by Collection
its transitive collapse exists there. By uniqueness of the collapse,
$x=\pi_E(e)\in M[G]$, and $x$ has hereditary size $<\sigma$ in
$M[G]$. Conversely, an enumeration in $M[G]$ which witnesses that a
set has hereditary size $<\sigma$ witnesses the same in $V[G]$. This
proves the equality and completes the induction.
\end{proof}

\begin{proposition}\label{prop:eligibility-exact}
Assume $(\ast)$. There is a first-order formula
$\operatorname{Wit}(R,\chi)$ over $H_\sigma$, with class
variables $R,\chi$, such that for every class $R$ over $H_\sigma$,
\[
 R\text{ is the code of a }\Gamma\text{-iteration}
 \quad\Longleftrightarrow\quad
 \exists\chi\,\operatorname{Wit}(R,\chi).
\]
The formula has no parameters; in particular it does not depend on
$\vec T$ or on the stationary sets $S^n_t$.
\end{proposition}
\begin{proof}
The formula $\operatorname{Wit}(R,\chi)$ requires $R$ to consist of ordinals
and expresses that $\chi$ codes
a set $A\subseteq\sigma$, a binary relation $E$ on $A$, a relation
$\operatorname{Sat}$, a function $j:H_\sigma\to A$ and elements
$h,s,r\in A$ with the following properties. The relation $E$ is
well-founded and extensional, $j$ is an embedding of
$(H_\sigma,\in)$ into $(A,E)$, and
\[
 \begin{gathered}
 \operatorname{pred}_E(j(a))=j[a]\quad(a\in H_\sigma),\\
 \operatorname{pred}_E(h)=j[H_\sigma],\qquad
 \operatorname{pred}_E(s)=j[\sigma].
 \end{gathered}
\]
The element $r$ represents the code $R$, with all its components, in
$(A,E)$; this is expressed by comparing the elements of $R$ and of
the set represented by $r$ one by one.
The relation $\operatorname{Sat}$ satisfies the Tarski conditions for
the satisfaction relation of $(A,E)$, and
\[
 \forall n\in\omega\,
 (\operatorname{Ax}_{\mathcal T}(n)\Rightarrow\operatorname{Sat}(n)),
\]
where $\operatorname{Ax}_{\mathcal T}$ is the recursive set of
G\"odel numbers of the axioms of $\mathcal T$, with the constants
interpreted by $j(\rho)$, $s$ and $r$.

All these conditions are first order over $H_\sigma$. In particular,
well-foundedness is expressed by
\[
 \neg\exists f\in H_\sigma\,
 \bigl(f:\omega\to A\ \land\
       \forall n<\omega\ f(n+1)\mathrel E f(n)\bigr),
\]
since every $\omega$-sequence of elements of $A\subseteq\sigma$
belongs to $H_\sigma$.

Suppose that $\operatorname{Wit}(R,\chi)$ holds, and let $\pi$ be the
transitive collapse of $(A,E)$. By induction on the rank,
\[
 \pi(j(a))=a\quad(a\in H_\sigma),\qquad
 \pi(h)=H_\sigma,\qquad\pi(s)=\sigma,
\]
and $\pi(r)=R$. Thus the collapse is a transitive model of
$\mathcal T$ which contains $H_\sigma$, $\sigma$ and $R$, and
Lemma~\ref{lem:eligibility-prefix} shows that $R$ is the code of a
$\Gamma$-iteration.

Conversely, let $R$ be the code of a $\Gamma$-iteration. Take a
sufficiently large regular $\theta$ and
$N\prec H_\theta$ with $|N|=\sigma$, $H_\sigma\subseteq N$, and
$H_\sigma,\sigma,R\in N$. The transitive collapse of $N$ is the
identity on $H_\sigma$ and therefore maps $R$ to itself.
The collapsed model and its satisfaction relation yield a class
$\chi$ with $\operatorname{Wit}(R,\chi)$.
\end{proof}

The code $R$ and the witness $\chi$ can be combined into one class.
Thus an existential quantifier over $\Gamma$-iterations can be
expressed by an existential quantifier over pairs of classes
$(R,\chi)$ with $\operatorname{Wit}(R,\chi)$, and a universal
quantifier over $\Gamma$-iterations by a universal quantifier over
such pairs. The model coded by $\chi$ is used only to verify that
$R$ codes an iteration and to compute the order of the conditions
enumerated by the last map $e_\delta$ of the code: by
Lemma~\ref{lem:eligibility-prefix}, the order computed in this model
is the true one. Thus the preorder on $\sigma$ induced by $e_\delta$,
which is a dense subset of the last forcing, is first-order definable
over $H_\sigma$ from $R$ and $\chi$, and it does not depend on the
choice of $\chi$. All other class quantifiers in the formulas below
range over all classes of the universe. By
Proposition~\ref{prop:native-full}, the quotients of $D$ are forcing
equivalent to such iterations.

\subsection{Summary of the properties of the forcing}

\begin{theorem}\label{lem:master}\label{thm:master-full}
For each $B=B_{\vec T}$, the forcing $D=D_{\vec T}$ defined in
Section~\ref{sec:master} has the following properties. Its
bookkeeping function handles the names for forcings in the classes
$\Gamma_i$ and for term forcings as stated in
Proposition~\ref{prop:schedule-full}, and its definition is uniform
in $\vec T$.
\begin{enumerate}[label=(M\arabic*)]
\item\label{M:closed} For every $j\geq1$, $D_{\geq j}$ is a complete
subforcing of $D$ and is $<\omega_j$-closed. Moreover, for every
$p\in D$ and every $q\leq p\restriction\bigcup_{k\geq j}I_k$ in
$D_{\geq j}$ there is $p'\leq p$ whose restriction to
$\bigcup_{k\geq j}I_k$ is $q$ and which agrees with $p$ on the other
coordinates.
\item\label{M:native} Let $i\geq1$. In $W_{i+1}$, the quotient
$D_{\geq i}/G_{D_{>i}}$ is forcing equivalent to a
$<\omega_i$-support iteration with the stationary
$\omega_{i+1}$-chain condition. For every $C\in\mathcal C_i$ and
every class $X$ over $H_{\omega_i}$ in $W_i$ there is
$C'\in\mathcal C_i$ with $C\subseteq C'$ and $X\in N_{C'}$. For all
$C\subseteq C'$ in $\mathcal C_i$, the quotient $P_{C'}/G_C$ is
forcing equivalent to a $\Gamma_i$-iteration. For every
$C\in\mathcal C_i$, every forcing $Q\in N_C$ which belongs to
$\Gamma_i$ in $N_C$ and every $q\in Q$, there are
$C\subseteq C_0\subseteq C'$ in $\mathcal C_i$ and a filter $H$
which contains $q$ and is $Q$-generic over $N_{C_0}$, such that
$N_{C'}=N_{C_0}[H]$.
\item\label{M:cores} Every finite family of classes over $\HC$ in
$B[G_D]$ belongs to $B[G_C]$ for some closed set $C\in B$ with
$|C\cap I_i|\leq\omega_{i+1}$ for every $i$.
\item\label{M:cc} $D$ has cardinality at most $\aleph_{\omega+2}$ and
is $\aleph_{\omega+2}$-cc. It preserves all cardinals and adds no reals.
\item\label{M:cert} For every $i\geq1$, the $\Gamma_i$-iterations
have codes which are classes over $H_{\omega_{i+1}}$, and there is a
first-order formula $\operatorname{Wit}(R,\chi)$ over
$H_{\omega_{i+1}}$ such that, in $W_{i+1}$ and in $N_C$ for every
$C\in\mathcal C_i$, a class $R$ is the code of a
$\Gamma_i$-iteration if and only if
$\exists\chi\,\operatorname{Wit}(R,\chi)$. This formula has no
parameters; in particular it mentions neither the sets $T_n$ nor the
clubs added by $C_{\vec T}$.
\item\label{M:term} Let $j\geq1$ and let $C$, $P$ and $\dot E$ be as
in Lemma~\ref{lem:term-bookkeeping-verified}. Then the term forcing
$\Term(P,\mathbb{A}_{\omega_j}(\dot E))$ is appended at a stage after
all coordinates on which $P$ and $\dot E$ depend. The generic
filter induced on $P$ is generic over $W_j$, which contains the
generic filter for the term forcing.
\end{enumerate}
\end{theorem}
\begin{proof}
The results of Section~\ref{sec:master} hold for every iteration
with partial memory, and in particular for $D$. The uniformity in
$\vec T$ was shown after Definition~\ref{def:coding-recursion}.

\ref{M:closed} is Lemma~\ref{lem:native-structural}. The first
assertion of \ref{M:native} is Proposition~\ref{prop:native-full},
and the remaining ones are Lemma~\ref{lem:coretruth}.
\ref{M:cc} is Lemmas~\ref{lem:size-cc-verified}
and~\ref{lem:cardinals-verified}. \ref{M:cert} is
Proposition~\ref{prop:eligibility-exact}, together with the fact
that $(\ast)$ holds for $\rho=\omega_i$ in the models in question.
\ref{M:term} is Proposition~\ref{prop:schedule-full}(2).

For \ref{M:cores}, let $X$ be a class over $\HC$ in $B[G_D]$. In
$W_2=B[G_{D_{>1}}]$, $X$ has a nice name over the quotient at level
one which consists of $\omega_1$ conditions. Code this name by a
sequence of ordinals of length $\omega_1$, and apply
Proposition~\ref{prop:localization-full} to this sequence. Together
with the closure of the supports of the conditions, this gives a
closed set $C\in B$ which satisfies the bounds at every level, such
that $X\in B[G_C]$. For a finite family of classes, take the union of
the corresponding sets.
\end{proof}

\section{Forcing and the reduction of complexity}\label{sec:truthdrop}

\subsection{The class \texorpdfstring{$\Gamma_i$}{Gamma-i}}

Fix $i\geq1$, and write $\rho=\omega_i$ and $\sigma=\rho^+$.
All models $N$ in which the class $\Gamma_i$ is considered in this
section satisfy the hypothesis $(\ast)$ of
Subsection~\ref{subsec:eligibility-justification}, that is,
\[
 \rho^{<\rho}=\rho,\qquad 2^\rho=\sigma.
\]
Each $P\in\Gamma_i$ is given with an iteration code as in
Subsection~\ref{subsec:eligibility-justification}.

\begin{lemma}\label{lem:gamma-closure-verified}
Every forcing in $\Gamma_i$ is $<\rho$-closed, is
$\sigma$-cc, has a dense subset of cardinality at most $\sigma$, and
preserves $(\ast)$. Moreover:
\begin{enumerate}
\item If $P,Q\in\Gamma_i$ in $N$, then in $N[G_P]$ the forcing $Q$
of $N$ densely embeds into the iteration which is given there by the
same names for the iterands. In particular $P\times Q$ is forcing equivalent
to a member of $\Gamma_i$ in $N$.
\item If $P\in\Gamma_i$ and $P\Vdash\dot Q\in\Gamma_i$, then
$P*\dot Q$ is forcing equivalent to a member of $\Gamma_i$.
\item If $P\in\Gamma_i$ and $p\in P$, then $P\restriction p$
is forcing equivalent to a member of $\Gamma_i$.
\end{enumerate}
\end{lemma}
\begin{proof}
By Corollary~\ref{cor:CDMMS}, every $\Gamma_i$-iteration and each of
its initial segments is $<\rho$-closed and $\sigma$-cc.
By the chain condition, the value of a condition at a coordinate is
given by a nice name with antichains of cardinality at most $\rho$.
By induction, a condition which is given in this way is coded by a
set of at most $\rho$ ordinals below $\sigma$. Since
$\sigma^\rho=\sigma$, the iteration has a dense subset of size at
most $\sigma$.

The iteration adds no sequences of ordinals of length less than
$\rho$. It preserves $\sigma$, and a nice name for a subset of
$\rho$ is coded by a set of at most $\rho$ ordinals below $\sigma$;
there are at most $\sigma^\rho=\sigma$ such names. Thus it preserves
both equations in $(\ast)$.

For (1), apply Lemma~\ref{lem:old-iteration-verified} with $R=P$
to the iteration for $Q$. Concatenate $P$ with a name for the
iteration defined in $N[G_P]$ from the same names. The same concatenation proves (2): strengthen a condition
to decide the support of its second component, using the
$<\rho$-distributivity of $P$. If the length of the second iteration
is variable, the $\sigma$-chain condition gives a bound below
$\sigma^+$, and we extend the iteration to that length by trivial
iterands. For (3) use the remark on cones following
Definition~\ref{def:small-iterations}.
\end{proof}

\subsection{The forcing relation for second-order formulas}

We now show that the forcing relation for second-order statements
over $H_\tau$ is definable by a second-order formula of the same
complexity.

\begin{lemma}\label{lem:forcing-prefix-verified}
Let $\tau>\omega$ be regular with $\tau^{<\tau}=\tau$ in $N$, and let
$Q\in N$ be $\tau$-cc with a dense subset which is given as a preorder
on $\tau$. For every second-order formula $\varphi$ over
$H_\tau^{N[G_Q]}$, the relation $q\Vdash_Q\varphi$ is expressed over
$H_\tau^N$ by a second-order formula with the same prefix of class
quantifiers, which is obtained recursively from $\varphi$. In
particular, the forcing relation of a first-order formula is first
order, the forcing relation of a $\Sigma^2_m$ formula is
$\Sigma^2_m$, and the forcing relation of a $\Pi^2_m$ formula is
$\Pi^2_m$. The parameters of the translated formula are classes
coding the dense preorder and the names for the class parameters of
$\varphi$.
\end{lemma}
\begin{proof}
Put $H=H_\tau^N$, and let $Q_0$ be the given dense preorder on
$\tau$.

\emph{Names for elements.} Every element $x$ of $H_\tau^{N[G_Q]}$ is
the image of a point under the transitive collapse of a well-founded
extensional relation on some $\gamma<\tau$. We code a name for $x$ by
a triple $a$ consisting of nice $Q_0$-names for the field and for the
relation, and a maximal antichain deciding the point. Such a name
involves fewer than $\tau$ statements, each of which is decided by an
antichain of size $<\tau$; hence $a\in H$. Every
$x\in H_\tau^{N[G_Q]}$ has a name of this form: by the chain
condition, $\gamma$ can be bounded below a given condition $p$, and
mixing with a name for the empty set below the conditions
incompatible with $p$ gives a name which is forced by the trivial
condition to be of the required form.

Let $\operatorname{Obj}(a)$ express that $a$ is such a code, and write
$x_a$ for the element named by $a$. The predicate
$\operatorname{Obj}$ is first-order definable over $H$. To see this,
note that a function $f:\omega\to\gamma$ in the extension has a nice
name for its graph which is an element of $H$, and that every
descending $\omega$-sequence witnessing ill-foundedness has such a
name; these names are used to express well-foundedness, and they do
not involve the codes of elements just defined. The Boolean values
$\|x_a=x_b\|$ and $\|x_a\in x_b\|$ are expressed by means of names for
isomorphisms between the corresponding relations below their
distinguished points: by the maximum principle such an isomorphism
has a name consisting of fewer than $\tau$ pairs, and conversely an
isomorphism shows that the collapses agree. Hence the Boolean values
of atomic formulas are first-order definable over $H$.

\emph{Boolean values.} We represent Boolean values by antichains
$A\in H$ of $Q_0$, modulo equality of the values $[A]$ which they
represent. Their order is first order:
\[
 [A]\leq[B]\quad\Longleftrightarrow\quad
 \forall a\in A\,\forall q\leq a\,
 \exists r\leq q\,\exists b\in B\ (r\leq b).
\]
Consequently equality and complements are first order. The
supremum of a family of Boolean values which is given by a class is
its least upper bound in this order; here all quantifiers range over
antichains in $H$. The supremum exists by the chain condition and the
completeness of the Boolean completion.

\emph{Names for classes.} A $Q$-name for a class over
$H_\tau^{N[G_Q]}$ is coded by a function $\chi$ which assigns a
Boolean value to every $a$ with $\operatorname{Obj}(a)$, such that
\[
 \|x_a=x_b\|\wedge\chi(a)
   =\|x_a=x_b\|\wedge\chi(b)
 \quad(\operatorname{Obj}(a),\operatorname{Obj}(b)).
\]
Let $\operatorname{Name}(\chi)$ express that $\chi$ is such a function,
where the values are determined up to equality of Boolean values;
this is first order. Every name in $N$ for a class over
$H_\tau^{N[G_Q]}$ gives rise to such a class $\chi$ over $H$.
Conversely, every $\chi$ with $\operatorname{Name}(\chi)$ determines
a name for a class, and the value of this name does not depend on the
choice of the codes $a$.

\emph{First-order formulas.} For a first-order formula $\theta$, the
Boolean value $\|\theta\|$ is defined by recursion on $\theta$, using
the Boolean operations and, for the quantifiers, suprema and infima
over all $a$ with $\operatorname{Obj}(a)$. Since $\theta$ is a fixed
formula, this gives a first-order definition of $\|\theta\|$, and
$p\Vdash\theta$ means $[p]\leq\|\theta\|$.

\emph{Class quantifiers.} Finally, the maximum principle and mixing
with a name for the empty class below the conditions incompatible
with $p$ give
\[
 \begin{aligned}
 p\Vdash\exists X\,\psi(X)
 &\ \Longleftrightarrow\
 \exists\chi\,[\operatorname{Name}(\chi)\land p\Vdash\psi(\chi)],\\
 p\Vdash\forall X\,\psi(X)
 &\ \Longleftrightarrow\
 \forall\chi\,[\operatorname{Name}(\chi)\Rightarrow p\Vdash\psi(\chi)].
 \end{aligned}
\]
Here $\chi$ ranges over all classes over $H$ in $N$. By recursion on
the finite prefix of class quantifiers we obtain a formula with the
same prefix.
\end{proof}

We apply this lemma with $\tau=\sigma$. The hypotheses
hold because $2^\rho=\sigma$ implies $\sigma^{<\sigma}=\sigma$, and
because every forcing in $\Gamma_i$ is $\sigma$-cc and has a dense
subset of size at most $\sigma$.
Note that a name for a subset of $H_\rho$ is coded by an
\emph{element} of $H_\sigma$, since it has hereditary size at most
$\rho$. By contrast, a name for an arbitrary class over
$H_\sigma$ of the extension is coded by a \emph{class} over
$H_\sigma$. This difference is the reason for the drop in complexity
in Lemma~\ref{lem:drop}.

\subsection{The translation \texorpdfstring{$T_i$}{Ti}}

We consider second-order formulas in prenex form. Using a fixed
pairing function, we contract every block of class quantifiers of the
same kind to a single class quantifier. The translation $T_i$ is
defined by recursion on the number of class quantifiers. The base
case consists of the formulas with at most one class quantifier, that
is, of the first-order formulas (with class parameters) and the
$\Sigma^2_1(H_\rho)$ and $\Pi^2_1(H_\rho)$ formulas.

\begin{definition}
Let $N$ be a model of $(\ast)$, and let the class parameters belong to
$N$.
\begin{enumerate}
\item If $\varphi$ is a first-order, $\Sigma^2_1(H_\rho)$ or
$\Pi^2_1(H_\rho)$ formula, let
$T_i(\varphi)$ be its first-order translation over $H_\sigma$:
relativize the first-order quantifiers to the definable set $H_\rho$,
and replace the class quantifier over $H_\rho$, if there is one, by a
first-order quantifier over the subsets of $H_\rho$, all of which
belong to $H_\sigma$. The class parameters of $\varphi$ are elements
of $H_\sigma$ and are used as such.
\item For formulas with more than one class quantifier put
\[
 T_i(\neg\varphi)=\neg T_i(\varphi).
\]
\item If $\exists X\,\varphi$ has at least two class quantifiers, let
$T_i(\exists X\,\varphi)$ assert that there are an iteration code
with last forcing $Q\in\Gamma_i$, a condition $q\in Q$, and a
$Q$-name $\dot X$ for a subset of $H_\rho$ such that
\[
                    q\Vdash_Q T_i(\varphi)(\dot X).
\]
\end{enumerate}
A universal class quantifier is treated as $\neg\exists\neg$. In (3),
the forcing relation is the usual forcing relation for the formula
$T_i(\varphi)$, which has already been defined.
\end{definition}

\begin{lemma}\label{lem:potential-invariance-verified}
Let $N$ be a model of $(\ast)$, let $P$ be forcing equivalent in
$N$ to a member of $\Gamma_i$, and let $G_P$ be $P$-generic over $N$.
Then for every formula
$\varphi$ as above with parameters in $N$,
\[
 N\models T_i(\varphi)
 \quad\Longleftrightarrow\quad N[G_P]\models T_i(\varphi).
\]
\end{lemma}
\begin{proof}
Since $P$ is forcing equivalent to a member of $\Gamma_i$, we may
assume $P\in\Gamma_i$. We argue by induction on the number of class
quantifiers of $\varphi$. In the base case, $T_i(\varphi)$ is
equivalent to a first-order, $\Sigma^2_1$ or $\Pi^2_1$ statement over $H_\rho$; since
$P$ is $<\rho$-closed, $H_\rho$ is the same in $N$ and $N[G_P]$, and
Fact~\ref{fact:closedabs} gives the equivalence. The step for
negation is immediate.

Consider $\exists X\,\varphi$ with at least two class quantifiers,
and put $\psi=T_i(\varphi)$. Suppose first that $Q,q,\dot X\in N$
witness $T_i(\exists X\,\varphi)$ in $N$. In $N[G_P]$ let
$e:Q\to\widehat Q$ be the dense embedding of
Lemma~\ref{lem:gamma-closure-verified}(1) into the iteration defined
in $N[G_P]$ from the same names. For any
$\widehat Q$-generic filter over $N[G_P]$ containing $e(q)$, let
$G_Q$ be the induced $Q$-generic filter and put $X=\dot X^{G_Q}$.
The filters $G_P,G_Q$ are product-generic over $N$. Thus
$N[G_Q]\models\psi(X)$, and $P$ is forcing equivalent in $N[G_Q]$
to a member of $\Gamma_i$ by the same lemma. By induction
$N[G_Q][G_P]\models\psi(X)$. Hence $\widehat Q$, $e(q)$ and the
$\widehat Q$-name for $X$ which corresponds to $\dot X$ under $e$
witness $T_i(\exists X\,\varphi)$ in $N[G_P]$.

Conversely, suppose that $T_i(\exists X\,\varphi)$ holds in
$N[G_P]$. Let $\dot Q,\dot q,\ddot X\in N$ be names for witnesses, and
choose $p\in G_P$ which forces that they are as required.
By Lemma~\ref{lem:gamma-closure-verified}, choose $S\in\Gamma_i$
forcing equivalent to $(P\restriction p)*(\dot Q\restriction\dot q)$.
Under this equivalence, $\ddot X$ corresponds to an $S$-name
$\dot X^*$ with $1_S\Vdash_S\psi(\dot X^*)$. Thus $S,1_S,\dot X^*$ witness
$T_i(\exists X\,\varphi)$ in $N$.
\end{proof}

\begin{lemma}\label{lem:truthlemma}
Let $C\in\mathcal C_i$. Then for every second-order formula $\varphi$
over $H_\rho$ and all parameters $\vec a\in N_C$,
\[
 W_i\models\varphi(\vec a)
 \quad\Longleftrightarrow\quad
 N_C\models T_i(\varphi)(\vec a).
\]
\end{lemma}
\begin{proof}
Recall that $W_i=B[G_{D_{\geq i}}]$ and
$N_C=W_{i+1}[G_C]\subseteq W_i$. We use the following facts.
\begin{enumerate}[label=(\roman*)]
\item $W_i$ is an extension of $N_C$ by a $<\rho$-closed forcing
(Lemma~\ref{lem:quotient-closed}).
\item For every $X\subseteq H_\rho$ in $W_i$ there is
$C'\in\mathcal C_i$ with $C\subseteq C'$ and $X\in N_{C'}$, and
the quotient from $N_C$ to $N_{C'}$ is forcing equivalent to
a member of $\Gamma_i$ (Lemma~\ref{lem:coretruth}(a)).
\item For every forcing $Q\in N_C$ which belongs to $\Gamma_i$ in
$N_C$ and every $q\in Q$ there are $C\subseteq C_0\subseteq C'$ in
$\mathcal C_i$ and a filter $H$ with $q\in H$ which is
$Q$-generic over $N_{C_0}$, such that $N_{C'}=N_{C_0}[H]$
(Lemma~\ref{lem:coretruth}(b)).
\end{enumerate}
We argue by induction on the number of class quantifiers of
$\varphi$, simultaneously for all $C\in\mathcal C_i$ and all
parameters. The base case follows from (i) and
Fact~\ref{fact:closedabs}, and the step for negation is immediate.

Suppose that $W_i\models\exists X\varphi(X,\vec a)$. By (ii), choose
$C'$ such that a witness $X$ belongs to $N_{C'}$. The induction
hypothesis gives
$N_{C'}\models T_i(\varphi)(X,\vec a)$. Use (ii) to choose
$R\in\Gamma_i$ representing the quotient, with $N_{C'}=N_C[G_R]$.
By the forcing theorem there
are $q\in G_R$ and an $R$-name $\dot X$ with
$q\Vdash_R T_i(\varphi)(\dot X,\vec a)$.
Thus $(R,q,\dot X)$ witnesses
$N_C\models T_i(\exists X\varphi)(\vec a)$.

Conversely, let $Q,q,\dot X\in N_C$ witness
$T_i(\exists X\varphi)(\vec a)$ in $N_C$, and choose $C_0$, $C'$ and
$H$ as in (iii). By (ii), applied to $C\subseteq C_0$, there are
$S\in\Gamma_i$ in $N_C$ and an $S$-generic filter $g$ over $N_C$ with
$N_{C_0}=N_C[g]$. Since $H$ is $Q$-generic over $N_C[g]$, the pair
$(g,H)$ is generic for $S\times Q$ over $N_C$; in particular $H$ is
$Q$-generic over $N_C$, and $g$ is $S$-generic over $N_C[H]$. Put
$X=\dot X^H$. Since $q\in H$,
\[
 N_C[H]\models T_i(\varphi)(X,\vec a).
\]
By Lemma~\ref{lem:gamma-closure-verified}(1), $S$ is forcing
equivalent in $N_C[H]$ to a member of $\Gamma_i$. Hence
Lemma~\ref{lem:potential-invariance-verified}, applied to the model
$N_C[H]$ and the formula $\varphi$, gives
\[
 N_{C'}=N_C[H][g]\models T_i(\varphi)(X,\vec a),
 \qquad\text{and hence}\qquad
 W_i\models\varphi(X,\vec a)
\]
by the induction hypothesis for $C'$.
Hence $W_i\models\exists X\varphi(X,\vec a)$.
\end{proof}

\begin{lemma}\label{lem:drop}
For $k\geq2$ the translation can be defined uniformly so that
\[
\begin{aligned}
 T_i(\Sigma^2_k(H_{\omega_i}))
       &\subseteq\Sigma^2_{k-1}(H_{\omega_{i+1}}),\\
 T_i(\Pi^2_k(H_{\omega_i}))
       &\subseteq\Pi^2_{k-1}(H_{\omega_{i+1}}).
\end{aligned}
\]
The translated formulas have no parameters for $\vec T$, for the
clubs added by $C_{\vec T}$, or for the bookkeeping function.
\end{lemma}
\begin{proof}
Since $\mathcal P(H_\rho)\subseteq H_\sigma$ and $H_\rho,\rho$
are definable in $H_\sigma$, the first-order formulas and the
$\Sigma^2_1(H_\rho)$ and $\Pi^2_1(H_\rho)$ formulas are translated
into first-order formulas over
$H_\sigma$, and their class parameters are elements of $H_\sigma$.

For $\exists X\varphi\in\Sigma^2_k(H_\rho)$, we express
$T_i(\exists X\varphi)$ by means of a single existential class
quantifier $\exists Z$, where $Z$ codes four classes
\[
 Z_0=R,\qquad Z_1=\chi,\qquad
 Z_2=\{q\},\qquad Z_3=\{a\},
\]
by a fixed pairing function. Here $R,\chi\subseteq H_\sigma$ are a
code of $Q$ and a witness for it, $q,a\in H_\sigma$, and $a$ codes
the name $\dot X$.
By Proposition~\ref{prop:eligibility-exact}, the condition
$\operatorname{Wit}(R,\chi)$ is first order. The conditions on $q$
and $a$ are also first order. The elements $q$ and $a$ are uniquely
determined by $Z$, independently of the classes quantified later, so
the first-order quantifiers which define them from $Z$ can be moved
inside the remaining class quantifiers.

Lemma~\ref{lem:forcing-prefix-verified} requires a dense subset of
$Q$ which is given as a preorder on $\sigma$. We use the preorder
induced by the last map $e_\delta$ of the code $R$. Its order is not
a component of $R$, but it is computed correctly in the model coded
by $\chi$ (Lemma~\ref{lem:eligibility-prefix}); hence it is
first-order definable over $H_\sigma$ from $R$ and $\chi$, as noted
after Proposition~\ref{prop:eligibility-exact}. By the induction
hypothesis and Lemma~\ref{lem:forcing-prefix-verified},
the statement $q\Vdash_Q T_i(\varphi)(\dot X)$ is therefore first order over
$H_\sigma$ for $k=2$, and $\Pi^2_{k-2}(H_\sigma)$ for $k>2$.
Together with the leading quantifier $\exists Z$ this gives a
$\Sigma^2_{k-1}(H_\sigma)$ formula in both cases. Taking negations
gives the assertion for $\Pi^2_k$.

All these formulas refer to the class $\Gamma_i$ of the model in
which they are evaluated. The translation is defined by recursion
on the formula, using only the definable cardinals $\rho,\sigma$ and
the class $\Gamma_i$; hence it is recursive and introduces no
parameters.
\end{proof}

\subsection{Reduction to \texorpdfstring{$B_{\vec T}$}{B-T}}

\begin{proposition}\label{prop:compile}
For every $\Pi^2_n$ formula $A$ on the reals there is a
$\Pi^2_1(H_{\omega_n})$ formula $\Phi_A$, obtained recursively from
$A$ and independent of $\vec T$, such that
\begin{equation}\label{eq:compile}
 B_{\vec T}[G_{D_{\vec T}}]\models A(a)
 \quad\Longleftrightarrow\quad
 B_{\vec T}\models\Phi_A(a)
\end{equation}
for all reals $a$. The corresponding statement holds for
$\Sigma^2_n$ formulas, with
$\Sigma^2_1(H_{\omega_n})$ in place of $\Pi^2_1(H_{\omega_n})$.
The only parameters of $\Phi_A(a)$ are the real parameters $a$.
\end{proposition}
\begin{proof}
Apply Lemmas~\ref{lem:truthlemma} and~\ref{lem:drop} successively for
$i=1,\ldots,n-1$, each time with $C=\varnothing$. At step $i$ we
pass from $W_i=B_{\vec T}[G_{D_{\geq i}}]$ to
$N_\varnothing=W_{i+1}=B_{\vec T}[G_{D_{\geq i+1}}]$. This model
satisfies $(\ast)$ for $\rho=\omega_i$, because $B_{\vec T}$
satisfies GCH and $D_{\geq i+1}$ is
$<\omega_{i+1}$-closed. At each step the complexity drops by one
and the domain of the first-order quantifiers increases from
$H_{\omega_i}$ to $H_{\omega_{i+1}}$.

The resulting statement is $\Pi^2_1(H_{\omega_n})$, and it is
evaluated in $W_n=B_{\vec T}[G_{D_{\geq n}}]$. The forcing
$D_{\geq n}$ is $<\omega_n$-closed, so
Fact~\ref{fact:closedabs} shows that the statement holds in $W_n$ if
and only if it holds in $B_{\vec T}$. For $n=1$ only this last step
is needed. The same argument applies to $\Sigma^2_n$ formulas.
The class $\Gamma_i$ and the formula $\operatorname{Wit}$ do not
depend on $\vec T$, and the translation is defined by recursion on
formulas. Hence $\Phi_A$ is obtained recursively from $A$ and has no
parameters other than $a$.
\end{proof}

\section{The formulas \texorpdfstring{$\Code_n$}{Code-n}}\label{sec:decoders}

\subsection{The case \texorpdfstring{$n=1$}{n=1}}

Put
\[
 \Code_1(t)\iff\exists C\subseteq\omega_1\,[C\text{ is club and }C\cap S^1_t=\varnothing].
\]
By Proposition~\ref{hyp:source}\ref{S:graph}, $\Code_1$ is $\Sig_1$
over $\HC$.

\subsection{Iterated almost disjoint coding}

Recall that $\langle B^r_\eta:\eta<\omega_{r+1}\rangle$ is the
canonical constructible almost disjoint family of subsets of
$\omega_r$. For $X\subseteq\omega_1$ write $Z_1(X)=X$ and
\[
 Z_{r+1}(X)=\{\eta<\omega_{r+1}:
       Z_r(X)\cap B^r_\eta\text{ is bounded in }\omega_r\}.
\]
The set $Z_r(X)$ belongs to every model considered below which
contains $X$, and its definition is absolute between these models.

\begin{lemma}\label{lem:generic-pair-full}
Let $V\subseteq V'$ be transitive models of $\mathsf{ZFC}$, let
$P\in V$ be a forcing, and put
$T=\Term^V(P,\mathbb{A}_\rho(\dot E))$ for a $P$-name $\dot E\in V$,
where $\rho$, $P$ and $\dot E$ satisfy the hypotheses of
Subsection~\ref{sec:term} in $V$.
Suppose $H\in V'$ is $T$-generic over $V$ and $G_P$ is $P$-generic
over $V'$. Then $(G_P,H)$ is $P\times T$-generic over $V$, and
$V'[G_P]$ contains a subset of $\rho$ which codes
$\dot E^{G_P}$ by almost disjoint coding.
\end{lemma}
\begin{proof}
Since $V[H]\subseteq V'$, the filter $G_P$ is $P$-generic over $V[H]$.
The two-step iteration $T*\check P$ is the product $T\times P$.
Apply the projection
$P\times T\longrightarrow P*\mathbb{A}_\rho(\dot E)$
from Theorem~\ref{thm:term}. The generic filter induced by this
projection gives the required set. Whether the intersection of
two given subsets of an ordinal is bounded is absolute between
transitive models.
\end{proof}

\begin{lemma}\label{lem:cascade}
Fix $j\geq1$ and work in
$W_j=B[G_{D_{\geq j}}]$.
Let $C$ be a set of coordinates of levels below $j$, with at most
$\omega_{k+1}$ coordinates of level $k$ for each $k<j$, such that
$C\cup\bigcup_{k\geq j}I_k$ is closed. Let $P$ be the
preorder on $\omega_j$ representing the quotient
$D_{C\cup\bigcup_{k\geq j}I_k}/G_{D_{\geq j}}$ as in
Lemma~\ref{lem:term-bookkeeping-verified}. Then $P$ is
$\omega_j$-cc. For $j=1$ it is forcing equivalent to the
trivial forcing.

If $\dot E$ is a nice $P$-name for a subset of $\omega_{j+1}$, then
the term forcing
$\Term(P,\mathbb{A}_{\omega_j}(\dot E))$ is appended in the
construction of $D$ at a coordinate of level $j$, after all
coordinates on which $P$ and $\dot E$ depend. Its generic filter and
the $P$-generic filter induced by $G_D$ give a generic filter for
$P*\mathbb{A}_{\omega_j}(\dot E)$.
\end{lemma}
\begin{proof}
The quotient is a finite iteration of the restrictions to $C$
of the quotients at levels $j-1,\ldots,1$. Each is $\omega_j$-cc
by Proposition~\ref{prop:native-full}, so $P$ is $\omega_j$-cc.

By Lemma~\ref{lem:term-bookkeeping-verified}, $P$ and $\dot E$ belong
to $B[G_H]$ for a suitable closed set $H$ of coordinates of levels at
least $j$. In this model
$\omega_j^{<\omega_j}=\omega_j$ and
$\omega_{j+1}^{<\omega_j}=\omega_{j+1}$ hold, since
neither the forcing $D_{>j}$ nor an initial segment of the quotient at
level $j$ adds sequences of ordinals of length $<\omega_j$.
Theorem~\ref{thm:term} therefore shows that the term forcing, defined
in $B[G_H]$, is $<\omega_j$-compact and $\omega_j$-linked and has a
dense subset of size at most $\omega_{j+1}$.

By Proposition~\ref{prop:schedule-full}(2), this term forcing is
appended at some stage. Let $H'$ be its generic filter.
The filter induced on $P$ is $P$-generic over $W_j$, and $W_j$
contains $H'$; here we use that the restriction to the closed set
$C\cup\bigcup_{k\geq j}I_k$, which contains all coordinates of levels
at least $j$, is a projection. By
Lemma~\ref{lem:generic-pair-full}, the two filters are generic for the
product over the model in which the term forcing was defined, and the
projection of Theorem~\ref{thm:term} gives a generic filter for
$P*\mathbb A_{\omega_j}(\dot E)$.
\end{proof}

\begin{proposition}\label{prop:cascade-full}
Fix $j\geq1$, put $W_j=B[G_{D_{\geq j}}]$, and let
\[
 C\subseteq\bigcup_{1\leq k<j}I_k,\qquad
 |C\cap I_k|\leq\omega_{k+1}\quad(1\leq k<j),
\]
be such that $C\cup\bigcup_{k\geq j}I_k$ is closed. Put
\[
 Q_j=D_{C\cup\bigcup_{k\geq j}I_k}/G_{D_{\geq j}}.
\]
Then for every $E\subseteq\omega_{j+1}$ in $W_j[G_{Q_j}]$, the final
model $B[G_D]$ contains a set $X\subseteq\omega_1$ with
$Z_{j+1}(X)=E$.
\end{proposition}
\begin{proof}
For $1\leq r\leq j$, put
$C_r=C\cap\bigcup_{1\leq k<r}I_k$, and define $W_r$ and $Q_r$ as
above, with $C_r$ and $r$ in place of $C$ and $j$. Since the memory
of a coordinate contains only coordinates of the same or higher
levels, $C_r\cup\bigcup_{k\geq r}I_k$ is closed. By
Lemma~\ref{lem:cascade}, we can identify $Q_r$ with a preorder
$P_r$ on $\omega_r$, by means of an enumeration of its conditions,
such that
\[
 |P_r|\leq\omega_r,\qquad P_r\text{ is }\omega_r\text{-cc};
\]
the forcing $Q_1$ is trivial.
Starting with $E_{j+1}=E$, we construct by downward recursion sets
$E_r\in W_r[G_{Q_r}]$, $E_r\subseteq\omega_r$, for $r=j,\ldots,1$,
such that
\[
 E_{r+1}=\{\eta<\omega_{r+1}:
 E_r\cap B^r_\eta\text{ is bounded in }\omega_r\}.
\]

At step $r$, a nice $P_r$-name $\dot E_{r+1}$ for $E_{r+1}$ is coded
by at most $\omega_{r+1}$ ordinals, since its antichains
have size $<\omega_r$. By
Lemma~\ref{lem:term-bookkeeping-verified}, $P_r$ and $\dot E_{r+1}$
belong to $B[G_{H_r}]\subseteq W_r$ for a closed set $H_r\in B$ of
coordinates of levels at least $r$, with the same order on $P_r$, and
the term forcing
$T_r=\Term^{B[G_{H_r}]}(P_r,\mathbb{A}_{\omega_r}(\dot E_{r+1}))$
is appended at a later stage, at level $r$.
Since the restriction to $C_r\cup\bigcup_{k\geq r}I_k$ is a
projection, the filter induced on $Q_r$ is generic over $W_r$. It
yields a $P_r$-generic filter over $W_r$, and $W_r$ contains the
generic filter for $T_r$.
Lemma~\ref{lem:generic-pair-full} gives the required set $E_r$.
Note that the enumeration which identifies $Q_r$ with $P_r$ need not
belong to $B[G_{H_r}]$.

For $r>1$ we have
\[
 W_r[G_{Q_r}]\subseteq W_{r-1}[G_{Q_{r-1}}],
\]
because $C_r\subseteq C_{r-1}\cup I_{r-1}$; in particular $W_{r-1}$
contains the generic filter on the coordinates of level $r-1$ in $C$
and the generic filters for the term forcings used so far. Thus the
same procedure applies to $E_r$ and $C_{r-1}$.
At each step we choose a stage of the bookkeeping which is later than
all stages used before.
For $r=1$ the forcing $Q_1$ is trivial, so $T_1$ is the almost
disjoint coding forcing itself. With $X=E_1$, the recursion gives
$Z_r(X)=E_r$ for $1\leq r\leq j+1$, hence $Z_{j+1}(X)=E$.
\end{proof}

\begin{lemma}\label{lem:final-source-stationarity}
If $t\notin T_n$, the set $S^n_t$ remains stationary in
$B_{\vec T}[G_{D_{\vec T}}]$.
\end{lemma}
\begin{proof}
The set $S^n_t$ is stationary in $B_{\vec T}$ by
Lemma~\ref{lem:source}.
First force with $D_{\geq n}$, which is $<\omega_n$-closed.
The statement that there is a club disjoint from $S^n_t$ is
$\Sig_1(H_{\omega_n})$, with $S^n_t$ as a class parameter.
Fact~\ref{fact:closedabs} shows that no such club is added.

The remaining quotients at the finitely many levels below $n$ form an
$\omega_n$-cc iteration. Such a forcing preserves stationary
subsets of $\omega_n$: given a name for
a club, for each $\alpha<\omega_n$ take a maximal antichain
deciding the first point of this club above $\alpha$. There are fewer
than $\omega_n$ possible values; let $f(\alpha)<\omega_n$ be a bound.
The set of limit ordinals closed under $f$ is a club in the ground
model which is forced to be contained in the given club, and it meets
$S^n_t$. For $n=1$ there is no quotient below level $n$.
\end{proof}

\subsection{Codes of models}

In this subsection we use only the following properties of the final
model $B[G_D]$.
\begin{enumerate}[label=(C\arabic*)]
\item All cardinals are preserved, and for every finite family of
classes over $\HC$ there is a closed set $C\in B$ with
$|C\cap I_k|\leq\omega_{k+1}$ for every $k$ such that the family
belongs to $B[G_C]$.
\item Let $r\geq2$, and let $E\subseteq\omega_r$ belong to a model
$W_{r-1}[G_{Q_{r-1}}]$, where $Q_{r-1}$ is as in
Proposition~\ref{prop:cascade-full} for $j=r-1$. Then there is a
class $X\subseteq\omega_1$ with $Z_r(X)=E$.
\item Each stationary set $S^n_t$ with $t\notin T_n$ is stationary
in the final model. For $t\in T_n$ the final model contains
a club $C\subseteq\omega_n$ disjoint from $S^n_t$ and the set
$P^n_t(C)$ of Lemma~\ref{lem:direct-padding}; both belong to
$B_{\vec T}$.
\end{enumerate}
Here (C1) is Theorem~\ref{lem:master}\ref{M:cores} and \ref{M:cc},
(C2) is Proposition~\ref{prop:cascade-full}, and (C3) is
Lemma~\ref{lem:final-source-stationarity} together with the definition
of $C_{\vec T}$. We refer to the models $W_{r-1}[G_{Q_{r-1}}]$ in
(C2) as the \emph{models of (C2) for $r$}. Note that if $C$ is as in
(C1) and $j\geq1$, then $C\cap\bigcup_{k<j}I_k$ satisfies the
hypotheses of Proposition~\ref{prop:cascade-full}, and $B[G_C]$ is
contained in the corresponding model $W_j[G_{Q_j}]$. All these models
have the same cardinals as $L$.

We now fix the coding of models by sets of ordinals.
Fix once and for all constructible functions which code finite
tuples of ordinals, finite sequences, and terms of the language
introduced below by ordinals. For every infinite cardinal $\tau$,
these functions map objects built from ordinals below $\tau$ to
ordinals below $\tau$.

The language $\mathcal L$ is defined in $\omega$ steps. Let
$\mathcal L_0$ consist of $\in$, of two function symbols $F_j,H_j$
for each $j<r$, and of the constants
\[
 b_0,\ldots,b_{r+2},\quad z_1,\ldots,z_k,\quad
 t_1,\ldots,t_l,\quad d_1,\ldots,d_m.
\]
Given $\mathcal L_m$, let $\mathcal L_{m+1}$ be obtained by adding a
function symbol $f_\varphi$ for every formula of $\mathcal L_m$ of
the form $\exists v\,\varphi(v,\vec u)$, with one argument for each
of the free variables $\vec u$. Put
$\mathcal L=\bigcup_{m<\omega}\mathcal L_m$. This language is
countable and recursive, and every formula
$\exists v\,\varphi(v,\vec u)$ of $\mathcal L$ has its function
symbol $f_\varphi$ in $\mathcal L$; we call the functions which
interpret these symbols Skolem functions. In addition we use the
constants
\[
 c_\alpha,\ g_\alpha\qquad(\alpha<\omega_r).
\]
Here the constants $c_\alpha$ are intended to denote the ordinals
$\alpha<\omega_r$, the constants $g_\alpha$ name all elements of the
model, the $b_j$ denote the cardinals $\omega_j$, the $z_i$ denote
subsets of $\omega_r$ which are of the form $Z_r(X)$, and the $t_i$
denote real parameters. The
$d_i$ are finitely many further constants, for objects such as
the order types and isomorphisms occurring in a witness below.

A \emph{code} of a model for this language, in which every element is
the value of a closed term, is a set of ordinals below $\omega_r$
which consists of three disjoint parts, coding the equality relation
between closed terms, the membership relation between closed terms,
and the satisfaction relation of the model. An instance of the
satisfaction relation is given by the G\"odel number of a formula of
the countable language and a finite sequence of closed terms which
are assigned to its free variables; thus the uncountably many
constants occur only in these assignments, and not in the recursive
enumeration of formulas. We write
$q_{\rm eq}(u,v)$, $q_{\rm mem}(u,v)$ and
$q_{\rm sat}(e,\vec u)$ for the ordinals which code the corresponding
instances of equality, membership and satisfaction. These functions,
and the function assigning to a function symbol and codes of terms the
code of the resulting term, are fixed in advance; in particular
they do not depend on the model which is coded.
Similarly, a \emph{code of an ordinal} consists of two disjoint
parts, coding the field and the strict order of a wellorder; the
field is coded separately so that the ordinal $0$ is included.

For $0<a<\omega_{j+1}$ let $f_{j,a}$ be the $<_L$-least
surjection from $\omega_j$ onto $a$. The $\Sigma_1$ definition of the
initial segments of $<_L$ in Section~\ref{sec:source-spec} gives a
functional
$\Sigma_1$ definition of the map $a\mapsto f_{j,a}$, uniformly in $j$:
one asserts the existence of a surjection $f$ and of the set of all
its $<_L$-predecessors, such that no $<_L$-predecessor of $f$ is a
surjection with the same domain and range. The fixed coding
functions have $\Sigma_1$ definitions of the same kind. The theories
below require these $\Sigma_1$ definitions to define total functions
over $L_\beta$ on the relevant domains. By the upward absoluteness
of $\Sigma_1$ formulas, the values computed over $L_\beta$ are then
the true values.

\begin{definition}
For each $r\geq1$ let $\mathsf T_r$ be the theory in the above
language which consists of $\mathsf{ZFC}^-$ and the following axioms.
\begin{enumerate}[label=(\roman*)]
\item $b_j=\omega_j$ for $0\leq j\leq r+2$, and
$b_{r+2}=\omega_{r+2}^L$. We write $\rho=b_r$ and $\beta=b_{r+1}$.
The set $L_\beta$ exists.
\item The definitions $\Phi_D$ and $\Phi_B$ of the canonical diamond
sequence and of the canonical almost disjoint family at $\rho$
are total and functional over $L_\beta$ on exactly the domains
specified in Definition~\ref{def:suitable}. The same holds for the
definitions of the fixed coding functions and their inverses, and
their restrictions to $\beta$ satisfy the corresponding identities for
the even and odd parts of subsets of $\beta$. For each $j\leq r$, the
definition of the map $a\mapsto f_{j,a}$, evaluated over the same
$L_\beta$, is total and functional for $0<a<b_{j+1}$.
\item Every $z_i$ is a subset of $\rho$, and every $t_i$ is a
subset of $\omega$.
\item For $j<r$ and $0<a<b_{j+1}$, $F_j(a)$ is a function from $b_j$
onto $a$, and, for $u\in a$,
\[
 H_j(a,u)\in b_j,\qquad F_j(a)(H_j(a,u))=u.
\]
\end{enumerate}
The axiom schemes of $\mathsf{ZFC}^-$ are taken for the formulas of
the language of set theory only. The finitely many axioms listed here
and the axiom schemes are given by a recursive set of G\"odel
numbers. The Skolem functions are additional functions on the model;
they need not be elements of the model or definable in it, and they
are subject only to the conditions on the satisfaction relation
stated in the proof of Lemma~\ref{lem:recursive-V-explicit}.
\end{definition}

Note that a transitive model $N$ of $\mathsf T_r$ with
$b_r^N=\omega_r$ is suitable for $r+1$ in the sense of
Definition~\ref{def:suitable}.

We shall define formulas $V_r(Y;\vec X,\vec t)$ expressing that the
set decoded from $Y$ is a code of a well-founded model of
$\mathsf T_r$ with certain properties, and we shall prove two things:
every $Y$ satisfying $V_r$ yields such a model (soundness), and every
such model which belongs to a model of (C2), and has cardinality
$\omega_r$ there, is obtained in this way (completeness). The models used for completeness are transitive
models of $\mathsf T_r$ of cardinality $\omega_r$ in which
$b_j=\omega_j$ holds for all $j\leq r$; they are obtained as
transitive collapses of the elementary submodels considered below. In
the soundness direction, the fact that the constants $b_j$, $j\leq r$,
are interpreted by the true cardinals $\omega_j$ will follow from
the formula $V_r$.

\begin{lemma}\label{lem:direct-V1-explicit}
There is a first-order formula
$V_1(Y;X_1,\ldots,X_k,t_1,\ldots,t_l)$ over $\HC$ with the following
properties. If $V_1(Y;\vec X,\vec t)$ holds, then $Y\subseteq\omega_1$
is a code of a well-founded extensional model, together with its
satisfaction relation and its Skolem functions, whose transitive
collapse $N$ satisfies $\mathsf T_1$, $b_1^N=\omega_1$, $z_i^N=X_i$
and $t_i^N=t_i$. Conversely, every transitive model of this kind
which has cardinality $\omega_1$ has a code $Y$ satisfying $V_1$.
\end{lemma}
\begin{proof}
The equality relation, the membership relation and the satisfaction
relation are given by the three parts of the class $Y\subseteq\omega_1$.
The finitely many conditions which express that $Y$ is the code of a
model, and the Tarski conditions for the satisfaction relation,
quantify over countable ordinals only, and hence
are first order over $\HC$. Every descending $\omega$-sequence
of closed terms is coded by an element of $\HC$, so well-foundedness
is first order as well. The model $N$ is the transitive collapse of
the quotient of the set of closed terms by the equality relation.

We require $c_0=\varnothing$ and $c_{n+1}=c_n\cup\{c_n\}$ for
$n<\omega$. Thus the constants $c_n$, $n<\omega$, are interpreted by
the natural numbers.
For each $0<\alpha<\omega_1$, we require
that $c_\alpha$ is an ordinal, that $F_0(c_\alpha)$ has domain
$\omega$ and range $c_\alpha$, and that $H_0$ gives the corresponding
preimages. We require, for $n<\omega$,
\[
 F_0(c_\alpha)(c_n)=c_{f_{0,\alpha}(n)}.
\]
The surjections $f_{0,\alpha}$ are hereditarily countable, so
these conditions are first order over $\HC$. By induction on
$\alpha$, we see that $c_\alpha$ is interpreted by $\alpha$. We
require all these constants to belong to $b_1$. Since $N$ satisfies
$b_1=\omega_1$, it follows that $b_1^N$ is the true $\omega_1$:
it cannot be smaller because it contains every countable ordinal;
if it were larger, $N$ would contain a surjection from
$\omega$ onto the true $\omega_1$.

The conditions $z_i\subseteq b_1$ and
$c_\alpha\in z_i\leftrightarrow\alpha\in X_i$ ensure $z_i^N=X_i$. The
analogous conditions for the constants $c_n$ ensure $t_i^N=t_i$.
Conversely, given a transitive model as in the statement, interpret
the constants $g_\alpha$ by an enumeration of its elements and the
constants $c_\alpha$ by the ordinals $\alpha$, put
$F_0(\alpha)=f_{0,\alpha}$ for $0<\alpha<\omega_1$, and choose
interpretations of the function symbols $f_\varphi$ by recursion on
the languages $\mathcal L_m$. This gives a code $Y$ satisfying $V_1$.
\end{proof}

We now define the corresponding formulas for $r\geq2$ by simultaneous
recursion. Finitely many classes which are quantified together are
always combined into one class by the fixed pairing function. A
variable $E$ is intended to range over codes of ordinals: we write
$|E|_r$ for the order type of the wellorder coded by $Z_{r-1}(E)$.
We write $\PiT_0$ for the class of first-order
formulas.

\begin{lemma}\label{lem:recursive-V-explicit}\label{spec:signed}
For $r\geq2$ there are formulas
\[
 O_r,C_r^+,C_r^-\in\PiT_{r-2},\qquad V_r\in\PiT_{r-1},
\]
defined uniformly in $r$, with the following properties in the final
model.
\begin{enumerate}[label=(\alph*)]
\item $O_r(E)$ holds if and only if $Z_{r-1}(E)$ is a code of a
wellorder of a subset of $\omega_{r-1}$, and
$\{|E|_r:O_r(E)\}=\omega_r$.
\item If $O_r(E)$ and $a=|E|_r$, then
\[
 \begin{aligned}
 a\in Z_r(X)&\ \Longleftrightarrow\ \exists H\,C_r^+(X,E;H),\\
 a\notin Z_r(X)&\ \Longleftrightarrow\ \exists H\,C_r^-(X,E;H).
 \end{aligned}
\]
In particular, exactly one of the two formulas on the right has a
witness $H$.
\item If $V_r(Y;\vec X,\vec t)$ holds, then $Z_r(Y)$ is a code
of a well-founded extensional model whose transitive collapse $N$
satisfies $\mathsf T_r$ and
\[
 b_j^N=\omega_j\ (j\leq r),\qquad
 z_\ell^N=Z_r(X_\ell),\qquad t_\ell^N=t_\ell.
\]
The same holds if finitely many further constants are required to
denote sets of the form $Z_r(X)$.
\item Conversely, let $N$ be a transitive model of
$\mathsf{ZFC}^-$, with interpretations of the constants of
$\mathcal L_0$, which satisfies the axioms (i)--(iii) of
$\mathsf T_r$ and the conditions displayed in (c). Suppose that $N$
belongs to a model of (C2) for $r$ and has cardinality $\omega_r$
in that model. Then the constants $c_\alpha,g_\alpha$ and the
function symbols can be interpreted in such a way that $N$ becomes a
model of $\mathsf T_r$ in which every element is the value of a
closed term, and there is $Y$ with $V_r(Y;\vec X,\vec t)$ such that
$Z_r(Y)$ is the code of this model.
\end{enumerate}
\end{lemma}
\begin{proof}
We argue by induction on $r$. For each $r$ we define $O_r$, then
$C_r^+$ and $C_r^-$ together with auxiliary formulas
$W_r^\theta$, and finally $V_r$. For $r=2$ the formulas of the
previous step are replaced by $V_1$ and by first-order statements
about subsets of $\omega_1$.

\emph{The formula $O_r$.}
$O_r(E)$ expresses that the relation coded by $Z_{r-1}(E)$ is a
wellorder. For $r=2$ we have $Z_1(E)=E$, and this is first order
over $\HC$. For $r>2$, membership of an ordinal below $\omega_{r-1}$
in $Z_{r-1}(E)$ is expressed by means of a code $E'$ of this ordinal
and a witness $H$ for $C_{r-1}^+(E,E';H)$ or $C_{r-1}^-(E,E';H)$. The
axioms of linear orders, which concern finitely many elements, are
required to hold for all codes of ordinals and all such witnesses.
Well-foundedness is expressed by requiring that there is no class
which codes a descending $\omega$-sequence of codes of ordinals
together with witnesses, in the sense of $C^+_{r-1}$, showing that
each of them precedes the previous one in the order coded by
$Z_{r-1}(E)$. Since the formulas of the previous step are
$\PiT_{r-3}$, we obtain $O_r\in\PiT_{r-2}$. Every
$a<\omega_r$ is the order type of a constructible wellorder of a
subset of $\omega_{r-1}$, whose code is a constructible subset of
$\omega_{r-1}$. If $r>2$, (C2) for $r-1$ gives $E$ such that
$Z_{r-1}(E)$ is this code; for $r=2$ the code itself is a subset of
$\omega_1$. This proves (a).

\emph{The formulas $C_r^\pm$ and $W^\theta_r$.}
In $C_r^\varepsilon(X,E;H)$, for $\varepsilon\in\{+,-\}$, the class
$H$ is required to code the following objects: a class $Y'$ which
satisfies $V_{r-1}$ and thus yields a transitive model $N$ of
$\mathsf T_{r-1}$; constants of this model which denote the sets
decoded from the given classes, an ordinal $a$ and an isomorphism
$\pi$; and, for each of the finitely many instances of the
satisfaction relation of $N$ which are needed below, a witness in the
sense of $C^+_{r-1}$ that this instance belongs to the code
$Z_{r-1}(Y')$. In terms of the model $N$, this means
\[
 \begin{gathered}
 \rho=\omega_{r-1},\qquad \beta=(\rho^+)^N,\\
 z=Z_{r-1}(X)\in N,\qquad e=Z_{r-1}(E)\in N,\\
 a<\beta,\qquad\pi:e\cong a,
 \end{gathered}
\]
where $\pi$ is an isomorphism between the wellorder coded by $e$ and
the ordinal $a$; the same applies if finitely many classes $X_i,E_i$
are given. The formula $V_{r-1}$ ensures that the constants denote the
sets $z$ and $e$. Since the witnesses for the instances of the
satisfaction relation are part of $H$, no existential class
quantifier is needed to express that these instances hold in $N$.
Hence the conditions on the model, on the constants and on the
instances of satisfaction are $\PiT_{r-2}$. For $r=2$, the model is
coded directly by a subset of $\omega_1$ and these conditions are
first order.

In addition, $C_r^\varepsilon(X,E;H)$ requires that $N$ satisfies
\[
 \sup(z\cap B^{r-1}_a)<\rho\quad(\varepsilon=+),\qquad
 \sup(z\cap B^{r-1}_a)=\rho\quad(\varepsilon=-).
\]
Since the definition of the almost disjoint family is total and
functional over $L_\beta^N$ and $a<\beta$, the upward absoluteness of
$\Sigma_1$ formulas shows that $N$ computes $B^{r-1}_a$ correctly.
Whether $z\cap B^{r-1}_a$ is bounded in $\rho$ is absolute between $N$
and the universe, and $a=|E|_r$ because $\pi\in N$. Hence a witness
for $C^+_r(X,E;H)$ shows $a\in Z_r(X)$, and a witness for
$C^-_r(X,E;H)$ shows $a\notin Z_r(X)$.
In the same way we obtain formulas $W_r^\theta\in\PiT_{r-2}$ such that
$\exists K\,W^\theta_r(\vec E;K)$ expresses that the relation $\theta$
holds between the ordinals $|E_1|_r,|E_2|_r,\ldots$, for each of the
following relations $\theta$: equality and order of ordinals, being
smaller than the true $\omega_j$ (for $j<r$), the graphs of the
fixed coding functions and of the functions $f_{j,a}$, and the
recursive relations on $\omega$ which describe the syntax of
$\mathcal L$ and the axioms of $\mathsf T_r$. These
relations are absolute between $N$ and the universe, in contrast to
arbitrary statements which hold in $N$.

For the converse implications in (b), let $X$ and $E$ be given, with
$O_r(E)$ and $a=|E|_r$. By (C1), $X$ and $E$ belong to $B[G_C]$ for a
closed set $C$ as in (C1). If $r>2$, we work in the corresponding
model of (C2) for $r-1$, which contains $B[G_C]$; if $r=2$, we work in
the final model. Let $\mathcal G$ be the set of the graphs of the
fixed coding functions, including the pairing function and the coding
of terms. With $\rho=\omega_{r-1}$, $z=Z_{r-1}(X)$, $e=Z_{r-1}(E)$
and the isomorphism $\pi:e\cong a$, choose
\[
 \begin{gathered}
 M\prec H_\vartheta,\quad |M|=\rho,\quad
 \rho\cup(a+1)\subseteq M,\\
 \rho,z,e,a,\pi,L_{\omega_r},
 \operatorname{Sat}(L_{\omega_r}),\mathcal G\in M,
 \end{gathered}
\]
for a sufficiently large $\vartheta$. Let $\bar\pi:M\to N$ be the
transitive collapse. Condensation gives
\[
 \bar\pi(\rho)=\rho,\qquad
 \beta=\bar\pi(\omega_r)=(\rho^+)^N>a,\qquad
 \bar\pi(L_{\omega_r})=L_\beta,
\]
and by elementarity the definitions mentioned in $\mathsf T_{r-1}$
are total over $L_\beta$. The model $N$ belongs to the model in
which we work and has cardinality $\rho$ there. For $r>2$, part (d)
of the induction hypothesis gives $Y'$ satisfying $V_{r-1}$ such that
$Z_{r-1}(Y')$ is a code of $N$, and
part (b) of the induction hypothesis gives the witnesses for the
required instances of the satisfaction relation. For $r=2$, a code
of $N$ is itself a subset of $\omega_1$, and we use
Lemma~\ref{lem:direct-V1-explicit}. This yields a witness $H$ for
$C^+_r$ or $C^-_r$, according to whether $a\in Z_r(X)$ or not. The
same argument applies to finitely many classes and to the formulas
$W^\theta_r$. This proves (b).

\emph{The formula $V_r$.}
Membership of ordinals in $Z_r(Y)$ is again expressed by means of
codes of ordinals and witnesses. The formula $V_r$ is a conjunction
of finitely many conditions, each of the following kind: for all
ordinals $\alpha_1,\ldots,\alpha_m<\omega_r$ which stand in a given
relation $\theta$ of the kind considered above, a given propositional
combination $A$ of statements of the form $\alpha_i\in Z_r(Y)$ or
$\alpha_i\in Z_r(X_\ell)$ is true. We
express such a condition by saying that there is no witness for a
violation of it. Here a witness for a violation consists of
\begin{itemize}
\item codes $E_1,\ldots,E_m$ of ordinals, that is, classes with
$O_r(E_i)$;
\item a class $K$ with $W^\theta_r(\vec E;K)$, which shows that the
ordinals $|E_i|_r$ stand in the relation $\theta$;
\item signs $\varepsilon_i\in\{+,-\}$, one for each of the
statements which occur in $A$, such that $A$ is false when these
statements are given the truth values $\varepsilon_i$, and classes
$H_i$ with $C_r^{\varepsilon_i}(Y,E_i;H_i)$ or
$C_r^{\varepsilon_i}(X_\ell,E_i;H_i)$, respectively.
\end{itemize}
Since there are only finitely many choices of signs, the statement
that $(\vec E,K,\vec H)$ is a witness for a violation is a finite
disjunction of conjunctions of $\PiT_{r-2}$ formulas; we denote it by
$\operatorname{Viol}_r(Y;\vec E,K,\vec H)$. The condition is expressed
by
\[
 \forall\vec E\,\forall K\,\forall\vec H\,
 \neg\operatorname{Viol}_r(Y;\vec E,K,\vec H),
 \qquad\operatorname{Viol}_r\in\PiT_{r-2},
\]
which is $\PiT_{r-1}$. By (a), every ordinal below $\omega_r$ has a
code. As shown in the proof of (b), every true instance of $\theta$
has a witness in the sense of $W^\theta_r$, and for every code $E$
exactly one of the two signs has a witness. Hence the displayed formula expresses
the intended condition on $Z_r(Y)$. The ordinals $q_{\rm eq}(u,v)$,
$q_{\rm mem}(u,v)$ and $q_{\rm sat}(e,\vec u)$ and the codes of terms
enter these conditions through the relation $\theta$, which contains
the graphs of the fixed coding functions.

In this form, $V_r$ requires the axioms of equality, including the
congruence conditions for membership, for the function symbols and for
satisfaction; the Tarski conditions for the satisfaction relation,
with existential quantifiers treated by means of the function symbols
$f_\varphi$; and the axioms of $\mathsf T_r$. Writing $S$ for the part of the code
which codes satisfaction, the conditions for the existential
quantifier are
\[
 \begin{aligned}
 S(\exists v\,\varphi(v,\vec u))
   &\Rightarrow S(\varphi(f_\varphi(\vec u),\vec u)),\\
 S(\varphi(w,\vec u))
   &\Rightarrow S(\exists v\,\varphi(v,\vec u))
   \quad(\text{for all closed terms }w).
 \end{aligned}
\]
Thus an element which satisfies an existential formula is always
given by the term $f_\varphi(\vec u)$, and no existential quantifier
over closed terms is needed.

For well-foundedness, $V_r$ requires that there is no class which
codes, for every $n<\omega$, the following classes: codes $E_n$ and
$E'_n$ of ordinals; a witness, in the sense of $W^\theta_r$, that
$|E_n|_r$ and $|E_{n+1}|_r$ are codes of closed terms $u_n$ and
$u_{n+1}$ and that $|E'_n|_r=q_{\rm mem}(u_{n+1},u_n)$; and a
witness, in the sense of $C^+_r$, that $|E'_n|_r\in Z_r(Y)$. These
countably many classes are combined into one class. The condition on
this class is $\PiT_{r-2}$, since a quantifier over $n<\omega$ can be
moved inside the class quantifiers. Hence the existence of such a
class is a $\Sig_{r-1}$ statement, and its negation is
$\PiT_{r-1}$.

Next, $V_r$ requires, in the same form,
\[
 c_0=\varnothing,\qquad c_{n+1}=c_n\cup\{c_n\}\ (n<\omega),
 \qquad c_\alpha\text{ is an ordinal}.
\]
Thus the constants $c_n$, $n<\omega$, are interpreted by the natural
numbers.
For $j<r$, $0<\alpha<\omega_{j+1}$ and $i<\omega_j$, it requires
\[
 \begin{gathered}
 c_\alpha\in b_{j+1},\qquad
 F_j(c_\alpha):b_j\twoheadrightarrow c_\alpha,\\
 u\in c_\alpha\Rightarrow
 \bigl[H_j(c_\alpha,u)\in b_j\ \land\
 F_j(c_\alpha)(H_j(c_\alpha,u))=u\bigr],\\
 F_j(c_\alpha)(c_i)=c_\gamma
 \quad(\gamma=f_{j,\alpha}(i)).
 \end{gathered}
\]
The last equation is required for all $\gamma$ and all witnesses, in
the sense of $W^\theta_r$, for $\gamma=f_{j,\alpha}(i)$.
The condition $\alpha<\omega_{j+1}$ is expressed by $W_r^\theta$ if
$j+1<r$, and it holds for all codes of ordinals if $j+1=r$; the
condition $i<\omega_j$ is treated in the same way. The witnesses and
the closed terms $u$ are universally quantified. Note that the term
$H_j(c_\alpha,u)$ provides the preimage of an element of $c_\alpha$;
therefore no existential quantifier over codes of ordinals is needed.

We now verify that these requirements have the intended effect. Let
$N$ be the transitive collapse of the model coded by $Z_r(Y)$. Then
$b_0^N=\omega$. If $b_j^N=\omega_j$, where $j<r$, then induction on
$0<\alpha<\omega_{j+1}$ gives
\[
 c_\alpha^N=\{c^N_{f_{j,\alpha}(i)}:i<\omega_j\}
          =\{f_{j,\alpha}(i):i<\omega_j\}=\alpha.
\]
Thus $b_{j+1}^N\geq\omega_{j+1}$. If the inequality were strict, then,
since $N\models b_{j+1}=(b_j)^+$, the model $N$ would contain a
surjection from $\omega_j$ onto $\omega_{j+1}$.
Hence $b_j^N=\omega_j$ for all $j\leq r$.

Finally, $V_r$ requires $z_i\subseteq b_r$ and, for all codes $E$ of
ordinals and all witnesses for the two sides, that
\[
 |E|_r\in Z_r(X_i)
 \quad\Longleftrightarrow\quad
 q_{\rm mem}(c_{|E|_r},z_i)\in Z_r(Y);
\]
here the ordinal $q_{\rm mem}(c_{|E|_r},z_i)$ is identified by
$W_r^\theta$. By (a) and (b), every ordinal below $\omega_r=b_r^N$ has
a code and the corresponding witnesses exist, so we obtain
$z_i^N=Z_r(X_i)$. The analogous conditions for the
constants $c_n$, $n<\omega$, ensure $t_i^N=t_i$.
All these conditions are $\PiT_{r-1}$. This proves (c).

For (d), let $N$ be as in the statement, and work in the given model
of (C2) for $r$. Interpret the constants $g_\alpha$ by an enumeration
of $N$ of length $\omega_r$ and the constants $c_\alpha$ by the
ordinals $\alpha<\omega_r=b_r^N$. For $j<r$ and $0<a<b_{j+1}^N$ put
$F_j(a)=f_{j,a}$, which belongs to $N$ by axiom (ii), and let
$H_j(a,u)$ be the least preimage of $u$ under $f_{j,a}$. Finally
choose interpretations of the function symbols $f_\varphi$, by
recursion on $m$ for the symbols of $\mathcal L_{m+1}$, such that
$f_\varphi(\vec u)$ satisfies $\varphi(v,\vec u)$ whenever some
element does. The resulting structure is a model of $\mathsf T_r$ in
which every element is the value of a closed term. Its code,
consisting of the equality, membership and satisfaction relations, is
a subset of $\omega_r$ which belongs to the model of (C2), and the
structure satisfies all conditions required by $V_r$. By (C2) there
is $Y$ in the final model such that $Z_r(Y)$ is this code. The
witnesses are sound by (b): a witness for a violation of one of the
conditions of $V_r$ would show that the structure violates this
condition, and a class as in the condition on well-foundedness would
yield a descending $\in$-sequence in $N$. Hence
$V_r(Y;\vec X,\vec t)$ holds.

This completes the definition of $O_r$, of $C_r^\pm$ and
$W_r^\theta$, and of $V_r$, and the proof of their properties.
The formulas are obtained recursively from $r$. The sets $Z_r(X)$ and
the true cardinals $\omega_j$ were used above only to describe the
meaning of the formulas; the formulas themselves refer to them only
through the witnesses and the constants, and they have no parameters.
Whenever countably many classes are needed, for sequences or for the
axiom schemes, they are combined into one class, using
$\omega_1\times\omega=\omega_1$.
\end{proof}

\subsection{The definition of \texorpdfstring{$\Code_n$ for $n\geq2$}{Code-n for n>=2}}

For each $r\geq1$ fix the following first-order formula
$\chi_r(\rho,\beta,z,t)$ of the language of set theory. It asserts
that there are sets $P,C$ such that
\[
 P=\{\eta<\beta:z\cap B^r_\eta
                        \text{ is bounded in }\rho\},
 \qquad C=\operatorname{Even}(P),
\]
where the almost disjoint family is computed over $L_\beta$ by means
of $\Phi_B$, and such that
$\psi_{\rm club}(\rho,\beta,C,t)$ holds. Let $q_r$ be the ordinal
which codes, in the coding fixed above, the instance of the
satisfaction relation for the formula $\chi_r$ and the terms
$b_r,b_{r+1},z_1,t_1$. Thus $q_r$ is determined by the fixed coding
of formulas and terms; it does not depend on the model which is
coded.

For $r=1$, let $D_1(Y,t)$ be the first-order statement that $q_1$
belongs to the part of $Y$ which codes satisfaction. For $r\geq2$
put
\[
 D_r(Y,t)\quad\Longleftrightarrow\quad
 \exists E,H,K\,
   \bigl[W_r^{=q_r}(E;K)
                 \wedge C_r^+(Y,E;H)\bigr],
\]
where $W_r^{=q_r}(E;K)$ is the formula $W^\theta_r$ which expresses
$|E|_r=q_r$. The three class variables can be combined into one.
Hence $D_1$
is first order and $D_r\in\Sig_{r-1}$ for $r\geq2$.
If $V_r(Y;X,t)$ holds and $N$ is the transitive collapse of the model
coded by $Z_r(Y)$, then
\[
 D_r(Y,t)
 \quad\Longleftrightarrow\quad
 N\models\chi_r(\rho,\beta,z,t),
\]
where $\rho=b_r^N=\omega_r$, $\beta=b_{r+1}^N$ and
$z=z_1^N=Z_r(X)$. Moreover, since the definition of the almost
disjoint family is total over $L_\beta$, the set $P$ defined in $N$
is $Z_{r+1}(X)\cap\beta$.

For $r\geq1$, put
\[
 \Code_{r+1}(t)\quad\Longleftrightarrow\quad
 \exists X\,\forall Y\,
       [V_r(Y;X,t)\longrightarrow D_r(Y,t)].
\]
For $r=1$, $V_1(Y;X,t)\to D_1(Y,t)$ is first order; for $r\geq2$,
$V_r$ is $\PiT_{r-1}$ and $D_r$ is $\Sig_{r-1}$. Thus
$\Code_{r+1}\in\Sig_{r+1}$, uniformly in $r$.

\begin{proposition}\label{prop:exact}
In the final model, for every $n\geq1$ and $t\in I$,
\[
                  \Code_n(t)\quad\Longleftrightarrow\quad t\in T_n.
\]
\end{proposition}
\begin{proof}
Let $n=1$. If $t\in T_1$, then $C_{\vec T}$ adds a club disjoint from
$S^1_t$, so $\Code_1(t)$ holds. If $t\notin T_1$, then $S^1_t$ is
stationary by (C3), so $\Code_1(t)$ fails.

Let $n=r+1\geq2$ and $t\in T_n$. By (C3) there is a club
$C\subseteq\omega_n$ disjoint from $S^n_t$ such that
$P=P_t^n(C)\subseteq\omega_n$ belongs to $B_{\vec T}$.
By (C2) there is $X$ with $Z_n(X)=P$. Suppose that
$V_r(Y;X,t)$ holds, and let $N$ be the corresponding transitive
model. Then $\rho=b_r^N=\omega_r$,
$\beta=b_{r+1}^N<\omega_n$, and $z=Z_r(X)\in N$. As noted above, the
set $\{\eta<\beta:z\cap B^r_\eta\text{ is bounded in }\rho\}$,
computed in $N$, is $P\cap\beta$; in particular $P\cap\beta\in N$.
Since $N$ is suitable for $n$, Lemma~\ref{lem:direct-padding} gives
$N\models\psi_{\rm club}(\rho,\beta,C\cap\beta,t)$.
Thus $N\models\chi_r(\rho,\beta,z,t)$, and $D_r(Y,t)$ follows. Since
this holds for every $Y$ with $V_r(Y;X,t)$, the set $X$ witnesses
$\Code_n(t)$.

Conversely let $n=r+1\geq2$ and $t\notin T_n$, and suppose that $X$
witnesses $\Code_n(t)$. By (C1), $X$ belongs to $B[G_C]$ for a closed
set $C$ as in (C1), and hence to the corresponding model
$W_r[G_{Q_r}]$. We work in this model. For every
$\delta<\omega_n$, choose an elementary submodel of a
sufficiently large $H_\theta$ which has size $\omega_r$ and contains
$\omega_r$ and $\delta+1$ as subsets and the set $z=Z_r(X)$,
the real $t$, $L_{\omega_n}$, its
satisfaction relation, and the graphs of the fixed coding functions,
including the pairing function and the coding of terms, as elements.
Its transitive collapse $N_\delta$ satisfies $\mathsf{ZFC}^-$ and
the axioms (i)--(iii) of $\mathsf T_r$, with
$b_r^{N_\delta}=\omega_r$ and
$\beta_\delta=b_{r+1}^{N_\delta}>\delta$.

The model $N_\delta$ belongs to $W_r[G_{Q_r}]$ and has cardinality
$\omega_r$ there. If
$r\geq2$, recall from the proof of
Proposition~\ref{prop:cascade-full} that $W_r[G_{Q_r}]$ is contained
in the model $W_{r-1}[G_{Q_{r-1}}]$ which is defined from the
coordinates of $C$ of levels below $r-1$; this is a model of (C2) for
$r$. Hence Lemma~\ref{lem:recursive-V-explicit}(d) gives
$Y_\delta$ with $V_r(Y_\delta;X,t)$ such that $Z_r(Y_\delta)$ codes
$N_\delta$. If $r=1$, a code of $N_\delta$ itself is such a set
$Y_\delta$, by Lemma~\ref{lem:direct-V1-explicit}. Since $X$
witnesses $\Code_n(t)$, we obtain $D_r(Y_\delta,t)$.

Put $P=Z_n(X)$ in the final model. As noted above,
$D_r(Y_\delta,t)$ means that
$\operatorname{Even}(P)\cap\beta_\delta$ is a club in
$\beta_\delta$ which is disjoint from $S^n_t$. The ordinals
$\beta_\delta$, $\delta<\omega_n$, are unbounded in $\omega_n$. By
Lemma~\ref{lem:club-from-segments}, $\operatorname{Even}(P)$ is a
club in $\omega_n$ disjoint from $S^n_t$. This contradicts (C3).
\end{proof}

\section{Comparing the forcings \texorpdfstring{$F_{\vec T}$}{F-T}}\label{sec:comparison}

Put $F_{\vec T}=C_{\vec T}*\dot D_{\vec T}$ and write
$M_{\vec T}=B_{\vec T}[G_{D_{\vec T}}]$. Throughout this
section, the stationary sets and the forcings $C_n(T)$ are those
defined in $L$, and all products are formed in $L$. By
Proposition~\ref{prop:compile}, for each $\PiT_n$ formula $A(a)$
there is a $\PiT_1(H_{\omega_n})$ formula $\Phi_A(a)$, independent
of $\vec T$, such that
\[
 M_{\vec T}\models A(a)\quad\Longleftrightarrow\quad
 B_{\vec T}\models\Phi_A(a).
\]
The corresponding assertion holds for $\Sig_n$ formulas. Note that
this equivalence holds for every $D_{\vec T}$-generic filter over
$B_{\vec T}$. Therefore we never have to compare the forcings
$D_{\vec R}$ and $D_{\vec T}$ for different sequences directly.

Whenever $\vec R\leq\vec T$ and a $C_{\vec T}$-generic filter $G$
over $L$ is given, we let $B_{\vec R}$ be the extension of $L$ by the
restriction of $G$ to the complete subforcing $C_{\vec R}$. Thus
$B_{\vec R}\subseteq B_{\vec T}$. The statements about the forcing
relations of $F_{\vec R}$ and $F_{\vec T}$ proved below do not depend
on this choice of the generic filters.

\begin{lemma}\label{lem:comparison-domain}
Suppose that $\vec R\leq\vec T$ and $R_m=T_m$ for all $m<k$. Then,
for generic filters chosen as above,
\[
 H_{\omega_k}^{B_{\vec R}}=H_{\omega_k}^{B_{\vec T}}.
\]
\end{lemma}
\begin{proof}
Write
\[
 P=\prod_{m<k}C_m(T_m),\qquad
 U=\prod_{m\geq k}C_m(R_m),\qquad
 R'=\prod_{m\geq k}C_m(T_m\setminus R_m).
\]
Thus $C_{\vec R}=U\times P$ and $C_{\vec T}=U\times R'\times P$.
By Lemma~\ref{lem:source}, $U\times R'$ is
$<\omega_k$-distributive over $L$, so $R'$ is
$<\omega_k$-distributive over $L[G_U]$. Since $|P|<\omega_k$,
Lemma~\ref{lem:small} gives the same distributivity over
$L[G_U][G_P]$. The models $B_{\vec R}$ and $B_{\vec T}$ therefore have the same
sequences of ordinals of length $<\omega_k$, and hence the same
$H_{\omega_k}$.
\end{proof}

\begin{proposition}\label{prop:ownaxis}
Suppose that $\vec R\leq\vec T$ and $R_m=T_m$ for all $m\ne n$. If
$A(a)$ is a $\PiT_n$ formula with real parameters and
$F_{\vec T}\Vdash A(a)$, then
$F_{\vec R}\Vdash A(a)$.
\end{proposition}
\begin{proof}
Let $B_{\vec R}\subseteq B_{\vec T}$ be as above. By
Lemma~\ref{lem:comparison-domain} with $k=n$, the two models have
the same $H_{\omega_n}$. Every class over $H_{\omega_n}$
in $B_{\vec R}$ is also a class in $B_{\vec T}$, and the first-order
part of $\Phi_A$ has the same meaning in both models for such a
class. Since $\Phi_A$ is $\PiT_1(H_{\omega_n})$,
\[
 B_{\vec T}\models\Phi_A(a)\quad\Longrightarrow\quad
 B_{\vec R}\models\Phi_A(a).
\]

Now assume $F_{\vec T}\Vdash A(a)$, and let $G'$ be an arbitrary
$C_{\vec R}$-generic filter over $L$. Extend $G'$ to a
$C_{\vec T}$-generic filter by forcing with
$C_n(T_n\setminus R_n)$. By the hypothesis and
Proposition~\ref{prop:compile}, $\Phi_A(a)$ holds in the larger model
$B_{\vec T}$, and by the implication just proved it holds in
$L[G']$. Since $G'$ was arbitrary, Proposition~\ref{prop:compile}
gives $F_{\vec R}\Vdash A(a)$.
\end{proof}

\begin{proposition}\label{prop:higheraxis}
Suppose that $R_m=T_m$ for all $m\leq n$, and let $A(a)$ be a
$\Sig_n$ or $\PiT_n$ formula with real parameters. If
$\vec R\leq\vec T$ and $B_{\vec R}\subseteq B_{\vec T}$ are as above,
then $\Phi_A(a)$ has the same truth value in $B_{\vec R}$ and in
$B_{\vec T}$. Consequently, whether or not $\vec R\leq\vec T$,
$F_{\vec R}$ and $F_{\vec T}$ force the same $\Sig_n$ and $\PiT_n$
statements with real parameters.
\end{proposition}
\begin{proof}
For $\vec R\leq\vec T$, Lemma~\ref{lem:comparison-domain} with
$k=n+1$ gives $H_{\omega_{n+1}}^{B_{\vec R}}=H_{\omega_{n+1}}^{B_{\vec T}}$.
In particular they have the same $H_{\omega_n}$ and the same
subsets of it, since $|H_{\omega_n}|=\omega_n$. Thus $\Phi_A(a)$
has the same truth value in both models.

Every $C_{\vec R}$-generic filter extends to a $C_{\vec T}$-generic
filter, and every $C_{\vec T}$-generic filter restricts to a
$C_{\vec R}$-generic filter. Proposition~\ref{prop:compile} now
shows that $F_{\vec R}$ and $F_{\vec T}$ force the same statements of
this kind. For incomparable
sequences, compare each with their coordinatewise intersection.
\end{proof}

\begin{proposition}\label{prop:determination}
For every $\vec T\in L$, the forcing $F_{\vec T}$ preserves cardinals
and adds no reals, and for every second-order formula $\varphi$ and
all reals $a\in L$, the trivial condition of $F_{\vec T}$ decides
$\varphi(a)$. Thus $F_{\vec T}$ satisfies \textup{(A1)}.
\end{proposition}
\begin{proof}
Every second-order formula is equivalent to a $\Sig_n$ or
$\PiT_n$ formula for some finite $n$. By
Proposition~\ref{prop:compile}, the truth value of such a formula
$A(a)$ in the extension by $D_{\vec T}$ is the truth value in
$B_{\vec T}$ of the formula $\Phi_A(a)$, whose only parameters are
the constructible reals $a$. Since $C_{\vec T}$ is cone homogeneous
by Lemma~\ref{lem:source}, Fact~\ref{fact:cone} shows that the
trivial condition of $C_{\vec T}$ decides $\Phi_A(a)$.
By Proposition~\ref{prop:compile}, which holds for every
$D_{\vec T}$-generic filter, the trivial condition of
$F_{\vec T}$ decides $A(a)$.

The forcing $C_{\vec T}$ preserves cardinals and adds no reals by
Lemma~\ref{lem:source}, and $D_{\vec T}$ preserves cardinals and adds
no reals by Theorem~\ref{lem:master}\ref{M:cc}. Hence the
two-step iteration $F_{\vec T}$ preserves cardinals and adds no
reals.
\end{proof}

\begin{remark}\label{rem:uniform}
The map $A\mapsto\Phi_A$ of Proposition~\ref{prop:compile} is
recursive. The statement $\Phi_A(a)$ is a statement about the
structure $(H_{\omega_n},\mathcal P(H_{\omega_n}),\in)$ of
$B_{\vec T}$, which is a set, and the forcing relation of
$C_{\vec T}$ for statements about such structures is definable in
$L$, uniformly in $\vec T$ and in $A$. Hence the relation
``$1_{F_{\vec T}}\Vdash A(a)$'' is definable in $L$.
\end{remark}

\subsection{Verification of the four hypotheses}

\begin{theorem}\label{thm:A1-A4}
For $\vec T\in L$ put $F_{\vec T}=C_{\vec T}*\dot D_{\vec T}$.
The family $\langle F_{\vec T}\rangle$ is definable in $L$ and
satisfies \textup{(A1)--(A4)} of Theorem~\ref{thm:abstract}. The
relation ``$1_{F_{\vec T}}\Vdash\varphi(a)$'' is definable in $L$, as
a relation between G\"odel numbers of second-order formulas,
sequences $\vec T$ and real parameters.
\end{theorem}
\begin{proof}
The definability of the family was shown after
Definition~\ref{def:coding-recursion}.
Proposition~\ref{prop:determination} gives (A1).
The parameter-free formulas $\Code_n$ of
Proposition~\ref{prop:exact}, which do not depend on $\vec T$,
give (A2).
Propositions~\ref{prop:ownaxis} and~\ref{prop:higheraxis}
give (A3) and (A4), respectively.
Remark~\ref{rem:uniform} gives the definability of the forcing
relation.
\end{proof}

\section{The choice of \texorpdfstring{$\vec T_*$}{T*} and the uniformizing functions}\label{sec:selection}

In this section we prove Theorem~\ref{thm:abstract}. The proof uses
only the hypotheses (A1)--(A4) and does not depend on the particular
definition of the forcings $F_{\vec T}$.
Theorem~\ref{thm:A1-A4} shows that these hypotheses hold for
$F_{\vec T}=C_{\vec T}*\dot D_{\vec T}$.

\subsection{The recursion in \texorpdfstring{$L$}{L}}

Fix a recursive enumeration $\langle A_{n,e}:e<\omega\rangle$ of the
$\PiT_n$ formulas $A_{n,e}(x,y)$. We work in $L$. For $n\geq1$,
$e<\omega$ and $x\in\R^L$, we say that the section of $A_{n,e}$ at
$x$ is \emph{nonempty at $\vec T$} if
\[
 1_{F_{\vec T}}\Vdash A_{n,e}(x,y)\qquad\text{for some }y\in\R^L.
\]
By (A1), every statement $A_{n,e}(x,y)$ with $x,y\in\R^L$ is decided
by the trivial condition of $F_{\vec T}$. Since $F_{\vec T}$ adds no
reals, it follows that the section of $A_{n,e}$ at $x$ is nonempty at
$\vec T$ if and only if the set $\{y:A_{n,e}(x,y)\}$ is nonempty in
the extension by $F_{\vec T}$.

Put $\vec T^0=\langle I:n\geq1\rangle$. Let $n\geq1$ and suppose that
$\vec T^{n-1}$ has been defined, with $T_n^{n-1}=I$. We define
sequences $\vec T^{n,\alpha}$, which agree with $\vec T^{n-1}$ at all
coordinates other than $n$, and a partial function $a_n$ on the
pairs $(e,x)$, by recursion on $\alpha$.
Set $\vec T^{n,0}=\vec T^{n-1}$ and begin with the empty function.
At stage $\alpha$:
\begin{enumerate}[label=(\roman*)]
\item if there is a pair $(e,x)$ to which no value has been assigned
yet and such that the section of $A_{n,e}$ at $x$ is nonempty at
$\vec T^{n,\alpha}$, take the $<_L$-least such
$(e,x)$ and the $<_L$-least real $y$ with
$1_{F_{\vec T^{n,\alpha}}}\Vdash A_{n,e}(x,y)$.
Put $a_n(e,x)=y$ and remove only $\ell(e,x,y)$
from $T_n$, that is,
$T_n^{n,\alpha+1}=T_n^{n,\alpha}\setminus\{\ell(e,x,y)\}$;
\item at a limit stage take the intersection of the sets
$T_n^{n,\alpha}$ defined so far and the union of the partial
functions;
\item the recursion stops if there is no pair as in (i).
\end{enumerate}
Write $\vec T^n$ for the sequence at which the recursion stops and
$a_n$ for the resulting partial function. Finally put
\[
 T_n^*=T_n^n,\qquad \vec T_*=\langle T_n^*:n\geq1\rangle.
\]

\begin{lemma}\label{lem:recursion}
For every $n$, the recursion stops before $\omega_2$.
It is definable in $L$, and $\vec T_*\in L$.
Moreover, if $a_n(e,x)=y$, then
\[
 1_{F_{\vec T^n}}\Vdash A_{n,e}(x,y),
\]
and $a_n(e,x)$ is defined for every pair $(e,x)$ such that the
section of $A_{n,e}$ at $x$ is nonempty at $\vec T^n$.
\end{lemma}
\begin{proof}
The least pair and the least value in (i) exist by the definition
of nonemptiness, since $<_L$ is a wellorder.
Every step of the form (i) assigns a value to a new pair $(e,x)$.
There are only $\omega_1$ such pairs in $L$, so this cannot happen
at every stage below $\omega_2$.

Suppose that the value $y$ has been assigned to $(e,x)$ at stage
$\alpha$. For every later stage $\alpha'$ of this recursion,
including limit stages and the last stage, we have
$\vec T^{n,\alpha'}\leq\vec T^{n,\alpha}$, and these two sequences
differ only at the coordinate $n$. Hence (A3) shows that
$A_{n,e}(x,y)$ is forced by $F_{\vec T^{n,\alpha'}}$.
When the recursion stops, there is no pair as in (i), which is the
last assertion.

Since the relation ``$1_{F_{\vec T}}\Vdash\varphi(a)$'' is definable
in $L$, the recursion is a definable recursion in $L$. Replacement
gives the sequence $\langle\vec T^n:n\geq1\rangle$ and the
sequence $\vec T_*$.
\end{proof}

\subsection{Verification}

\begin{proof}[Proof of Theorem~\ref{thm:abstract}]
Fix $n$. The sequences $\vec T^n$ and $\vec T_*$ have the same
entries at all coordinates $m\leq n$. By (A4),
\[
 1_{F_{\vec T^n}}\Vdash A_{n,e}(x,y)
 \quad\Longleftrightarrow\quad
 1_{F_{\vec T_*}}\Vdash A_{n,e}(x,y)
\]
for all $e<\omega$ and all reals $x,y\in L$. In particular the same
sections are nonempty at $\vec T^n$ and at $\vec T_*$.

In an extension by $F_{\vec T_*}$ define
\[
 U_{n,e}(x,y)\quad\Longleftrightarrow\quad
 A_{n,e}(x,y)\ \wedge\ \neg\Code_n(\ell(e,x,y)).
\]
This is $\PiT_n$. By (A2) and the injectivity of $\ell$,
\[
 \neg\Code_n(\ell(e,x,y))
 \quad\Longleftrightarrow\quad\ell(e,x,y)\notin T^*_n
 \quad\Longleftrightarrow\quad a_n(e,x)=y.
\]
Indeed, at the beginning of the recursion for $n$ the set $T_n$
contains every code, and at each step of the form (i) exactly the code
$\ell(e,x,a_n(e,x))$ is removed; the recursions for the other
coordinates do not change $T_n$.
By Lemma~\ref{lem:recursion} and (A4), $A_{n,e}(x,a_n(e,x))$ holds in
the extension whenever $a_n(e,x)$ is defined.
Conversely, if the section $\{y:A_{n,e}(x,y)\}$ is nonempty in the
extension, then the section of $A_{n,e}$ at $x$ is nonempty at
$\vec T_*$ and hence at $\vec T^n$, so $a_n(e,x)$ is defined.
Consequently $U_{n,e}$ is a subset of $A_{n,e}$ and has
exactly one value at every $x$ with nonempty section.

Thus $\PiT_n$-$\UP$ holds for every $n\geq1$.
\end{proof}

\begin{proof}[Proof of Theorem~\ref{thm:main}]
Apply Theorem~\ref{thm:abstract} to the family
$F_{\vec T}=C_{\vec T}*\dot D_{\vec T}$, whose hypotheses
are verified in Theorem~\ref{thm:A1-A4}.
The sequence $\vec T_*$ belongs to $L$, and
$F_{\vec T_*}$ is the required cardinal-preserving set forcing
which adds no reals.
\end{proof}

\section{Projective determinacy}\label{sec:PD}

We first carry out the construction over a ground model of the form
$L[A]$, where $A\subseteq\omega_1$ is used as an additional predicate.
We then show that for $L(\mathbb R)^{M_\omega}$ the predicate $A$ can
be chosen with parameter-free $\Delta^2_1$ definitions, which allows
us to eliminate it.

\subsection{The construction over \texorpdfstring{$L[A]$}{L[A]}}

In this subsection fix the ground model
\[
 \mathcal W=L[A]\models\mathsf{ZFC},\qquad A\subseteq\omega_1^{\mathcal W}.
\]
Throughout this section, all ground-model objects are computed in
$\mathcal W$, and $\omega_n$, $\aleph_\alpha$ denote the cardinals of
$\mathcal W$.
The forcing extensions used in the construction preserve these cardinals.
Write $\Sigma^2_n[A]$ and $\Pi^2_n[A]$ for the pointclasses over
$(\HC,\in,A)$, with class quantifiers ranging over all subsets of
$\HC$ in the extension under consideration. The predicate $A$ always
denotes the same set of the ground model.

\begin{lemma}\label{lem:PD-source}
The stationary sets, the almost disjoint families and the sets
$P^n_t(C)$ of Section~\ref{sec:source} have analogues in
$\mathcal W$, defined with the canonical wellorder $<_{L[A]}$ in place
of $<_L$. For $n=1$, the diamond sequence is first-order definable
over $(\HC,\in,A)$, by a definition which is absolute for extensions
preserving $\HC$. For $n\geq2$, the definitions corresponding to
$\Phi_D$ and $\Phi_B$ are $\Sigma_1(A)$ over $L_{\omega_n}[A]$, and
Lemma~\ref{lem:direct-padding} holds with $L_\beta[A]$ in
place of $L_\beta$, for suitable models $N$ containing $A$.
\end{lemma}
\begin{proof}
Condensation for the hierarchy $L_\alpha[A]$ gives GCH in
$\mathcal W$ and
\[
 H_{\omega_n}^{\mathcal W}=L_{\omega_n}[A]\qquad(n\geq1).
\]
Indeed, for $n\geq2$ consider elementary submodels of size
$\omega_{n-1}$ which contain $\omega_{n-1}$ as a subset
and $A$ as an element. Their transitive collapses map $A$ to itself,
since $A\subseteq\omega_1$, and this shows that every member of
$H_{\omega_n}^{\mathcal W}$ belongs to $L_{\omega_n}[A]$.
For $x\in\HC^{\mathcal W}$, take a countable elementary submodel which
contains $A$ and $\omega_1$ as elements and the transitive closure of
$x$ as a subset. Its transitive collapse shows that
$x\in L_\beta[A\cap\delta]$ for countable ordinals $\beta,\delta$, and
$L_\beta[A\cap\delta]$ belongs to a countable level of the hierarchy
$L_\alpha[A]$; this proves the equality for $n=1$. Counting the
levels of the hierarchy gives CH, and for $n\geq2$ the submodels
considered above give $2^{\omega_{n-1}}=\omega_n$.

We use the good $\Sigma_1$ wellorder $<_{L[A]}$ throughout. For
$n=1$, the diamond sequence is defined as follows. For a countable
limit ordinal $\delta$, let $(u,c)$ be the least pair in
$L[A\cap\delta]$ such that
\[
 u\subseteq\delta,\quad c\subseteq\delta\text{ is club},\quad
 \forall\xi\in c\ (u\cap\xi\ne d_\xi),
\]
and put $d_\delta=u$; if there is no such pair, or if $\delta$ is not
a limit ordinal, put $d_\delta=\varnothing$. Every pair in
$L[A\cap\delta]$ which satisfies this condition occurs in a countable
level of that hierarchy. Indeed, with
$\lambda=\max(\omega,|\delta|^{L[A\cap\delta]})$, condensation
with $\delta$ fixed puts the pair below
$(\lambda^+)^{L[A\cap\delta]}\leq\omega_1^{\mathcal W}$.
The recursion through $\delta$ is hereditarily countable.
That a pair is the least one satisfying the condition, or that
there is no such pair, is verified by quantifying over the
countable levels of $L[A\cap\xi]$ at each $\xi\leq\delta$. These
levels and the sets $A\cap\xi$ are first-order definable over
$(\HC,\in,A)$, so the recursion is first order over this structure.

This is a diamond sequence. Otherwise, let $(X,C)$ be the
$<_{L[A]}$-least counterexample, and take a countable elementary
submodel of a sufficiently large limit
level $L_\theta[A]$ which contains $A$, $\omega_1$,
$\langle d_\xi:\xi<\omega_1\rangle$, $X$ and $C$. Let $\delta$ be its
intersection with $\omega_1$. The transitive collapse maps $A$ to
$A\cap\delta$, the diamond sequence to
$\langle d_\xi:\xi<\delta\rangle$, and $(X,C)$ to
$(X\cap\delta,C\cap\delta)$. The image of $(X,C)$ is the least pair
in $L[A\cap\delta]$ which satisfies the condition displayed above: all
its predecessors in the canonical wellorder belong to the collapsed
level, and the condition is absolute for this level. Hence
$d_\delta=X\cap\delta$, contradicting $\delta\in C$. The stationary
sets $S^1_t$ and $S^1_*$ are defined from this sequence as in
Section~\ref{sec:source}.

For $n\geq2$ use the recursion from
Lemma~\ref{lem:canonical-diamond} with $<_{L[A]}$. In the
condensation argument the elementary submodels contain $A$ as an
element and $\omega_{n-1}$ as a subset, so their collapses fix $A$.
The witnesses for $\Phi_D$ and for the analogous definition
$\Phi_B$ belong to $H_{\omega_n}^{\mathcal W}=L_{\omega_n}[A]$.

For the sets $P^n_t(C)$, choose the club
$E\subseteq\operatorname{Lim}(C)$ in $L[A,C]$ and the codes $w_d$ of
wellorders using $<_{L[A]}$. As before, the odd part of $P\cap\beta$
gives $\beta\in E$, and the proof of Lemma~\ref{lem:direct-padding}
applies with $L_\beta[A]$ in $N$. The collapses used in the converse fix $A$
and compute $\Phi_D$ and $\Phi_B$ correctly, so
Lemma~\ref{lem:club-from-segments} applies.
\end{proof}

\begin{lemma}\label{lem:PD-models}
The formulas of Section~\ref{sec:decoders} can be modified so that
the models which they describe contain $A$ as an element, with the
same bounds on their complexity: there are formulas
\begin{gather*}
 V_1^A\text{ (first order)},\\
 O_r^A,\ C_r^{A,+},\ C_r^{A,-}\in\PiT_{r-2}[A],\qquad
 V_r^A\in\PiT_{r-1}[A]\qquad(r\geq2)
\end{gather*}
with the properties stated in Lemmas~\ref{lem:direct-V1-explicit}
and~\ref{lem:recursive-V-explicit}, and such that in the models
described by $V^A_r$ an additional constant $a$ is interpreted by $A$.
\end{lemma}
\begin{proof}
Add a constant $a$ to the language of the theories $\mathsf T_r$,
replace $L_\beta$ by $L_\beta[a]$, and replace the statements about
cardinals of $L$ by the corresponding statements about $L[a]$. The
surjections $f_{j,\alpha}$ and the other definable functions are
defined using $<_{L[a]}$.

Let $r=1$, and let $Y$ be a code with equality relation $=_Y$ and
membership relation $E_Y$. We can require the constants $c_\xi$ to
denote the ordinals $\xi$ by the first-order conditions
\[
 E_Y(u,c_\xi)\iff\exists\eta<\xi\ (u=_Yc_\eta)
                         \qquad(\xi<\omega_1).
\]
Indeed, for the transitive collapse $N$, induction gives
$c^N_\xi=\xi$. If all these constants belong to $b_1$ and
$N\models b_1=\omega_1$, then $b_1^N$ is the true $\omega_1$:
a larger value would give a surjection from $\omega$ onto the
true $\omega_1$ in $N$. Now require
\[
 N\models a\subseteq b_1,\qquad
 E_Y(c_\xi,a)\iff A(\xi)\quad(\xi<\omega_1).       \tag{P1}
\]
These are first-order conditions over $(\HC,\in,A)$, and they imply
$a^N=A$. This gives the formula $V_1^A$ and the analogue of
Lemma~\ref{lem:direct-V1-explicit}.

For $r\geq2$ we follow the proof of
Lemma~\ref{lem:recursive-V-explicit}. By induction, the models
described by the formulas of the previous step contain $A$, and
hence they compute the functions defined from $<_{L[A]}$ correctly.
In the definition of $V_r^A$, the condition (P1) is expressed in the
same way as the conditions which ensure $z_i^N=Z_r(X_i)$, by means
of codes of ordinals and witnesses. Here $\xi$
ranges over countable ordinals and $A(\xi)$ is an atomic formula, so
no additional class quantifier is needed, and the complexity of the
formulas is computed as before. For the completeness assertions, the
elementary submodels are chosen so as to contain $A$ as an element
and $\omega_r$ as a subset; their transitive collapses map $A$ to
itself and satisfy (P1).
\end{proof}

\begin{theorem}\label{thm:PD-relative}
There is a cardinal-preserving forcing $\mathbb P_A\in \mathcal W$ adding no reals
and forcing $\PiT_n[A]$-$\UP$ for every finite $n\geq1$.
\end{theorem}
\begin{proof}
Let $I$ be the set of codes $\ell(e,x,y)$ with
$x,y\in\R^{\mathcal W}$, and carry out the constructions of
Sections~\ref{sec:source}--\ref{sec:names}
in $\mathcal W$, using the stationary sets and almost disjoint
families of Lemma~\ref{lem:PD-source} and the wellorder
$<_{L[A]}$ wherever $<_L$ was used. The definitions of the conditions,
the supports, the memories, the term forcings and the subforcings are
unchanged. The proofs use only GCH, the closure and linkedness of the
iterands, and the stationary sets $S^n_*$, all of which are available
in $\mathcal W$. Thus Theorem~\ref{thm:master-full} and the results on
which it rests hold over this ground model.

Write
\[
 F^A_{\vec T}=C^A_{\vec T}*\dot D^A_{\vec T},\qquad
 B^A_{\vec T}=\mathcal W[G_{C^A_{\vec T}}].
\]
The proof of Proposition~\ref{prop:compile} goes through with $A$ as
an additional predicate. Fact~\ref{fact:closedabs} allows ground-model
predicates, and Lemma~\ref{lem:forcing-prefix-verified} uses the
canonical name for $A$. Hence each $\PiT_n[A]$ formula
$\varphi$ has a $\PiT_1(H_{\omega_n},A)$ translation $\Phi^A_\varphi$,
independent of $\vec T$, such that
\[
 B^A_{\vec T}[G_{D^A_{\vec T}}]\models\varphi(x,y)
 \iff B^A_{\vec T}\models\Phi^A_\varphi(x,y).
\]
Lemmas~\ref{lem:PD-source} and~\ref{lem:PD-models}, together with
the proof of Proposition~\ref{prop:exact}, give formulas
\[
 \Code_n^A\in\Sigma^2_n[A],\qquad
 F^A_{\vec T}\Vdash(\Code_n^A(t)\iff t\in T_n).       \tag{P2}
\]

The forcing $C^A_{\vec T}$ is cone homogeneous, and
Fact~\ref{fact:cone} applies to statements whose parameters are $A$
and reals of $\mathcal W$. The proof of
Proposition~\ref{prop:determination} therefore gives (A1), and (P2)
gives (A2). For (A3), let $\vec R\leq\vec T$ differ only at the
coordinate $n$. By Lemma~\ref{lem:comparison-domain} over $L[A]$ with
$k=n$, the models $B^A_{\vec R}$ and $B^A_{\vec T}$ have the same
$H_{\omega_n}$ and the same predicate $A$. Hence a
$\PiT_1(H_{\omega_n},A)$ statement which holds in the larger model
holds in the smaller one, and we can argue as in
Proposition~\ref{prop:ownaxis}. For (A4), use
Lemma~\ref{lem:comparison-domain} over $L[A]$ with $k=n+1$, as in
Proposition~\ref{prop:higheraxis}. The forcing relation is
definable in $\mathcal W$ as in Remark~\ref{rem:uniform}.

The proof of Theorem~\ref{thm:abstract} can now be carried out in
$\mathcal W$, with the set $I$ defined above and with
$<_{L[A]}$-least choices in place of $<_L$-least choices. There are
$\omega_1$ pairs $(e,x)$, the recursion for each $n$ stops before
$\omega_2$, and Replacement gives the
sequence $\vec T_*\in \mathcal W$. For
$\mathbb P_A=F^A_{\vec T_*}$, the uniformizing functions are given by
\[
 U_{n,e}(x,y)\iff
 \varphi_{n,e}(x,y)\land
                  \neg\Code_n^A(\ell(e,x,y)),
\]
where $\langle\varphi_{n,e}:e<\omega\rangle$ enumerates the
$\PiT_n[A]$ formulas. The verification in Section~\ref{sec:selection}
uses only (A1)--(A4), so it gives the conclusion, including the
preservation of cardinals and the fact that no reals are added.
\end{proof}

\subsection{Eliminating the predicate}

\begin{lemma}\label{lem:PD-eliminate}
Suppose that membership and nonmembership in $A\subseteq\omega_1$
have parameter-free $\Sigma^2_1$ definitions.
Then every $\Sigma^2_n[A]$ or $\PiT_n[A]$ formula, $n\geq1$, is
equivalent to a $\Sigma^2_n$ or $\PiT_n$ formula, respectively. This
holds in every model in which the two definitions are correct, and
also for formulas with additional free class variables.
\end{lemma}
\begin{proof}
Write the definitions as $\exists X\,\theta_+(\xi,X)$ and
$\exists X\,\theta_-(\xi,X)$, where $\theta_+$ and $\theta_-$
are first order.
By the axiom of choice, witnesses for all $\xi<\omega_1$ can be
collected into one class. Thus there is a first-order formula
$D(B,Z)$ asserting $B\subseteq\Ord\cap\HC$ and
\[
 \forall\xi<\omega_1\bigl[
  (\xi\in B\to\theta_+(\xi,Z_\xi))\land
  (\xi\notin B\to\theta_-(\xi,Z_\xi))\bigr],
\]
such that $\exists Z\,D(B,Z)$ holds exactly for $B=A$.
For a formula $\varphi(A)$ which begins with an existential class
quantifier use
\[
 \exists B\exists Z\,[D(B,Z)\land\varphi(B)];
\]
for a formula which begins with a universal class quantifier use
\[
 \forall B\forall Z\,[D(B,Z)\to\varphi(B)].
\]
The additional class quantifiers can be merged with the first block
of class quantifiers of $\varphi$.
\end{proof}

\subsection{The ground model \texorpdfstring{$L(\mathbb R)^{M_\omega}$}{L(R) of M-omega}}

Assume the large-cardinal hypothesis of Theorem~\ref{thm:PD}.
Let $\mathcal M=M_\omega$, $\mathcal R=\mathbb R^{\mathcal M}$, and $\mathcal W=L(\mathcal R)^{\mathcal M}$.
At first we use only $\mathcal W\models\mathsf{ZF}$. In the following
definition of $A$ we identify reals with subsets of $\omega$, using
the fixed recursive coding when a relation is given on
$\omega^\omega$.
In first-order definitions over $\mathcal W$, the set $\mathcal R$ is
allowed as a parameter,
as in \cite[Section~7, printed p.~66]{SteelOutline}.

\begin{lemma}\label{lem:PD-ground}
There is $A\subseteq\kappa=\omega_1^{\mathcal W}$ such that
\[
 \mathcal W=L[A]\models\mathsf{ZFC}+\mathsf{CH}+\mathsf{PD}.
\]
Both membership and nonmembership in $A$, for countable ordinals, are
$\Sigma_1$ over $\mathcal W$ in the parameter $\mathcal R$, without
further parameters.
\end{lemma}
\begin{proof}
Let $\prec$ be the wellorder of $\mathcal R$ given by the order of
construction in $\mathcal M$.
Steel states that this order is $\Delta_1$ over $\mathcal W$ in the
paragraph after Corollary~7.21
\cite[printed p.~78]{SteelOutline}. Its graph is a set in $\mathcal W$
by Separation. Acceptability gives
\[
 \mathcal P(\omega)^{\mathcal M}\subseteq J^{\mathcal M}_{\omega_1^{\mathcal M}};
\]
see \cite[Definition~2.3 and the following paragraph]{SteelOutline}.
The hierarchy is continuous, and each $J^{\mathcal M}_\beta$ for
$\beta<\omega_1^{\mathcal M}$ is countable in $\mathcal M$, since it
is obtained by closing under countably many finitary rudimentary
functions. Thus each
proper initial segment of $\prec$ is countable in
$\mathcal M$. A bijection witnessing this is coded by a real and hence
belongs to $\mathcal W$. Since $\mathcal M$ and $\mathcal W$ have the
same reals, they have the same $\HC$ and the same $\omega_1$.

Every proper initial segment of $\prec$ is therefore countable
in $\mathcal W$. In $\mathsf{ZF}$ it follows that the order type of
$\prec$ is at most $\omega_1^{\mathcal W}$, and by Cantor's theorem it
is uncountable. Thus the unique
increasing enumeration
\[
 e:\kappa\longrightarrow \mathcal R
\]
belongs to $\mathcal W$. Neither the axiom of choice nor CH in
$\mathcal W$ has been used so far.
Put
\[
 A=\{\omega\cdot\xi+n:\xi<\kappa,\ n\in e(\xi)\}.       \tag{P3}
\]
Then $A\in \mathcal W$ and $L[A]\subseteq \mathcal W$. Conversely
every real $e(\xi)$ can be computed from $A$, so $e$ and its range
$\mathcal R$ belong to $L[A]$.
By the minimality of $L(\mathcal R)$ we obtain $\mathcal W=L[A]$.
Hence the axiom of choice holds in $\mathcal W$,
and $e$ witnesses CH. Every level of projective determinacy
holds in $\mathcal M$ by the projective determinacy
theorem~\cite{MartinSteelPD},
using finitely many Woodin cardinals and a measurable cardinal below
a later Woodin cardinal. Since $\mathcal W$ has the same reals as
$\mathcal M$, it has the same projective sets and the same strategies
for games on $\omega$, so $\mathcal W\models\mathsf{PD}$.

For the computation of the complexity, let $E(\xi,r)$ assert
\[
 r\in \mathcal R\ \land\ \exists f\,
 [f:(\xi,\in)\cong(\{s\in \mathcal R:s\prec r\},\prec)].
\]
The statement that $f$ is such an isomorphism uses only
quantifiers bounded by $\mathcal R,\xi,f$ and the $\Delta_1$
definition of $\prec$. By Collection, the class of $\Sigma_1$ formulas
is closed under bounded quantification; hence
$E$ is $\Sigma_1$ over $\mathcal W$ in the parameter $\mathcal R$.
For $\xi<\kappa$ it holds exactly if $r=e(\xi)$. For
$\alpha=\omega\cdot\xi+n<\kappa$,
\[
 \begin{aligned}
 \alpha\in A&\iff\exists r\in \mathcal R\,[E(\xi,r)\land n\in r],\\
 \alpha\notin A&\iff\exists r\in \mathcal R\,[E(\xi,r)\land n\notin r].
 \end{aligned}                                                   \tag{P4}
\]
Division of ordinals by $\omega$ with remainder is
$\Delta_1$-definable in $\mathsf{ZF}$.
Both formulas in (P4) are therefore $\Sigma_1$, and they are
complementary, because for every $\xi<\kappa$ there is exactly one
$r$ with $E(\xi,r)$. If we substitute the
definition of $\prec$, then neither $e$, $\kappa$ nor the
graph of $\prec$ occurs as a parameter.
\end{proof}

\subsection{Definability of \texorpdfstring{$A$}{A} in the forcing extensions}

\begin{lemma}\label{lem:PD-translation}
Suppose $\mathcal V\models\mathsf{ZFC}+\mathsf{CH}$, put
$\mathcal R_{\mathcal V}=\mathbb R^{\mathcal V}$, and assume
\[
 \mathcal L=L(\mathbb R)^{\mathcal V}\models\mathsf{ZFC}+\mathsf{CH},\qquad
 \HC^{\mathcal L}=\HC^{\mathcal V}.
\]
For every $\Sigma_1$ formula $\varphi(\mathcal R_{\mathcal V},a)$ there
is a $\Sigma^2_1$ formula $\widehat\varphi(a)$, obtained recursively
from $\varphi$, such that
\[
 \mathcal L\models\varphi(\mathbb R^{\mathcal V},a)
 \iff \mathcal V\models\widehat\varphi(a)\qquad(a\in\HC^{\mathcal V}).
\]
The formula $\widehat\varphi$ has no parameter for
$\mathbb R^{\mathcal V}$.
\end{lemma}
\begin{proof}
Write $\varphi$ as $\exists z\,\psi(z,\mathcal R_{\mathcal V},a)$ with bounded $\psi$,
and put $H=\HC^{\mathcal V}$ and $\tau=\omega_1^{\mathcal V}=\omega_1^{\mathcal L}$.
Suppose that $\varphi$ holds in $\mathcal L$. Choose a sufficiently
large level
$L_\theta(\mathcal R_{\mathcal V})$ which contains a witness, $H$ and
all levels of the hierarchy up to a level containing these sets. In
$\mathcal L$, take an elementary submodel of
size $\tau$ which contains $\mathcal R_{\mathcal V},H,a$ as elements
and every member of $H$. Its transitive collapse is the identity on
$\operatorname{TC}(\{\mathcal R_{\mathcal V}\})$.
Condensation for the hierarchy $L_\alpha(\mathcal R_{\mathcal V})$
gives a level $L_\beta(\mathcal R_{\mathcal V})$ of size
at most $\tau$, with $\beta<\tau^+$, which contains
$\mathcal R_{\mathcal V},H,a$ and a witness for $\psi$.

This instance of condensation is proved by collapsing the
hierarchy together with the satisfaction relations of its levels.
At a successor step, elementarity shows that every element of the next
level has a definition over the preceding level, with a formula and
parameters in the submodel, and that every set defined in this way
belongs to the submodel. Hence the collapse of a successor level is
exactly the set of all subsets of the collapsed preceding level which
are definable over it. At limit steps the collapsed levels are
unions. By induction, the collapsed levels are the levels of the
hierarchy $L_\alpha(\mathcal R_{\mathcal V})$.

The level $L_\beta(\mathcal R_{\mathcal V})$, the sequence of the
preceding levels and their satisfaction relations
have hereditary size at most $\tau$. Hence they can be coded,
together with a witness for $\psi$, by a single class
$X\subseteq H$. We let $\widehat\varphi(a)$ assert that there is a
class $X$ with the following properties. $X$ codes a
well-founded extensional relation $E_X$ on a set $D\subseteq H$,
with a distinguished copy of $H$ and a distinguished element
representing $\mathcal R_{\mathcal V}$, whose
members are exactly the copies of the reals. The copy
of $H$ consists of the elements $\langle0,u\rangle$ for $u\in H$,
and the $E_X$-predecessors of $\langle0,u\rangle$ are exactly the
elements $\langle0,v\rangle$ with $v\in u$. The coded hierarchy
starts with $\operatorname{TC}(\{\mathcal R_{\mathcal V}\})$, its
successor levels consist exactly of the definable subsets of the
preceding levels, its limit levels are unions, and the coded
satisfaction relations satisfy the Tarski conditions, where formulas
are represented by their G\"odel numbers.
Finally, the last level contains $a$ and a
witness for $\psi(z,\mathcal R_{\mathcal V},a)$.

All these conditions are first order over $H$, because elements,
indices of levels, parameters and G\"odel numbers are represented by
elements of $H$.
Well-foundedness is first order as well, since it is expressed by
\[
 \neg\exists f\in H\,[f:\omega\to D\ \land\
          \forall n<\omega\ f(n+1)\mathrel{E_X}f(n)],
\]
and every countable sequence of members of $H$ belongs to $H$.
All components are coded by the one class $X$, so $\widehat\varphi$
is $\Sigma^2_1$. We have shown that $\widehat\varphi(a)$ holds in
$\mathcal V$ if $\varphi$ holds in $\mathcal L$.

Conversely, suppose that $X$ is as in $\widehat\varphi(a)$. The
transitive collapse of the coded structure contains the true
$\operatorname{TC}(\{\mathcal R_{\mathcal V}\})$ as its first level.
By the conditions on successor and limit levels, induction shows that
its levels are the levels $L_\gamma(\mathcal R_{\mathcal V})$.
The witness therefore belongs to $\mathcal L$, and since bounded
formulas are absolute,
$\mathcal L\models\psi(z,\mathcal R_{\mathcal V},a)$.
\end{proof}

\begin{lemma}\label{lem:PD-recovery}
For the set $A$ in (P3), there are parameter-free $\Sigma^2_1$ formulas
$\sigma_+(\alpha)$ and $\sigma_-(\alpha)$ such that, in every
set-forcing extension $\mathcal W'$ of $\mathcal W$ preserving $\HC$,
\[
 \begin{aligned}
 \mathcal W'\models\sigma_+(\alpha)&\iff\alpha\in A,\\
 \mathcal W'\models\sigma_-(\alpha)&\iff\alpha\notin A
 \end{aligned}
 \qquad(\alpha<\omega_1^{\mathcal W}).
\]
\end{lemma}
\begin{proof}
Let $\sigma_+$ and $\sigma_-$ be the translations of the two formulas
in (P4) given by Lemma~\ref{lem:PD-translation}.
Every such $\mathcal W'$ has the same reals and ordinals as $\mathcal W$, so
\[
 L(\mathbb R)^{\mathcal W'}=L(\mathbb R)^{\mathcal W}=\mathcal W.
\]
The enumeration $e\in\mathcal W$ still witnesses CH in $\mathcal W'$.
Hence Lemma~\ref{lem:PD-translation} applies in $\mathcal W'$, with
$\mathcal L=\mathcal W$.
The two formulas in (P4) refer to the same order
and the same enumeration, and hence define $A$ and its complement.
Note that by Lemma~\ref{lem:PD-translation} a new class in
$\mathcal W'$ cannot be a witness for $\sigma_+(\alpha)$ or
$\sigma_-(\alpha)$ unless the corresponding formula in (P4) holds in
$\mathcal W$.
\end{proof}

\begin{proof}[Proof of Theorem~\ref{thm:PD}]
Take $\mathcal W=L[A]$ from Lemma~\ref{lem:PD-ground} and the forcing
$\mathbb P_A$ of Theorem~\ref{thm:PD-relative}. In the extension by
$\mathbb P_A$, every $\PiT_n[A]$ relation has a $\PiT_n[A]$
uniformization. By
Lemma~\ref{lem:PD-recovery}, the hypothesis of
Lemma~\ref{lem:PD-eliminate} holds in this extension, so $A$ can be
eliminated without increasing the complexity. Thus the extension
satisfies $\PiT_n$-$\UP$ for every finite $n\geq1$.

The forcing preserves all cardinals of $\mathcal W$ and adds no reals.
It preserves CH, and it changes neither the projective sets nor the
set of strategies for games on $\omega$;
hence it preserves $\mathsf{PD}$. This proves the theorem.
\end{proof}

\subsection{The two uniformization patterns}

\begin{proposition}\label{prop:Sigma-ground}
Both $L$ and $\mathcal W=L(\mathbb R)^{M_\omega}$ satisfy
$\Sigma^2_n$-$\UP$ for every finite $n\geq1$.
\end{proposition}
\begin{proof}
First work in $\mathcal W=L[A]$ from Lemma~\ref{lem:PD-ground}, with
$A$ as an additional predicate, and put $H=\HC^{\mathcal W}$ and
$\kappa=\omega_1^{\mathcal W}$.
Every $X\subseteq H$ belongs to
\[
 H_{\kappa^+}^{\mathcal W}=L_{\kappa^+}[A].
\]
The initial segments of the canonical wellorder
$\prec\,=\,<_{L[A]}$ on
$\mathcal P(H)^{\mathcal W}$ have size at most $\kappa$.
This wellorder is $\Sigma^2_1[A]$-good in the following sense. There
is a first-order formula $G_A(T,J,C)$ over $(H,\in,A)$ such that
\[
 \exists C\,G_A(T,J,C)
 \quad\Longleftrightarrow\quad
 \{J_i:i\in D_J\}=\{S\subseteq H:S\prec T\},             \tag{P5}
\]
where the class $J$ codes the family of classes $J_i$, for $i$ in an
index set $D_J\subseteq H$. Moreover, for every $T\subseteq H$ there
are such $J$ and $C$. The enumeration $\langle J_i:i\in D_J\rangle$
of the predecessors of $T$ is arbitrary.

To see this, choose a limit ordinal $\beta<\kappa^+$ such that
$L_\beta[A]$ contains $T,H,A,J$ and the set of all
$<_{L[A]}$-predecessors of $T$; such an ordinal exists because
$<_{L[A]}$ is a good $\Sigma_1(A)$ wellorder.
Let $C$ code, by classes over $H$, the level $L_\beta[A]$, the
sequence of the preceding levels with their satisfaction relations,
and an enumeration of the predecessors of $T$ which
are subsets of $H$. The formula $G_A$ requires that the coded
structure contains copies of the true $H$ and $A$, and that the
elements representing $T$ and the classes $J_i$ have exactly the
copies of the elements of $T$ and of $J_i$ as members. It requires
that every successor level consists exactly of the subsets of the
preceding level which are definable over it with the predicate $A$.
Finally, it requires a witness for the $\Sigma_1$ definition of
$\operatorname{Init}_A(t,I)$, the analogue for $L[A]$ of
$\mathrm{Init}$ in Section~\ref{sec:source-spec}, and that the
enumerated classes are exactly the elements of $I\cap\mathcal P(H)$.
These are first-order conditions, as in
Lemma~\ref{lem:PD-translation}. The transitive collapse of a
well-founded structure satisfying them is a level of the hierarchy
$L_\alpha[A]$, and by the upward absoluteness of $\Sigma_1(A)$
formulas the set $I$ is the set of all $<_{L[A]}$-predecessors of
$T$.
This proves (P5), and the existence of $J$ and $C$ for every $T$ and
for any given enumeration $J$ of the
predecessors of $T$. Lemma~\ref{lem:PD-eliminate} applies to formulas
with the free class variables $T,J$, so the relation in (P5) is also
parameter-free $\Sigma^2_1$ in $\mathcal W$. Equivalently, using the
canonical bijection between $H$ and
$\kappa$, we obtain a good wellorder of
$\mathcal P(\omega_1)^{\mathcal W}$.

Let $\mathcal E(x,y)$ be a $\Sigma^2_n[A]$ relation. Combine the
class quantified by the leading
existential quantifier with $y$ into one class $T$, using a fixed
recursive coding. Write
\[
 \mathcal E(x,y)\iff\exists T\,
       [\operatorname{out}(T)=y\land \psi_*(x,T)],
\]
where $\psi_*$ includes the statement that $T$ codes a pair of this
form, is first order if
$n=1$, and belongs to $\PiT_{n-1}[A]$ if $n\geq2$.
We uniformize $\mathcal E$ by choosing the $\prec$-least $T$
satisfying $\psi_*(x,T)$:
\[
 \begin{split}
 U(x,y)\iff\exists T,J,C\,[&G_A(T,J,C)
        \land\operatorname{out}(T)=y\land \psi_*(x,T)\\
        &\land\forall i\in D_J\ \neg \psi_*(x,J_i)].
 \end{split}                                                    \tag{P6}
\]
For $n=1$, the part of this formula after the existential class
quantifiers is first order. For $n\geq2$, the statement
$\forall i\in D_J\ \neg\psi_*(x,J_i)$
is $\Sigma^2_{n-1}[A]$: the quantifier over $i$ is first order,
and by the axiom of choice the witnesses for all $i\in D_J$ can be
collected into one class.
Its conjunction with $\psi_*\in\PiT_{n-1}[A]$, preceded by the
existential class quantifiers, is $\Sigma^2_n[A]$.
Thus (P6) defines a $\Sigma^2_n[A]$ uniformization of $\mathcal E$.

By Lemma~\ref{lem:PD-recovery} and
Lemma~\ref{lem:PD-eliminate}, $A$ can be eliminated without
increasing $n$.
This proves the assertion for $\mathcal W$.
For $L$ use the same argument without the predicate $A$; in this case
nothing has to be eliminated.
\end{proof}

\begin{lemma}\label{lem:dual-UP}
For each finite $n\geq1$, boldface $\Sigma^2_n$ uniformization
and boldface $\PiT_n$ uniformization cannot both hold in a model
of $\mathsf{ZFC}$.
\end{lemma}
\begin{proof}
There is a universal $\PiT_n$ relation $\mathcal U(e,x,y)$, where the
real $e$ codes a formula and its real parameter. To handle arbitrary
first-order formulas, a class coding their satisfaction relation is
added to the last block of class quantifiers: it is characterized by
the first-order Tarski conditions, which are added as a conjunct if
the last block is existential and as the hypothesis of an implication
if it is universal. Since the satisfaction relation exists and is
unique, this does not increase the number of alternations.

Assume $\boldsymbol\Pi^2_n$ uniformization and let
$\mathcal F(e,x,y)$ be a $\boldsymbol\Pi^2_n$ uniformization of
$\mathcal U$. If $\mathcal U_e$ is the graph of a function, then
$\mathcal F_e=\mathcal U_e$. Now
\[
 \mathcal D(e,y)\iff\neg \mathcal F(e,e,y)
\]
is $\boldsymbol\Sigma^2_n$ and every section is nonempty, since
$\mathcal F(e,e,\cdot)$ has at most one value. If $\boldsymbol\Sigma^2_n$
uniformization also held, $\mathcal D$ would have a uniformization
which is a total function with graph
$\mathcal K\in\boldsymbol\Sigma^2_n$. Since this function is total,
\[
 \neg \mathcal K(e,y)\iff\exists z\,[\mathcal K(e,z)\land z\ne y],
\]
so $\mathcal K$ is also $\boldsymbol\Pi^2_n$. Choose $e_0$ with
$\mathcal U_{e_0}=\mathcal K$. Then $\mathcal F_{e_0}=\mathcal K$,
whereas the value
$y$ of this function at $e_0$ satisfies both $\mathcal K(e_0,y)$ and
$\neg \mathcal F(e_0,e_0,y)$, a contradiction.
\end{proof}

\begin{proof}[Proof of Corollary~\ref{cor:patterns}]
Proposition~\ref{prop:Sigma-ground} gives $\Sigma^2_n$-$\UP$ in the
two ground models, and
Theorems~\ref{thm:main} and~\ref{thm:PD} give $\PiT_n$-$\UP$ in
their extensions. Lightface uniformization gives boldface
uniformization by pairing the real parameter with the input and
restricting the resulting uniformizing function to that parameter.
Lemma~\ref{lem:dual-UP} gives the failures stated in the corollary.
The forcings add no reals and preserve $\omega_1$, hence
preserve $\HC$ and all projective truth.
\end{proof}

\paragraph{Use of generative AI.}
I used ChatGPT Astra 6 as a co-worker to help develop the forcing construction
and to draft parts of the proofs. The main ideas of this work, in particular the
strategy of removing codes and try to make an inverse iteration, where we iteratively shrink universes of set theory, are old ideas of mine, which I originally
developed in the plain projective context (where they typically failed to work). The
second-order context leaves enough room between what we code and how
the reals behave for it to work here. I also used Claude Fable 5 to improve
the language, to revise the exposition, to avoid typos. Both were used to spot errors, check proofs and to search for relevant literature. I take
responsibility for the contents of this paper.

\bigskip
\noindent
\small
TU Wien, Institut f\"ur Diskrete Mathematik und Geometrie\\
\texttt{stefan.hoffelner@tuwien.ac.at}

\end{document}